\documentclass[reqno,12pt]{amsart}

\usepackage[colorlinks=true]{hyperref}
\hypersetup{urlcolor=blue, citecolor=red}

\usepackage{amssymb,graphicx,color}
\usepackage[latin1]{inputenc}

\usepackage[capitalize,nameinlink,noabbrev]{cleveref}
\crefformat{equation}{#2(#1)#3}
\crefformat{enumi}{#2(#1)#3}

\theoremstyle{plain} 
\newtheorem{theorem}{Theorem}[section]
\newtheorem{lemma}{Lemma}[section]
\newtheorem{proposition}{Proposition}[section]

\newtheorem{stdhyp}{Standing Hypothesis}[section]
\theoremstyle{definition} 
\newtheorem{remark}{Remark}[section]

\newcommand{\dontshow}[1]{}
\newcommand{\tr}{{\operatorname{tr}}}

\begin{document}
\numberwithin{equation}{section}


\title[Euler method for Random ODEs with semi-martingales]{Strong order-one convergence of the Euler method for random ordinary differential equations driven by semi-martingale noises}

\author[P. E. Kloeden]{Peter E. Kloeden}
\author[R. M. S. Rosa]{Ricardo M. S. Rosa}

\address[Peter E. Kloeden]{Mathematics Department, University of T\"ubingen, Germany}
\address[Ricardo M. S. Rosa]{Instituto de Matem\'atica, Universidade Federal do Rio de Janeiro, Brazil}

\email[P. E. Kloeden]{kloeden@math.uni-frankfurt.de}
\email[R. M. S. Rosa]{rrosa@im.ufrj.br}

\date{\today}

\makeatletter
\@namedef{subjclassname@2020}{\textup{2020} Mathematics Subject Classification}
\makeatother
\subjclass[2020]{60H35 (Primary), 65C20, 60H25, 35R60 (Secondary)}

\keywords{random ordinary differential equations, Euler method, strong convergence, It\^o process, finite-variation process, semi-martingale, fractional Brownian motion}.

\begin{abstract}
    It is well known that the Euler method for a random ordinary differential equation $\mathrm{d}X_t/\mathrm{d}t = f(t, X_t, Y_t)$ driven by a stochastic process {$\{Y_t\}_{t\in I}$, on a time interval $I,$} with $\theta$-H\"older sample paths {is of strong order} $\theta$ with respect to the time step, provided $f=f(t, x, y)$ is sufficiently regular and with suitable bounds. This order is known to increase to $1$ in some special cases. Here, it is proved that, in many more typical cases, further structures on the noise can be exploited so that the strong convergence is of order 1. In fact, we prove so for any semi-martingale noise. This includes It\^o diffusion processes, point-process noises, transport-type processes with sample paths of bounded variation, and time-changed Brownian motion. The result follows from estimating the global error as an iterated integral over both large and small mesh scales, and by switching the order of integration to move the critical regularity to the large scale. The work is complemented with numerical simulations showing the optimality of the strong order 1 convergence in those cases, and with an example with fractional Brownian motion noise with Hurst parameter $0 < H < 1/2,$ which is not a semi-martingale and for which the order of convergence is $H + 1/2$, hence lower than the attained order 1 in the semi-martingale case, but still higher than the order $H$ of convergence expected from previous works.
\end{abstract}

\maketitle

\tableofcontents

\section{Introduction}
\label{secintro}

A \textbf{random ordinary differential equation (RODE or Random ODE)} is a differential equation of the form
\begin{equation}
  \label{rodeeq}
  \begin{cases}
    \displaystyle \frac{\mathrm{d}X_t}{\mathrm{d} t} = f(t, X_t, Y_t), \\
    \left. X_t \right|_{t = 0} = X_0,
  \end{cases}
\end{equation}
for $t$ on a time interval which we take to be $I=[0, T]$, with $T > 0$, where the initial condition $X_0$ is a given random variable, and the noise $\{Y_t\}_{t\in I}$ is a given stochastic process, both defined on a probability space $(\Omega, \mathcal{F}, \mathbb{P}).$ The function $f=f(t, x, y)$ is such that the differential equation in \cref{rodeeq} can be either a scalar or a system of equations and the noise can also be either a scalar or a vector-valued process.

Random ODEs appear in applications where the parameters either vary randomly or are not fully known. As such, they have been playing an important role in modeling practical problems, including uncertainty quantification (see e.g. \cite{Asai2016,BogdanoffGoldbergBernard1961,FreidlinWentzell1992,HanKloeden2017,MerdanBekiryaziciKesemenKhaniyev2017,NeckelRupp2013,SalakoShen2020,StrasserTheisMarr2012}). They are less known than the more popular It\^o-diffusion stochastic differential equations (SDE), but, depending on the situation, they may actually lead to more realistic models.

As in many other models, numerical approximations are of crucial importance. One popular method for approximating Random ODEs is the Euler method. It consists in approximating the solution on a uniform time mesh $t_j = t_j^N = j\Delta t_N$, $j = 0, \ldots, N$, with fixed time step $\Delta t_N = T/N$, for a given $N\in \mathbb{N},$ via
\begin{equation}
  \label{emscheme}
  \begin{cases}
    X_{t_j}^N = X_{t_{j-1}}^N + \Delta t_N f(t_{j-1}, X_{t_{j-1}}^N, Y_{t_{j-1}}), & j = 1, \ldots, N, \\
    X_0^N = X_0.
  \end{cases}
\end{equation}

An approximation $\{X_{t_j}^N\}_{j = 0, 1, \ldots, N}$ is said to converge to $\{X_t\}_{t\in I}$ with \emph{strong order $\theta>0$} when there exists a constant $C \geq 0$ such that
\begin{equation}
    \label{strongordertheta}
    \max_{j=0, \ldots, N}\mathbb{E}\left[ \left\| X_{t_j} - X_{t_j}^N \right\| \right] \leq C \Delta t_N^\theta, \qquad \forall N \in \mathbb{N},
\end{equation}
where $\mathbb{E}[\cdot]$ indicates the expectation of a random variable on the underlying probability space, and $\|\cdot\|$ is a suitable norm in the appropriate phase space.

There are other important notions of convergence, such as pathwise convergence, weak convergence, mean-square convergence, and $p$-th mean convergence (see e.g. \cite{AsaiKloeden2016,GruneKloeden2001,HanKloeden2017,HighamKloeden2021,JentzenKloeden2011,KloedenJentzen2007}), but here we focus only on the strong convergence.

Under certain regularity conditions on $f$ and when the noise $\{Y_t\}_{t\in I}$ has $\theta$-H\"older continuous sample paths, the result in \cite[Theorem 3]{WangCaoHanKloeden2021} proves the order of convergence of the Euler scheme, in the mean square sense, to be of order $\theta$ with respect to the time step. This implies the strong convergence \cref{strongordertheta} with the same order $\theta$. In the case of fractional Brownian motion noise with Hurst parameter $H$, \cite[Theorem 2]{WangCaoHanKloeden2021} proves the mean square convergence to be of order $H$, for any $0 < H \leq 1.$

In some special cases, the order of convergence of the Euler method is known to be higher than the H\"older exponent of the noise. For example, for It\^o diffusion noises, the RODE can be seen as part of a system of SDEs, with the Euler method for the RODE component coinciding with the Milstein method for the SDE, which is known to be of strong order 1 (see \cite[Section 5]{WangCaoHanKloeden2021}, in particular \cite[Section 5.2, Example 12 and Remark 13]{WangCaoHanKloeden2021}).

\begin{figure}[htb]
    \centering{\includegraphics[width=0.9\textwidth]{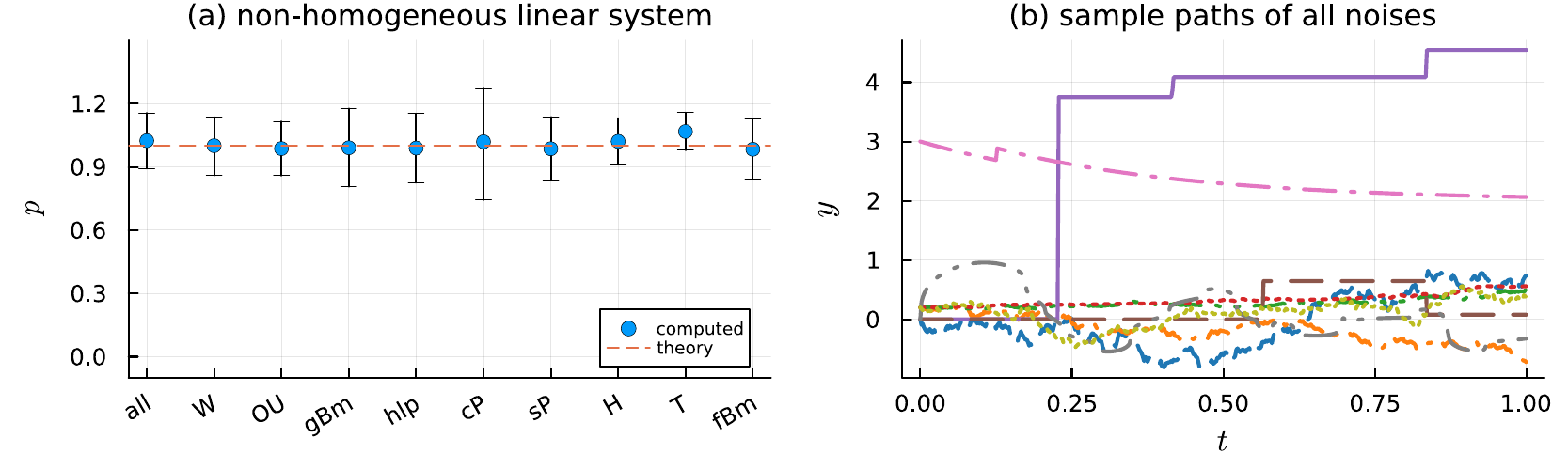}}
    \caption{(a) Estimated order of convergence of the strong error of the Euler method for the linear system \cref{allnoisesRODEsystemintro} as it depends on the noise, starting with all the noises combined (all) and continuing with the scalar equation, each with a different type of noise process, namely Wiener (W), Ornstein-Uhlenbeck (OU), geometric Brownian motion (gBm), a homogeneous linear It\^o diffusion (hlp), compound Poisson process (cP), step point process (sP), Hawkes process (H), transport process (T), and also a fractional Brownian motion with Hurst parameter $1/2 < H < 1.$ (fBm); (b) Sample paths of all the noises used in the linear system \cref{allnoisesRODEsystemintro}.}
    \label{figallnoisesintro}
\end{figure}

In many more classical examples, though, the actual order of convergence seems to be higher than what the theory predicts. Our aim here is to show that, indeed, it is possible to exploit further properties of the noise that yield a strong order 1 convergence for essentially any noise that is a semi-martingale (see \cref{thmsemimartingale} for the precise statement). This includes noises such as point processes, transport processes, It\^o processes, and time-changed Brownian motions, some of which are not even continuous and for which the previous theories did not apply.

This strong order 1 convergence is not only a major improvement from previous results, but it is also remarkable that it follows with a unified proof for all semi-martingale noises. It also has potential extensions to pathwise convergence for RODEs in the case of non-globally Lipschitz functions, at least when driven by finite-variation process. The result should also extend to random partial differential equations. It is not clear, though, whether this technique could also yield improvements in the order of convergence for higher order methods.

\begin{figure}[htb]
    \centering{\includegraphics[width=0.9\textwidth]{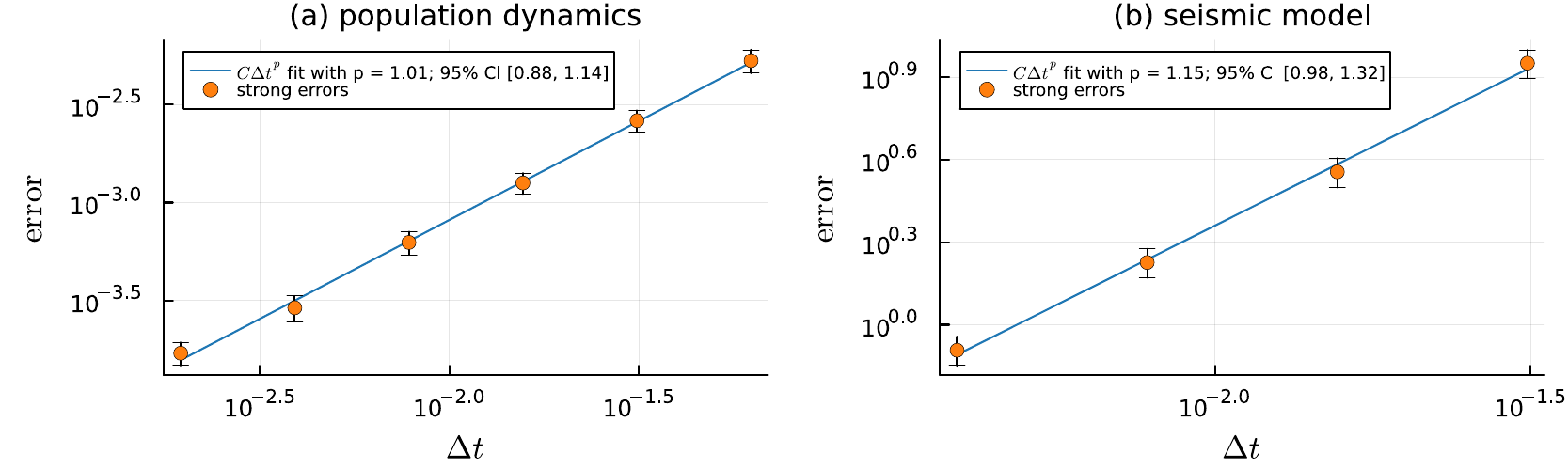}}
    \caption{Estimated order of convergence of the strong error of the Euler method for (a) the population dynamics model \cref{eqpopdyn} and (b) the seismic model \cref{eqearthquake}.}
    \label{figallcombinedpopdynrisk}
\end{figure}

\begin{figure}[htb]
    \centering{\includegraphics[width=0.9\textwidth]{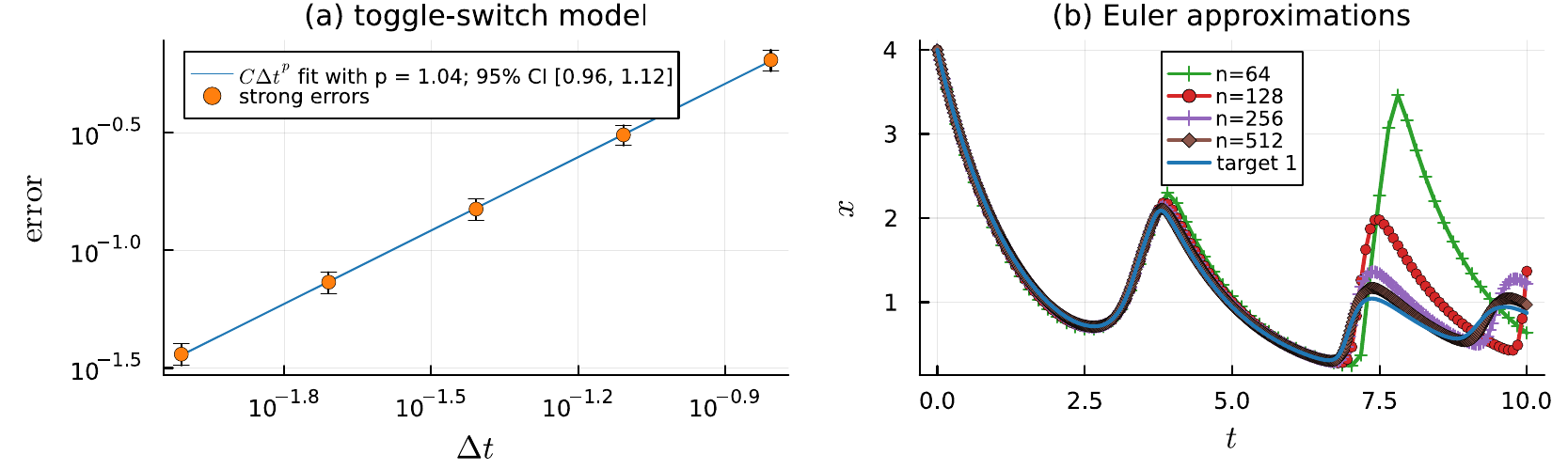}}
    \caption{(a) Estimated order of convergence of the strong error of the Euler method for the toggle-switch model \cref{toggleswitchsystemintro}; (b) Euler approximations of the component $\{X_t\}_{t\in I}$ of a sample path.}
    \label{figtoggleswitchintro}
\end{figure}

We numerically verify, in a series of examples, that the strong order 1 convergence is in fact optimal. Indeed, this work was motived by the fact that the observed convergence in these and many other examples was higher than what the theory predicted. The examples are presented in more details in \cref{secnumex} and in \cref{secsupplementary} and in the github repository \cite{RODEConvEM2023}. Below is a summary of these examples:
\begin{enumerate}
    \item Two types of non-homogeneous linear equations of the form
    \begin{equation}
        \label{allnoisesRODEsystemintro}
        \frac{\mathrm{d}X_t}{\mathrm{d} t} = - \left\|Y_t\right\|^2 X_t + Y_t,
    \end{equation} 
    where $\{Y_t\}_{t\in I}$ is either (a) a univariate process with one of many types of semi-martingale processes or (b) a vector-valued process combining all these noises as independent components.
    \item A logistic model of population dynamics with random coefficients, loosely inspired by \cite[Section 15.2]{HanKloeden2017}, and where the specific growth is the sine of a geometric Brownian motion process $\{G_t\}_{t\in I}$ and with an extra point-process random term $\{H_t\}_{t\in I}$ representing harvest:
    \begin{equation}
        \label{eqpopdyn}
        \frac{\mathrm{d}X_t}{\mathrm{d}t} = \gamma (1 + \varepsilon \sin(G_t)) X_t (1 - \frac{X_t}{r}) - \alpha H_t \frac{X_t}{r + X_t}.
    \end{equation}
    \item A toggle-switch model of gene expression (similar to \cite[Section 7.8]{Asai2016}, originated from \cite{VerdCrombachJaeger2014}, see also \cite{StrasserTheisMarr2012}) driven by a combination of a compound Poisson point process $\{A_t\}_{t\in I}$ and an It\^o diffusion process $\{B_t\}_{t\in I},$
    \begin{equation}
        \label{toggleswitchsystemintro}
        \begin{cases}
            \frac{\displaystyle \mathrm{d}X_t}{\displaystyle \mathrm{d} t} = \left( A_t + \frac{\displaystyle X_t^4}{\displaystyle a^4 + X_t^4}\right)\left(\frac{\displaystyle b^4}{\displaystyle b^4 + Y_t^4}\right) - \mu X_t, \\
            \frac{\displaystyle \mathrm{d}Y_t}{\displaystyle \mathrm{d} t} = \left( B_t + \frac{\displaystyle Y_t^4}{\displaystyle c^4 + Y_t^4}\right)\left(\frac{\displaystyle d^4}{\displaystyle d^4 + X_t^4}\right) - \nu Y_t.
        \end{cases}
    \end{equation}
    \item \label{egearthquase} A mechanical structure seismic model 
    \begin{equation}
        \label{eqearthquake}
        \ddot X_t + 2\zeta_0\omega_0\dot X_t + \omega_0^2 X_t = - \ddot M_t,
    \end{equation}
    driven by a random disturbance in the form of a transport process $M_t = \sum_{i=1}^k \gamma_i (t - \tau_i)_+^2 e^{-\delta_i (t - \tau_i)}\cos(\omega_i (t - \tau_i)),$ simulating seismic ground-motion excitations as inspired by the Bogdanoff-Goldberg-Bernard model in \cite{BogdanoffGoldbergBernard1961}.
    \item An actuarial risk model
    \begin{equation}
        \label{eqrisk}
        \frac{\mathrm{d}X_t}{\mathrm{d}t} = R_t X_t + R_t (C_t + O_t) + \nu O_t + \gamma
    \end{equation}
    for the surplus $U_t = X_t + C_t + O_t$ of an insurance company, inspired by \cite{GerberShiu1998} and \cite{BrigoMercurio2006}, where $C_t = \sum_{i=1}^{N_t} C_i$ defines a compound Poisson process modeling the insurance claims, $\{R_t\}_{t\in I}$ is a geometric Brownian motion process modeling the interest rate, and $\{O_t\}_{t\in I}$ is an Ornstein-Uhlenbeck process, used to transform an additive SDE into a RODE.
    \item A Fisher-KPP partial differential equation
    \begin{equation}
        \label{eqfisherkpp}
        \frac{\partial u}{\displaystyle \partial t} = \mu\frac{\partial^2 u}{\partial x^2} + \lambda u\left(1 - \frac{u}{u_m}\right), \quad (t, x) \in (0, \infty) \times (0, 1),
    \end{equation}
    endowed with random boundary conditions representing migration,
    \begin{equation}
        \label{eqfisherkppbcs}
        \frac{\partial u}{\partial x}(t, 0) = - Y_t, \quad \frac{\partial u}{\partial x}(t, 1) = 0,
    \end{equation}  
    as inspired by the works \cite{SalakoShen2020} and \cite{FreidlinWentzell1992} (see also \cite{Fisher1937} and \cite{KPP1937}), and where the noise $\{Y_t\}_{t\in I}$ is an Ornstein-Uhlenbeck process modulated by a decaying self-exciting Hawkes point process.
    \item \label{eglinearfbm} The linear equation
    \begin{equation}
        \label{eqlinearfbm}
        \frac{\mathrm{d}X_t}{\mathrm{d} t} = -X_t + B^H_t,
    \end{equation}
    where $\{B^H_t\}_{t\in I}$ is a fractional Brownian motion (fBm) with Hurst parameter $0 < H < 1$, which is a Wiener process when $H=1/2$ and which is not a semi-martingale for $H\neq 1/2.$
\end{enumerate}

\begin{figure}[htb]
    \centering{\includegraphics[width=0.9\textwidth]{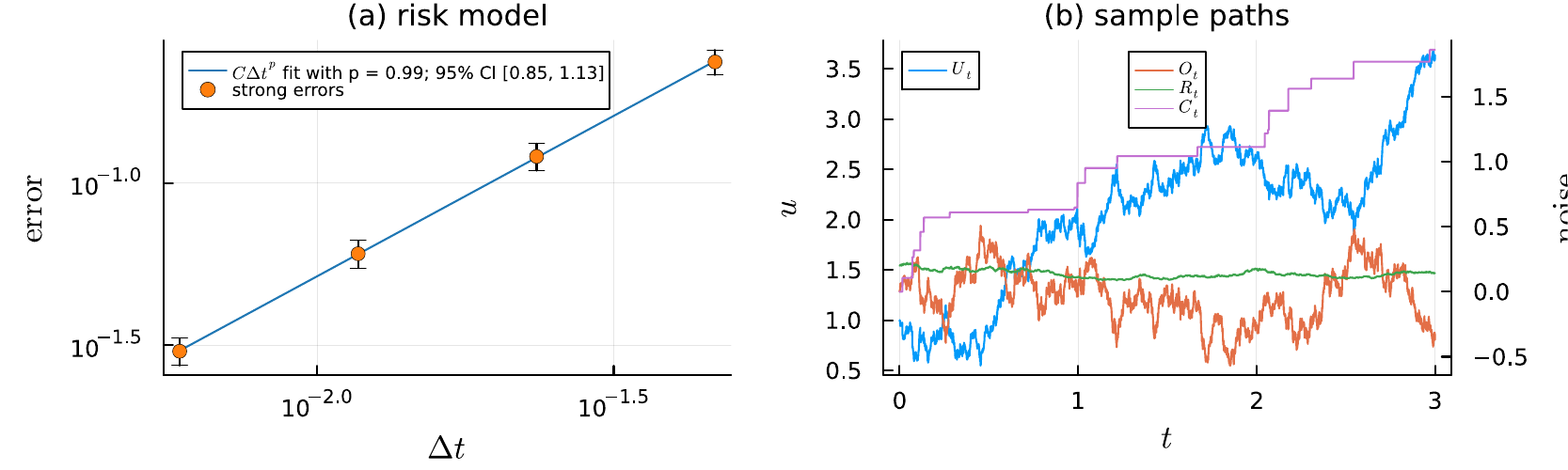}}
    \caption{(a) Estimated order of convergence of the strong error of the Euler method for the insurance risk model \cref{eqrisk}; (b) Sample paths of the surplus $\{U_t\}_{t\in I},$ the Ornstein-Uhlenbeck process $\{O_t\}_{t\in I},$ the interest rate {$\{I_t\}_{t\in I},$} and the insurance claims $\{C_t\}_{t\in I}.$}
    \label{figcombinedrisk}
\end{figure}

\begin{figure}[htb]
    \centering{\includegraphics[width=0.9\textwidth]{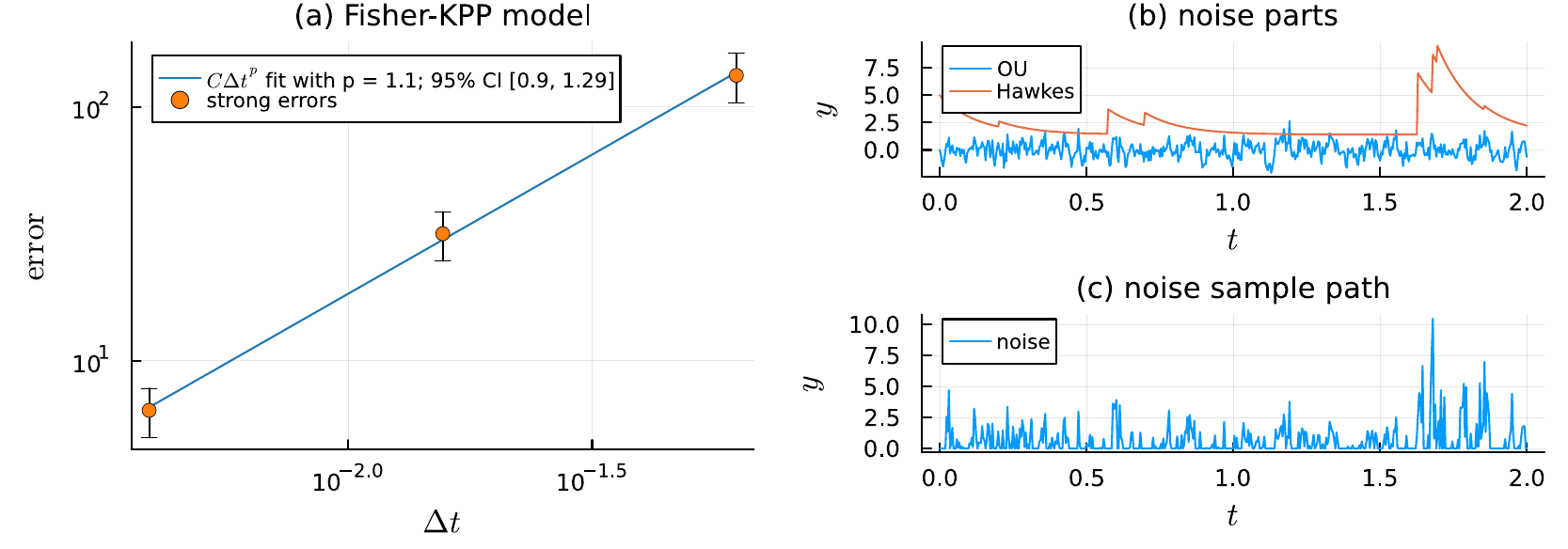}}
    \caption{(a) Estimated order of convergence of the strong error of the Euler method for the Fisher-KPP model \cref{eqfisherkpp}-\cref{eqfisherkppbcs}; Sample paths of (b) the Ornstein-Uhlenbeck and Hawkes processes, and of (c) the noise obtained as the product of these processes.}
    \label{figcombinedfisherkpp}
\end{figure}

In all but the last example, the numerically obtained convergence was found to be of strong order approximately 1, while previous results would either not be applicable or yield a much lower order of convergence (see \cref{figallnoisesintro,figallcombinedpopdynrisk,figtoggleswitchintro,figcombinedrisk,figcombinedfisherkpp}).

\begin{figure}[htb]
    \centering{\includegraphics[width=0.56\textwidth]{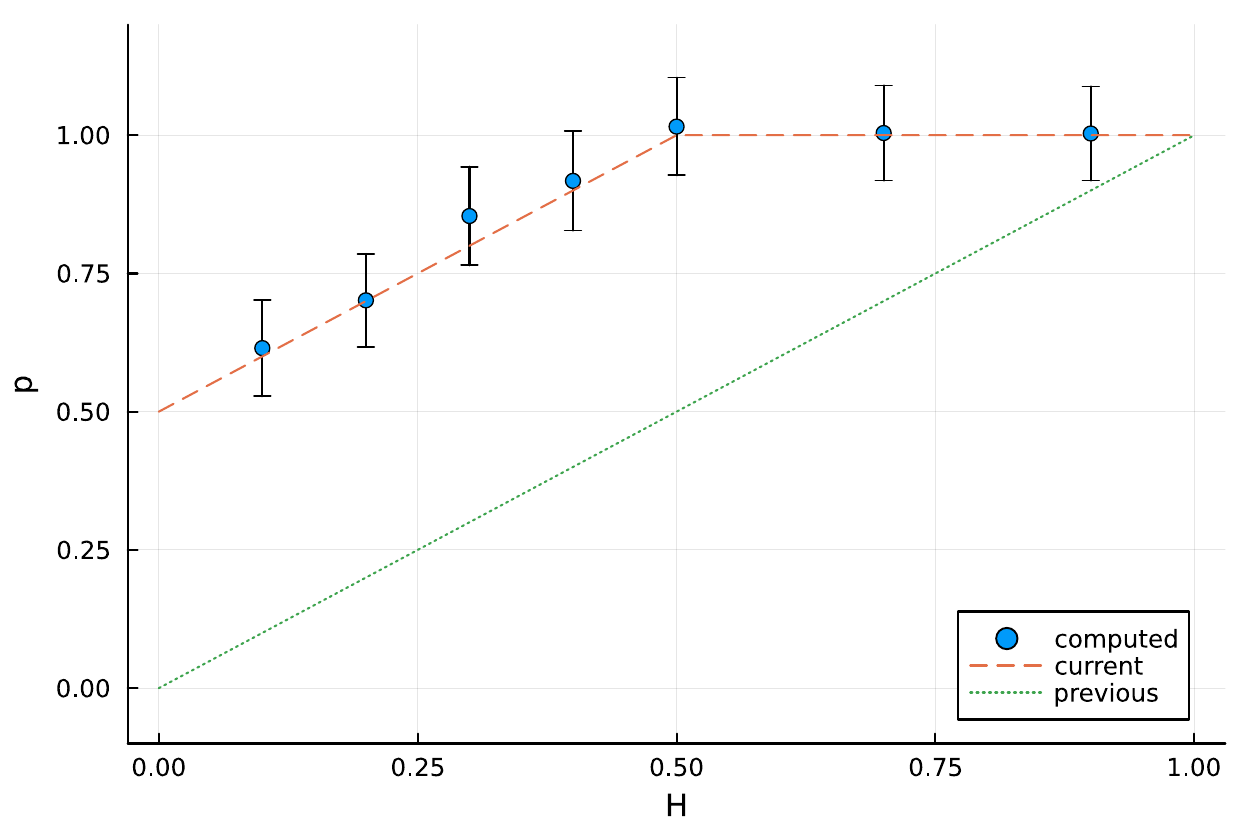}}
    \caption{Estimated order $p$ of strong convergence of the Euler method for \cref{eqlinearfbm} at different values of the Hurst parameter $H$ (scattered plot) along with the theoretical value $p=\min\{H + 1/2, 1\}$ proved here (dashed line), and the previously known theoretical value $p=H$ (dotted).}
    \label{figorderdepHfBmintro}
\end{figure}

In the last example, with the fractional Brownian motion, the order was found to be $\min\{H+1/2, 1\},$ as illustraded in \cref{figorderdepHfBmintro}, well above the previously expected order of convergence $H.$

Our aim is not only to report on these discrepancies but also to prove and explain why we have such better convergence rates. 

Details about the numerical setup for estimating the order of strong convergence are given in \cref{secnumex}. Details concerning the first example \cref{allnoisesRODEsystemintro} are given in \cref{subseclinearrodenumerics}. A more detailed description of the numerical setup and the other examples can be found in the open-access github repository \cite{RODEConvEM2023}.

Concerning the proof of strong order 1 convergence, the starting point is not to bound the local error and, instead, work with an explicit formula for the global error (see \cref{lemglobalerrorintegralformula}), similarly to the way it is done for approximations of stochastic differential equations. In this way, we obtain the bound (see \cref{Etjbasicbound})
\begin{multline}
    \label{Etjbasicboundintro}
        \|X_{t_j} - X_{t_j}^N\| \leq \left( \|X_0 - X_0^N\| + c_L \int_0^{t_j} \|X_s - X_{\tau^N(s)}\| \;\mathrm{d}s \right. \\
        \left. \left\|\int_0^{t_j} \left( f(s, X_{\tau^N(s)}^N, Y_s) - f(\tau^N(s), X_{\tau^N(s)}^N, Y_{\tau^N(s)}) \right)\;\mathrm{d}s\right\|\right) e^{c_L t_j},
\end{multline}
for $j = 1, \ldots, N,$ where $\tau^N(t) = \max_{t_j \leq t} t_j$ is the running time of the Euler method, i.e. a piecewise constant function jumping along the mesh points $t_j$ (see \cref{tauNt}). The constant $c_L$ is the Lipschitz constant of $f=f(t, x, y)$ with respect to $x,$ uniformly in $(t, y),$ a typical assumption for this kind of strong convergence result.

The only problematic, noise-sensitive term in \cref{Etjbasicboundintro} is the last one. The classical analysis is to use an assumed $\theta$-H\"older regularity of the noise sample paths and bound the local strong error by $C\Delta t_N^{\theta}.$ Instead, we look at the global noise error and show that the steps of the process given by $F_t = f(t, X_{\tau^N(t)}^N, Y_t)$ can be controlled in a suitable global way. In order to give the main idea, let us consider a scalar equation with a scalar noise and assume that the sample paths of $\{F_t\}_{t\in I}$ satisfy
\[
    F_s - F_\tau = \int_\tau^s \;\mathrm{d}F_\xi.
\]
For a semi-martingale noise, the integral above can be a combination of a Lebesgue-Stieltjes integral, a stochastic integral, and an integral with respect to a jump measure, but we leave it like that just for the sake of explaining the main idea. With that, we bound the global error term using the Fubini Theorem, taking into consideration that, despite the inner integral ranging from $\xi$ to the ``future'' time $\tau^N(\xi) + \Delta t_N$, the integrand is still non-anticipative, hence
\begin{multline*}
    \int_0^{t_j} \left( f(s, X_{\tau^N(s)}^N, Y_s) - f(\tau^N(s), X_{\tau^N(s)}^N, Y_{\tau^N(s)}) \right)\;\mathrm{d}s = \int_0^{t_j} \int_{\tau^N(s)}^s \;\mathrm{d}  F_\xi\;\mathrm{d}s \\
    = \int_0^{t_j} \int_{\xi}^{\tau^N(\xi) + \Delta t_N} \;\mathrm{d}s \;\mathrm{d} F_\xi  = \int_0^{t_j} (\tau^N(\xi) + \Delta t_N - \xi) \;\mathrm{d} F_\xi.
\end{multline*}
Then, since $0 \leq \tau^N(\xi) + \Delta t_N - \xi \leq \Delta t_N$, we find that
\begin{multline*}
    \mathbb{E}\left[\left\| \int_0^{t_j} \left( f(s, X_{\tau^N(s)}^N, Y_s) - f(\tau^N(s), X_{\tau^N(s)}^N, Y_{\tau^N(s)}) \right)\;\mathrm{d}s\right\|\right] \\
    \leq \Delta t_N\mathbb{E}\left[\left\| \int_0^{t_j} \;\mathrm{d} F_\xi \right\|\right],
\end{multline*}
which yields the strong order 1 convergence. The main work is to actually show that $\int_0^{t_j} \;\mathrm{d} F_\xi$ is strongly bounded.

In the fractional Brownian motion case, which is not a semi-martingale, the argument above does not apply. Instead, we essentially have (see \cref{stepfBm} and \cref{BHtintegralformulastep}), in the case $0 < H < 1/2,$
\[
    F_s - F_\tau = \int_\tau^s g(\xi) (s-\tau)^{H-1/2}\;\mathrm{d}W_\xi + \text{higher order term},
\]
for a suitable function $g=g(\xi)$. Notice the kernel $(s-\tau)^{H-1/2}$ is unbounded when $0<H<1/2.$ In this case, disregarding the higher order term,
\begin{multline*}
    \int_0^{t_j} \left( f(s, X_{\tau^N(s)}^N, Y_s) - f(\tau^N(s), X_{\tau^N(s)}^N, Y_{\tau^N(s)}) \right)\;\mathrm{d}s \\ 
    = \int_0^{t_j} \int_{\tau^N(s)}^s g(\xi) (s-\tau^N(s))^{H-1/2} \;\mathrm{d} W_\xi\;\mathrm{d}s = (\Delta t_N)^{H+1/2} \int_0^{t_j} g(\xi)\;\mathrm{d} W_\xi,
\end{multline*}
which yields a strong convergence of order $H + 1/2,$ better than the previously known rate of $H.$ For $1/2 \leq H < 1,$ there are other terms of order 1, so that we obtain order $\min\{H+1/2, 1\},$ for $0 < H < 1.$

\section{Pathwise solutions and main assumptions}
\label{secpathwisesolution}

For the notion and main results on pathwise solutions for RODEs, we refer the reader to \cite[Section 2.1]{HanKloeden2017} and \cite[Section 3.3]{NeckelRupp2013}. We start with a fundamental set of conditions that imply the existence and uniqueness of pathwise solutions of the RODE \cref{rodeeq} in the sense of Carath\'eodory. In what follows, $\|\cdot\|$ denotes de Euclidian norm for vectors and the corresponding operator norm for matrices.
\begin{stdhyp}
    \label{standinghypotheses1}
    Let $d, k\in \mathbb{N}$ and $I=[0, T]$, where $T > 0.$
    We consider a function $f=f(t, x, y)$ defined on $I\times \mathbb{R}^d\times\mathbb{R}^k$ and with values in $\mathbb{R}^d,$ a initial condition $X_0$, and an $\mathbb{R}^k$-valued stochastic process {$\{Y_t\}_{t\in I}.$} We make the following standing assumptions:
    \begin{enumerate}
        \item \label{standinghypothesisfcontinuous} The function $f = f(t, x, y)$ is continuous in $(t, x, y)\in I \times \mathbb{R}^d \times \mathbb{R}^k.$
        \item \label{standinghypothesisLipschitzbound} The mapping $x \mapsto f(t, x, y)$ is globally Lipschitz continuous in $x$, uniformly in $t\in I$ and locally in $y\in \mathbb{R}^k,$ i.e.
        \begin{equation}
            \label{Ltassumptionbasic}
            \|f(t, x_1, y) - f(t, x_2, y)\| \leq L(\|y\|) \|x_1 - x_2\|,
        \end{equation}
        for every $t\in I,$ $y\in\mathbb{R}^k,$ and $x_1, x_2 \in\mathbb{R}^d,$ where $L=L(\eta)$ is a monotonic nondecreasing real-valued function defined for $\eta \geq 0.$
        \item \label{standinghypothesisnoise} The stochastic vector process $\{Y_t\}_{t\in I}$ is an adapted semi-martingale on a right-continuous, filtered, and complete probability space $(\Omega, \mathcal{F}, (\mathcal{F}_t)_{t\geq 0}, \mathbb{P}),$ i.e. where $\mathcal{F}_0$ contains all the null sets of probability distribution $\mathbb{P},$ and with a right-continuous filtration $\mathcal{F}_t = \cap_{s > t} \mathcal{F}_s$ in the $\sigma$-algebra $\mathcal{F}.$
        \item \label{standinghypothesisic} The initial condition is a random variable $X_0$ on $(\Omega, \mathcal{F}, \mathbb{P}).$
    \end{enumerate}
\end{stdhyp}

Under these assumptions, for almost every sample value $\omega$ in $\Omega$, the integral equation
\begin{equation}
    \label{integralrodeform}
    X_t(\omega) = X_0(\omega) + \int_0^t f(s, X_s(\omega), Y_s(\omega)) \;\mathrm{d}s
\end{equation}
has a unique solution $t\mapsto X_t(\omega)$, for the realizations $X_0 = X_0(\omega)$ of the initial condition and $t\mapsto Y_t(\omega)$ of the noise process (see \cite[Theorem 1.1]{CoddingtonLevinson1985}). Each pathwise solution of \cref{integralrodeform} is absolutely continuous. Moreover, the mapping $(t, \omega) \mapsto X_t(\omega)$ is (jointly) measurable (see \cite[Section 2.1.2]{HanKloeden2017} and also \cite[Lemma 4.51]{AB2006}) and, hence, {gives} rise to a well-defined stochastic process $\{X_t\}_{t\in I}.$ Equation \eqref{integralrodeform} is the integral form of equation \eqref{rodeeq}, so that $\{X_t\}_{t\in I}$ is a solution of \eqref{rodeeq} in the pathwise sense.

\begin{remark}
    Since $\{Y_t\}_{t\in I}$ is a semi-martingale, the sample paths $t \mapsto Y_t(\omega)$ are almost surely continuous except for a set of measure zero (see \cref{secsubmartingale}). Then, since $f=f(t, x, y)$ is continuous, the integrand $s \mapsto f(s, X_s(\omega), Y_s(\omega))$ is (Riemann) integrable, and each pathwise solution is absolutely continuous.
\end{remark}

\begin{remark}
    Since the paths of a semi-martingale are c\`adl\`ag (see the definition in \cref{secsemimartingale}), the noise process is in fact progressively measurable with respect to the filtration $(\Omega, \mathcal{F}, (\mathcal{F}_t)_{t\geq 0}, \mathbb{P}).$ In turn, the solution process $\{X_t\}_{t\in I}$ is adapted and hence also progressively measurable to the filtration $(\Omega, \mathcal{F}, (\tilde{\mathcal{F}}_t)_{t\geq 0}, \mathbb{P}),$ where each $\tilde{\mathcal{F}}_t$ is the smallest $\sigma$-algebra generated by the union of $\mathcal{F}_t$ and the smallest $\sigma$-algebra that makes the initial condition $X_0$ measurable.
\end{remark}

Adding and subtracting $f(s, 0, Y_s)$ and using assumption \eqref{standinghypothesisLipschitzbound} of the \cref{standinghypotheses1} yield the pathwise bound
\[
    \|X_t\| \leq \|X_0\| + \int_0^t \left(\|f(s, 0, Y_s)\| + L(\|Y_s\|)\|X_s\|\right)\;\mathrm{d}s,
\]
almost surely. Using Gronwall's inequality, we obtain the following result.
\begin{lemma}
    Under the \cref{standinghypotheses1}, almost every solution $\{X_t\}_{t\in I}$ of \cref{integralrodeform} satisfies
    \begin{equation}
        \label{XtboundLXMt}
        \|X_t\| \leq \left(\|X_0\| + \int_0^t \|f(s, 0, Y_s)\|\;\mathrm{d}s\right) e^{\int_0^t L(\|Y_s\|)\;\mathrm{d}s}, \quad \forall t\in I.
    \end{equation}
\end{lemma}

This is a useful pathwise bound on its own and is particularly useful for proving pathwise convergence. For the strong convergence, however, we need further assumptions. We assume that $X_0$ has finite expectation, that $\|f(t, 0, Y_t)\|$ has integrable expectation, and that $L$ is bounded, i.e.
\begin{equation}
    \label{EX0strongbound}
    \mathbb{E}[\|X_0\|] < \infty,
\end{equation}
\begin{equation}
    \label{CMmeanbound}
    \int_0^T \mathbb{E}\left[\|f(t, 0, Y_t)\|\right] \;\mathrm{d}t < \infty,
\end{equation}
and there exists a constant $c_L \geq 0$ such that
\begin{equation}
    \label{LtLXbound}
    \sup_{\eta \geq 0} L(\eta) \leq c_L < \infty,
\end{equation}
so that
\begin{equation}
    \|f(t, x_1, y) - f(t, x_2, y)\| \leq c_L \|x_1 - x_2\|,
\end{equation}
for every $t\in I,$ $y\in\mathbb{R}^k,$ and $x_1, x_2 \in\mathbb{R}^d.$ With such assumptions, it is straightforward, from taking the expectation of \cref{XtboundLXMt}, to obtain the following result.

\begin{lemma}
    \label{lemstrongbound}
    Under the \cref{standinghypotheses1}, assume that \cref{EX0strongbound}, \cref{CMmeanbound} and \cref{LtLXbound} hold. Then
    \begin{equation}
        \label{EXtstrongbound}
        \mathbb{E}[\|X_t\|] \leq \left(\mathbb{E}[\|X_0\|] + \int_0^t \mathbb{E}[\|f(s, 0, Y_s)\|]\;\mathrm{d}s\right) e^{c_L t}, \quad t\in I.
    \end{equation}
\end{lemma}

\begin{remark}
    In a number of cases, depending on the structure of the equation, we might not need all the previous assumptions, especially the restrictive assumption \eqref{LtLXbound}, corresponding to a global Lipschitz condition in $x,$ uniformly in $t$ and $y$. Indeed, when some bounded, positively invariant region exists and is of interest, we may truncate the nonlinear term to achieve the desired global conditions. Within the invariant region, the equation with the truncated term coincides with the original equation, and so do their Euler approximations, for a sufficiently refined mesh. In other situations, higher order moments of the noise and/or of the initial condition may compensate the lack of some of these assumptions and still yield bounds for the expectation of the norm of the solution. In the examples described in the Introduction, the linear models \eqref{egearthquase} and \eqref{eglinearfbm} are the only ones that satisfy the global condition \eqref{LtLXbound}. All the other examples fall within the framework of having a bounded, positively invariant region. An example in which neither applies but bounds on higher order moments compensate is, for example, a linear equation with an unbounded multiplicative noise, which can be found in the github repository \cite{RODEConvEM2023}. Finally, it is worth mentioning that the global condition \eqref{LtLXbound} is the classical assumption used for the strong convergence of approximatinos of stochastic differential equations \cite{AsaiKloeden2016,GruneKloeden2001,HanKloeden2017,HighamKloeden2021,JentzenKloeden2011,KloedenJentzen2007}. 
\end{remark}

\section{Bound on the pathwise approximations}
\label{secpathwiseapproximation}

We also need bounds for the Euler approximation \cref{emscheme} analogous to \cref{XtboundLXMt} and \cref{EXtstrongbound}. Notice that
\begin{align*}
    \|X_{t_j}^N\| & \leq \|X_{t_{j-1}}^N\| + \Delta t_N \|f(t_{j-1}, X_{t_{j-1}}^N, Y_{t_{j-1}})\| \\
    & \leq \|X_{t_{j-1}}^N\| + \Delta t_N (\|f(t_{j-1}, 0, Y_{t_{j-1}})\| + L(\|Y_{t_{j-1}}\|)\|X_{t_{j-1}}\|) \\
    & \leq \left(1 + L(\|Y_{t_{j-1}}\|)\Delta t_N\right)\|X_{t_{j-1}}^N\| + \Delta t_N \|f(t_{j-1},0,Y_{t_{j-1}})\|.
\end{align*}
We bound $(1 + a) \leq e^a$, so that
\[
    \|X_{t_j}^N\| \leq e^{L(\|Y_{t_{j-1}}\|)\Delta t_N}\|X_{t_{j-1}}^N\| + \Delta t_N \|f(t_{j-1},0,Y_{t_{j-1}})\|.
\]
Iterating this relation we obtain the following bound.
\begin{equation}
    \label{XNtboundLXMtpre}
    \|X_{t_j}^N\| \leq \left(\|X_0^N\| + \Delta t_N \sum_{i=0}^{j-1} \|f(t_i,0,Y_{t_i})\|\right)e^{\Delta t_N \sum_{i=0}^{j-1} L(\|Y_{t_i}\|)},
\end{equation}
for $j = 0, \ldots, N.$ The right hand side above depends on the time mesh. A straighforward bound independent of the time mesh is obtained from \cref{XNtboundLXMtpre} and the fact that semi-martingales are pathwise locally bounded, almost surely:
\begin{equation}
    \label{XNtboundLXMt}
    \|X_{t_j}^N\| \leq \left(\|X_0^N\| + t_j \sup_{0\leq s\leq t_j} \|f(t_s,0,Y_s)\|\right)e^{t_j \sup_{0\leq s\leq t_j} L(\|Y_s\|)},
\end{equation}
for $j = 0, \ldots, N.$ With that, we obtain the following result, upon taking the expectation and using the assumed bounds.
\begin{lemma}
    \label{lemXNtmeanboundLXMt}
    Under the \cref{standinghypotheses1}, assume that \cref{EX0strongbound}, \cref{CMmeanbound} and \cref{LtLXbound} hold. Then, the Euler approximation \cref{emscheme} for a pathwise solution of the random ordinary differential equation \cref{rodeeq} satisfies
    \begin{equation}
        \label{XNtmeanboundLXMt}
        \mathbb{E}\left[\|X_{t_j}^N\|\right] \leq \left(\mathbb{E}\left[\|X_0^N\|\right] + \mathbb{E}\left[\sum_{i=0}^{j-1} \|f(t_i,0,Y_{t_i})\|\Delta t_N \right] \right)e^{c_L t_j},
    \end{equation}
    for all $j = 0, \ldots, N.$
\end{lemma}

\section{Semi-martingale noises}
\label{secsemimartingale}

For the theory of semi-martingales, we refer the reader to \cite{AitSahaliaJacod2014,HeWangYan1992,KaratzasShreve2014,Metivier1982,Protter2005}. Starting with the approach of \cite{Protter2005}, the main idea is that \emph{semi-martingales} are stochastic processes {$\{Y_t\}_{t\in I}$} which are c\`adl\`ag, adapted, and ``good integrators", in the sense that the sums
\[
    I_Y(H) = H_0 Y_0 + \sum_{i=1}^n H_i(Y_{t_i}-Y_{t_{i-1}})
\]
are continuous functions from the space of \emph{simple predictable} processes
\[
    H_t = H_0 \chi_{0}(t) + \sum_{i=1}^n H_j \chi_{(T_{i-1}, T_i]},
\]
to the space of a.s.~finite-valued random variables, where $0=T_0 < T_1 < \ldots < T_n,$ are stopping times, $H_j \in \mathcal{F}_{T_j},$ and $|H_j| < \infty$ almost surely, with the former space endowed with the topology of uniform convergence in $(t, \omega)\in [0,\infty)\times \Omega,$ and the latter, with the topology of convergence in probability. We recall the definitions and some facts about (scalar or vector-valued) semi-martingale processes.

A function defined on a real interval is called \emph{c\`adl\`ag} when it is right continuous and with left limits (``continue \`a droite, limite \`a gauche", in French). A \emph{c\`adl\`ag process} is a stochastic process $\{Y_t\}_{t\in I}$ whose sample paths are almost surely c\`adl\`ag paths. Thus, $\{Y_t\}_{t\in I}$ is a c\`adl\`ag process when for almost every $\omega\in \Omega$ and for all $t\in I$, we have that $\lim_{s\rightarrow t^+} Y_s(\omega) = Y_t(\omega)$ (right continuous) and the limit $\lim_{s \rightarrow t^-} Y_s(\omega)$ exists (left limit).

For a c\`adl\`ag process such as semi-martingales, we define the left-limit process $Y_{t^{-}} = \lim_{t \rightarrow t^-} Y_s$ and the jump process $\Delta Y_t = Y_t - Y_{t^{-}},$ with the usual convention that $Y_{0^-} = 0,$ and thus $\Delta Y_0 = Y_0.$

Almost surely, the sample paths of the jump process are zero except possibly at a countable number of points, denoted $J_{Y}$, i.e. for almost every $\omega\in \Omega$, there exists a countable set $J_{Y}(\omega) \subset I$ for which $\Delta Y_t(\omega)$ is zero on each connected component of $I\setminus J_{Y}(\omega)$. In fact, for any given $\varepsilon > 0$, the set of times for which the jump $\|\Delta Y_t\|$ is larger than or equal to $\varepsilon$ is finite. As a consequence, each sample path is continuous except on the countable null set $J_{Y}.$

Another important process associated with a semi-martingale is the \emph{quadratic variation} $\{[Y, Y]_t\}_{t\in I}$. It is defined either as a suitable limit of the summations of quadratic steps $(Y_{T_j^n} - Y_{T_{j-1}^n})^2$ over random partitions made of stopping times $T_j^n$ or, equivalently, as a suitable integral \cite[Section II.6]{Protter2005}. In the vector-valued case, denoting $Y_t = (Y_t^{(l)})_{l=1, \ldots, k}$, with the vectorial components $Y_t^{(l)},$ the quadratic process becomes matrix-valued, i.e. $[Y, Y]_t = ([Y^{(l_1)}, Y^{(l_2)}]_t)_{l_1,l_2=1,\ldots, k},$ with the appearance of the quadratic variations and covariations. The equivalent integral definition takes, more precisely, the form of an ``integration by parts" formula for two semi-martingales,
\[
    [Y^{(l_1)}, Y^{(l_2)}]_t = Y_t^{(l_1)} Y_t^{(l_2)} - \int_0^t Y_{s^-}^{(l_1)} \;\mathrm{d}Y_s^{(l_2)} - \int_0^t Y_{s^-}^{(l_2)} \;\mathrm{d}Y_s^{(l_1)}.
\]

The quadratic variation is a c\`adl\`ag, increasing, and adapted process \cite[Theorem II.22]{Protter2005}. The \emph{continuous part} {$\{[Y, Y]_t^c\}_{t\in I}$} of the quadratic variation is defined by
\[
    [Y, Y]_t^c = [Y, Y]_t - \sum_{0\leq s \leq t} \Delta Y_s\Delta Y_s^\tr,
\]
which itself is a semi-martingale of finite-variation type, and where $Y_t^\tr$ is the transpose of $Y_t.$ Notice that $[Y, Y]_0 = Y_0^2,$ while, for the continuous part, $[Y, Y]_0^c = [Y, Y]_0 - \Delta Y_0\Delta Y_0^\tr = Y_0 Y_0^\tr - Y_0 Y_0^\tr = 0.$

The Bichteler-Dellacherie Theorem \cite[Theorem III.47]{Protter2005} yields a characterization of a semi-martingale as a c\`adl\`ag process $\{Y_t\}_{t\in I}$ decomposable as a sum of a finite-variation (FV) process and a local martingale process, which we can write in standard form as $Y_t = Y_0 + F_t + Z_t$, with an FV process $\{F_t\}_{t\in I}$ and a local martingale process $\{Z_t\}_{t\in I}$ with $F_0 = Z_0 = 0.$

A \emph{Finite-Variation process} (\emph{FV process}) is a c\`adl\`ag process $\{F_t\}_{t \in I}$ whose sample paths are almost surely of finite variation $V({\{F_t\}_{t\in I}}; K)$ on compact subsets $K\subset I,$ i.e.
\[ V(\{F_t\}_{t\in I}; K) = \sup_{t_0 < \ldots < t_n \in K, n\in \mathbb{N}} \sum_{j=1}^n \|F_{t_j} - F_{t_{j-1}}\| < \infty.
\]
This total variation induces a positive measure $\|\mathrm{d}F_t\|,$ which is the infinitesimal of this total variation, so that $\|\mathrm{d}F_t\|([a, b]) = V({\{F_t\}_{t\in I}}; [a, b]),$ for every $[a, b]\subset I.$

The class of FV processes includes the compound Poisson process, renewal-reward processes, inhomogeneous Poisson process, Hawkes self-exciting process, and many more. It also includes transport processes $Y_t = g(t, Y_0)$ where $Y_0$ is a random variable and $t \mapsto g(t, y_0)$ is a c\`adl\`ag function of finite variation, for every sample value $y_0=Y_0(\omega)$.

A \emph{local martingale} is an adapted, c\`adl\`ag process $\{Z_t\}_{t\in I}$ such that there exists a sequence of increasing stopping times $T_n$ with $T_n \rightarrow \infty$ almost surely, for which $\{Z_{t \wedge T_n}1_{T_n > 0}\}_{t \in I}$ is a uniformly integrable martingale, for each $n\in\mathbb{N}.$ Here, $t \land T = \min\{t, T\},$ the function $1_{T_n > 0}$ is the indicator function of the set $\{\omega; \;T_n(\omega) > 0\},$ and uniform integrability means
\[
  \lim_{k\rightarrow \infty}\sup_{t \in I} \int_{\{T_n > 0\}\cap\{\|Z_{t\wedge T_n}\| \geq k\}}\|Z_{t\wedge T_n}\|\;\mathrm{d}\mathbb{P} = 0.
\]

The class of local martingales includes time-changed Brownian motions and It\^o diffusion processes. For more information on It\^o diffusion, see also \cite{Oksendal2003}.

With the decomposition $Y_t = Y_0 + F_t + Z_t$ described above, we recall from \cite[Lemma 18.7]{Metivier1982} that
\[
    [Y, Y]_t = [Z, Z]_t + \sum_{0 \leq s \leq t} \Delta F_s \Delta F_s^\tr.
\]
Each step is bounded by the variation of $\{F_t\}_{t\in I},$ i.e.
\begin{multline*}
    \sum_{0 \leq s \leq t} \|\Delta F_s\|^2 = \sum_{0 \leq s \leq t} \|F_s - F_{s^-}\|^2 \leq \sum_{0 \leq s \leq t} (\|F_s\| + \|F_{s^-}\|)\|F_s - F_{s^-}\| \\
    \leq 2\sup_{t\in I} \|F_t\| \sum_{0 \leq s \leq t} \|F_s - F_{s^-}\|.
\end{multline*}
Since the number of jumps is countable, we find, by an approximation procedure (taking sequences $s_j$ converging from below to each jump point $s$), that
\[
    \sum_{0 \leq s \leq t} \|\Delta F_s\|^2 \leq \sup_{t\in I} \|F_t\| \|V(\{F_t\}_{t\in I}; I)\| \leq \|V(\{F_t\}_{t\in I}; I)\|^2.
\]
We also have
\begin{equation}
    \label{estimatejumpfv}
    \sum_{s\in J_{Y}} \|\Delta F_s\| \leq V(\{F_t\}_{t\in I}; I) < \infty.
\end{equation}

In the more general case of a semi-martingale, the process $\{Y_t\}_{t\in I}$ has finite quadratic variation \cite[Section II.6]{Protter2005} on the bounded interval $I=[0, T],$ and we find, using that $[Y, Y]_0 = Y_0 Y_0^\tr$ and $\Delta [Y, Y]_t = \Delta Y_t \Delta Y_t^\tr$ (see \cite[Theorem II.22]{Protter2005}), that
\begin{equation}
    \label{estimatejumpsemimartingale}
    \sum_{s\in J_{Y}} \|\Delta Y_s\|^2 \leq \|[Y, Y]_T\| < \infty.
\end{equation}

By Doob's maximal quadratic inequality (see e.g. \cite[Theorem I.20]{Protter2005}), if $Z_t = Z_{t \land T},$ $0\leq t \leq T,$ is a square-integrable martingale, then
\begin{equation}
    \label{expectationZtsquaredfinite}
    \mathbb{E}\left[\sup_{t\in I} \|Z_t\|^2\right] \leq 4\mathbb{E}[\|Z_T\|^2] < \infty,
\end{equation}

Since $F_0 = 0,$ the finite variation part is bounded by
\[
    \sup_{t\in I} \|F_t\| \leq V(\{F_t\}_{t\in I}; I),
\]
which implies, in particular,
\[
    \mathbb{E}\left[\sup_{t\in I} \|F_t\|\right] \leq \mathbb{E}\left[V(\{F_t\}_{t\in I}; I)\right], \qquad \mathbb{E}\left[\sup_{t\in I} \|F_t\|^2\right] \leq \mathbb{E}\left[V(\{F_t\}_{t\in I}; I)^2\right].
\]
For $Y_t = Y_0 + Z_t + F_t,$ we find
\begin{equation}
    \label{expectationYtfinite}
    \mathbb{E}\left[\sup_{t\in I} \|Y_t\|\right] \leq \mathbb{E}[\|Y_0\|] + \mathbb{E}\left[V(\{F_t\}_{t\in I}; I)\right] + 2\mathbb{E}\left[\|Z_T\|^2\right]^{1/2},
\end{equation}
and
\begin{equation}
    \label{expectationYtsquaredfinite}
    \mathbb{E}\left[\sup_{t\in I} \|Y_t\|^2\right] \leq 3\mathbb{E}[\|Y_0\|^2] + 3\mathbb{E}\left[V(\{F_t\}_{t\in I}; I)^2\right] + 12\mathbb{E}[\|Z_T\|^2].
\end{equation}

For the quadratic variation part of a local martingale $\{Z_t\}_{t\in I},$ the Burkholder-David-Gundy inequalities (see e.g. \cite[Theorem IV.48]{Protter2005} and \cite{MARINELLI20161}) relate the moments of the maximum of the local martingale with moments of its quadratic variation,
\[
    \frac{1}{\gamma_p}\mathbb{E}\left[\|[Z, Z]_t\|^{p/2}\right] \leq \mathbb{E}\left[\max_{0\leq s \leq t} \|Z_s\|^p\right] \leq \tilde{\gamma}_p\mathbb{E}\left[\|[Z, Z]_t\|^{p/2}\right],
\]
for $p \geq 1$ and suitable constants $\gamma_p, \tilde{\gamma}_p > 0.$ Thus, we obtain the bound
\begin{equation}
    \label{expectationofquadraticvariationofZ}
    \mathbb{E}\left[\|[Z, Z]_t\|\right] \leq 4\gamma_2\mathbb{E}\left[\|Z_T\|^2\right] < \infty,
\end{equation}
for any $t\in I,$ and then
\begin{equation}
    \label{expectationofquadraticvariation}
    \mathbb{E}\left[\|[Y, Y]_t\|\right] \leq 4\gamma_2\mathbb{E}\left[\|Z_T\|^2\right] + \mathbb{E}\left[\|V(\{F_t\}_{t\in I}; I)\|^2\right].
\end{equation}

This, together with \cref{estimatejumpsemimartingale}, also yields
\begin{equation}
    \label{expectationsquaredjumps}
    \mathbb{E}\left[\sum_{s\in J_{Y}} \|\Delta Y_s\|^2 \right] \leq 4\gamma_2\mathbb{E}\left[\|Z_T\|^2\right] + \mathbb{E}\left[\|V(\{F_t\}_{t\in I}; I)\|^2\right].
\end{equation}

The polarization identity
\[ [Y^{(l_1)}, Y^{(l_2)}] = \frac{1}{2}\left([Y^{(l_1)} + Y^{(l_2)}, Y^{(l_1)} + Y^{(l_2)}] - [Y^{(l_1)}, Y^{(l_1)}] - [Y^{(l_2)}, Y^{(l_2)}]\right)
\]
implies that each quadratic covariation is a difference of quadratic variations, which are increasing processes, so that each quadratic covariation is also of finite variation (see \cite[Corollary 1 of Theorem II.22]{Protter2005}). In particular, the continuous part of the quadratic local martingale term $\{Z_t\}_{t\in I}$ is of finite variation and it follows from the definition of the continuous part of the quadratic variation, \cref{expectationofquadraticvariation}, and \cref{expectationofquadraticvariationofZ} that
\begin{equation}
    \label{eqquadraticvarcovarfiniteexpectation}
    \mathbb{E}\left[V(\{[Z, Z]_t^c\}_{t\in I}; I)\right] \leq 5\gamma_2\mathbb{E}\left[\|Z_T\|^2\right] < \infty.
\end{equation}

A fundamental tool for integration and differentiation is the chain rule. In the case of an FV process, given a continuously differentiable function $f:\mathbb{R}^k \rightarrow \mathbb{R}^d$, the following change of variables formula holds (see \cite[Theorems II.31 and II.33]{Protter2005}):
\begin{equation}
    f(Y_t) - f(Y_0) = \int_{0^+}^t Df(Y_{s^-}) \;\mathrm{d}Y_s + \sum_{0 < s \leq t} \left( f(Y_s) - f(Y_{s^{-}}) - Df(Y_{s^-})\Delta Y_s\right).
\end{equation}
The summation on the right hand side is an at most countable summation since the summand vanishes for $s \notin J_{Y}$. The summation can be written as an integral with respect to a jump measure (see \cref{changeofvariablesformulasemimartingale} for the full expression), but we keep it as a summation for simplicity.

The formula above is written in vectorial form, with $f$ and $Y_t$ of the form $f(y)=(f_i(y_1, \ldots, y_k))_{i=1, \ldots, d}$ and $Y_t = (Y_t^{(l)})_{l=1, \ldots, k}$, and
\[
    Df(y)\Delta Y_t = \left( \nabla f_i(y)\Delta Y_t\right)_{i=1, \ldots, d} = \left( \sum_{l=1, \ldots, k} \partial y_l f_i(y) \Delta Y_t^{(l)} \right)_{i=1, \ldots, d}.
\]

When $f=f(t, y)$ depends also on $t$, one can apply the above formula to the extended process $\{(t, Y_t)\}_{t\in I}$. The first component is continuous, hence it has no jumps, and the corresponding summation term vanishes. This leads to the following change of variables formula for FV processes:
\begin{multline}
    \label{changeofvariablesformulacadlagfv}
    f(t, Y_t) - f(0, Y_0) = \int_0^t \partial_s f(s, Y_{s^-})\;\mathrm{d}s + \int_{0^+}^t D_y f(s, Y_{s^-}) \;\mathrm{d}Y_s \\
    + \sum_{0 < s \leq t} \left( f(s, Y_s) - f(s, Y_{s^{-}}) - D_y f(s, Y_{s^-})\Delta Y_s\right).
\end{multline}

In the more general case of a semi-martingale noise $\{Y_t\}_{t\in I}$, the change of variables formula \cref{changeofvariablesformulacadlagfv} takes the form (see e.g. \cite[Theorems II.32 and II.33]{Protter2005}, \cite[Section 1.3.3]{AitSahaliaJacod2014}, \cite[Theorem 3.6]{KaratzasShreve2014}, or \cite[Theorem 9.35]{HeWangYan1992})
\begin{equation}
    \label{changeofvariablesformulasemimartingale}
    \begin{aligned}
        f(t, Y_t) - f(0, Y_0) & = \int_0^t \partial_s f(s, Y_{s^-})\;\mathrm{d}s + \int_{0^+}^t D_y f(s, Y_{s^-}) \;\mathrm{d}Y_s \\
        & \quad + \sum_{0 < s \leq t} \left( f(s, Y_s) - f(s, Y_{s^{-}}) - D_y f(s, Y_{s^-})\Delta Y_s\right) \\
        & \quad + \frac{1}{2}\int_0^t D_{yy}f(s, Y_{s^-})\;\mathrm{d}[Y, Y]_s^c,
    \end{aligned}
\end{equation}
where each component of $D_{yy}f(s, Y_{s^-})\;\mathrm{d}[Y, Y]_s^c$ is given by
\[ \left(D_{yy}f(s, Y_{s^-})\;\mathrm{d}[Y, Y]_s^c\right)^{(i)} = \sum_{l_1, l_2 = 1, \ldots, d} \partial_{y_{l_1}}\partial_{y_{l_2}} f_i(s, Y_{s^-})\;\mathrm{d}[Y^{(l_1)}, Y^{(l_2)}]_s^c.
\]

Since $[Y, Y]_t^c$ is of finite variation, the last integral in \cref{changeofvariablesformulasemimartingale} is a Lebesgue-Stieltjes integral.

We actually apply the formula \cref{changeofvariablesformulasemimartingale} to a function of the form $f=f(t, x, y)$, but with $x$ constant, so the above formula suffices.

Finally, a fundamental identity is the It\^o isometry, which, in the case of a vector-valued local martingale $\{Z_t\}_{t\in I}$ and an adapted matrix-valued process $\{A_t\}_{t\in I},$ takes the form (e.g. apply \cite[Chapter 3, Proposition 2.17]{KaratzasShreve2014} to the expanded coordinate terms)
\begin{multline}
    \label{eqitoisometrylocalmartingale}
    \mathbb{E}\left[\left\|\int_0^t A_s \;\mathrm{d}Z_s\right\|^2\right] = \mathbb{E}\left[ \int_0^t A_s^\tr A_s \odot \mathrm{d}[Z, Z]_s\right] \\ = \mathbb{E}\left[ \int_0^t \sum_{l_1, l_2 = 1}^k \sum_{i=1}^d A_s^{(il_1)}A_s^{(il_2)}\;\mathrm{d}[Z^{(l_1)}, Z^{(l_2)}]_s\right],
\end{multline}
provided each component $\{A_t^{il}\}_{t\in I}$ is square integrable for the integration with respect to $\{[Z^{(l)}, Z^{(l)}]_t\}_{t\in I}$, and where $\odot$ denotes the Hadamard product of two matrices of same dimensions. Recall that the steps of the cross variation matrix $[Z, Z]_s$ are positive semi-definite, so the right hand side is nonnegative.

\section{Integral formula for the global pathwise error}

In this section, we derive the following integral formula for the global error:
\begin{lemma}
    \label{lemglobalerrorintegralformula}
    Under the \cref{standinghypotheses1}, the Euler approximation \cref{emscheme} for a pathwise solution of the random ordinary differential equation \cref{rodeeq} satisfies almost surely the global error formula
    \begin{equation}
        \label{globalerrorintegralformula}
        \begin{aligned}
            X_{t_j} - X_{t_j}^N & = X_0 - X_0^N + \int_0^{t_j} \left( f(s, X_s, Y_s) - f(s, X_{\tau^N(s)}, Y_s) \right)\;\mathrm{d}s  \\ 
            & \qquad + \int_{0}^{t_j} \left( f(s, X_{\tau^N(s)}, Y_s) - f(s, X_{\tau^N(s)}^N, Y_s) \right)\;\mathrm{d}s \\
            & \qquad + \int_0^{t_j} \left( f(s, X_{\tau^N(s)}^N, Y_s) - f(\tau^N(s), X_{\tau^N(s)}^N, Y_{\tau^N(s)}) \right)\;\mathrm{d}s,
        \end{aligned}
    \end{equation}
    for $j = 1, \ldots, N,$ where $\tau^N$ is the piecewise constant function along the time mesh which gives the running time of the Euler method:
    \begin{equation}
        \label{tauNt}
        \tau^N(t) = \max_{t_j \leq t}\{t_j\} = \left\lfloor\frac{t}{\Delta t_N}\right\rfloor\Delta t_N = \left\lfloor\frac{tN}{T}\right\rfloor\frac{T}{N}.
    \end{equation}
\end{lemma}

\begin{proof}
    Under the \cref{standinghypotheses1}, the solutions of \cref{rodeeq} are pathwise solutions in the sense of \cref{integralrodeform}. With that in mind, we first obtain an expression for a single time step, from time $t_{j-1}$ to $t_j = t_{j-1} + \Delta t_N$.
    
    Subtracting the Euler step from the pathwise solution on a single time step yields
    $$
    X_{t_j} - X_{t_j}^N = X_{t_{j-1}} - X_{t_{j-1}}^N + \int_{t_{j-1}}^{t_j} \left( f(s, X_s, Y_s) - f(t_{j-1}, X_{t_{j-1}}^N, Y_{t_{j-1}}) \right)\;\mathrm{d}s.
    $$

    Adding and subtracting appropriate terms yield
    \begin{equation}
        \label{singlestep}
        \begin{aligned}
            X_{t_j} - X_{t_j}^N  = & X_{t_{j-1}} - X_{t_{j-1}}^N \\
            & + \int_{t_{j-1}}^{t_j} \left( f(s, X_s, Y_s) - f(s, X_{t_{j-1}}, Y_s) \right)\;\mathrm{d}s \\ 
            & + \int_{t_{j-1}}^{t_j} \left( f(s, X_{t_{j-1}}, Y_s) - f(s, X_{t_{j-1}}^N, Y_s) \right)\;\mathrm{d}s \\
            & + \int_{t_{j-1}}^{t_j} \left( f(s, X_{t_{j-1}}^N, Y_s) - f(t_{j-1}, X_{t_{j-1}}^N, Y_{t_{j-1}}) \right)\;\mathrm{d}s.
        \end{aligned}
    \end{equation}

    Now we iterate the time steps \cref{singlestep} to find that
    \begin{align*}
        X_{t_j} - X_{t_j}^N & = X_0 - X_0^N + \sum_{i=1}^{j} \left(\int_{t_{i-1}}^{t_i} \left( f(s, X_s, Y_s) - f(s, X_{t_{i}}, Y_s) \right)\;\mathrm{d}s \right. \\ 
        & \qquad + \int_{t_{i-1}}^{t_i} \left( f(s, X_{t_{i-1}}, Y_s) - f(s, X_{t_{i-1}}^N, Y_s) \right)\;\mathrm{d}s \\
        & \qquad \left. + \int_{t_{i-1}}^{t_i} \left( f(s, X_{t_{i-1}}^N, Y_s) - f(t_{i-1}, X_{t_{i-1}}^N, Y_{t_{i-1}}) \right)\;\mathrm{d}s \right).
    \end{align*}

    Using the running time $\tau^N$ defined by \cref{tauNt}, the above expression becomes \cref{globalerrorintegralformula}, as desired.
\end{proof}

\section{Basic bound for the global error}

Here we derive a bound for the global error which is the basis for the convergence bound. For that, we use the following classical discrete version of the Gronwall Lemma, which can be deduced from e.g. \cite[Lemma 1]{Sugiyama1969}.

\begin{lemma}[Discrete Gronwall Lemma]
    \label{lemdiscretegronwall}
    Let $(e_j)_j$ be a (finite or infinite) sequence of positive numbers starting at $j=0$ and satisfying
    \begin{equation}
        \label{integralgronwall}
        e_j \leq \sum_{i=0}^{j-1} a_i e_i + b_i,
    \end{equation}
    for every $j$, and where $a_i, b_i \geq 0$, with $b_i$ nondecreasing. Then, for all $j,$ we have
    \begin{equation}
        \label{estimategronwall}
        e_j \leq b_j e^{\sum_{i=0}^{j-1} a_i}.
    \end{equation}
\end{lemma}

With that, we are ready to prove our basic bound for the global error.
\begin{proposition}[Basic strong bound]
    \label{propbasicstrongestimate}
    Under the \cref{standinghypotheses1}, assume that \cref{LtLXbound} holds. Then, for every $j = 0, \ldots, N$,
    \begin{multline}
        \label{expectedestimateglobalerrorintegral}
        \mathbb{E}\left[\|X_{t_j} - X_{t_j}^N\|\right] \leq \Bigg( \mathbb{E}\left[\|X_0 - X_0^N\|\right] + \int_0^{t_j} c_L \mathbb{E}\left[\|X_s - X_{\tau^N(s)}\|\right] \;\mathrm{d}s \\
        + \mathbb{E}\left[\left\|\int_0^{t_j} \left( f(s, X_{\tau^N(s)}^N, Y_s) - f(\tau^N(s), X_{\tau^N(s)}^N, Y_{\tau^N(s)}) \right)\;\mathrm{d}s\right\|\right]\Bigg) e^{c_L t_j}.
    \end{multline}
\end{proposition}

\begin{proof}
    We bound the first two integrals in \cref{globalerrorintegralformula}. For the first one, we use \cref{Ltassumptionbasic} together with the assumed \cref{LtLXbound}, so that
    $$
        \|f(s, X_s, Y_s) - f(s, X_t, Y_s)\| \leq c_L \|X_s - X_t\|,
    $$
    for $t, s \in I$, and, in particular, for $t = \tau^N(s)$. Hence,
    $$
        \left\|\int_0^{t_j} \left( f(s, X_s, Y_s) - f(s, X_{\tau^N(s)}, Y_s) \right)\;\mathrm{d}s \right\| \leq c_L \int_0^{t_j} \|X_s - X_{\tau^N(s)}\| \;\mathrm{d}s.
    $$
    
    For the second term, we use again \cref{Ltassumptionbasic}, and find that
    \begin{align*}
        \left\|\int_0^{t_j} \left( f(s, X_{\tau^N(s)}, Y_s) - f(s, X_{\tau^N(s)}^N, Y_s) \right)\;\mathrm{d}s \right\| & \leq c_L \int_0^{t_j} \|X_{\tau^N(s)} - X_{\tau^N(s)}^N\| \;\mathrm{d}s \\
        & \leq c_L \sum_{i=0}^{j-1} \|X_{t_i} - X_{t_i}^N\| \Delta t.
    \end{align*}
    
    Using these two bounds in \cref{globalerrorintegralformula} yields
    \begin{multline*}
        \|X_{t_j} - X_{t_j}^N\| \leq \|X_0 - X_0^N\| + c_L \int_0^{t_j} \|X_s - X_{\tau^N(s)}\| \;\mathrm{d}s + c_L \sum_{i=0}^{j-1} \|X_{t_i} - X_{t_i}^N\| \Delta t \\
        + \left\|\int_0^{t_j} \left( f(s, X_{\tau^N(s)}^N, Y_s) - f(\tau^N(s), X_{\tau^N(s)}^N, Y_{\tau^N(s)}) \right)\;\mathrm{d}s\right\|.
    \end{multline*}
    This can be cast in the form of \cref{integralgronwall}, with a constant coefficient $a_i = c_L \Delta t.$ Then, using the discrete Gronwall \cref{lemdiscretegronwall}, we obtain from \eqref{estimategronwall} that
    \begin{multline}
        \label{Etjbasicbound}
            \|X_{t_j} - X_{t_j}^N\| \leq \left( \|X_0 - X_0^N\| + c_L \int_0^{t_j} \|X_s - X_{\tau^N(s)}\| \;\mathrm{d}s \right. \\
            + \left. \left\|\int_0^{t_j} \left( f(s, X_{\tau^N(s)}^N, Y_s) - f(\tau^N(s), X_{\tau^N(s)}^N, Y_{\tau^N(s)}) \right)\;\mathrm{d}s\right\|\right) e^{c_L t_j},
    \end{multline}
    for $j=1, \ldots, N$, where $\tau^N$ is given by \cref{tauNt}. Taking the expectation of \cref{Etjbasicbound}, we obtain the desired result \cref{expectedestimateglobalerrorintegral}.
\end{proof}

The first term in the right hand side of \cref{expectedestimateglobalerrorintegral} usually vanishes since in general we take $X_0^N = X_0$, but it suffices to assume that $X_0^N$ approximates $X_0$ to first order $\Delta t_N$, which is useful for lower order approximations or for discretizations of (random) partial differential equations.

The third term in \cref{expectedestimateglobalerrorintegral} is the more delicate one that will be handled differently in the next section.

As for the second term, which only concerns the solution itself, not the approximation, we use the following crude but useful general result.

\begin{lemma}
    Under the \cref{standinghypotheses1}, assume that \cref{EX0strongbound}, \cref{CMmeanbound} and \cref{LtLXbound} hold. Then, we find that, for $j=1, \ldots, N,$
    \begin{equation}
        \label{estimatesecondterminglobalstrongerrorintegral}
        \mathbb{E}\left[\int_0^{t_j}\left\|X_s - X_{\tau^N(s)}\right\| \;\mathrm{d}s\right] \leq \Delta t_N \left(\mathbb{E}[\|X_0\|] + \int_0^{t_j} \mathbb{E}[\|f(s, 0, Y_s)\|]\;\mathrm{d}s\right) e^{c_L t_j}.
    \end{equation}
\end{lemma}

\begin{proof}
    Using that the paths of $\{X_t\}_{t\in I}$ are almost surely absolutely continuous, adding and subtracting $f(\xi, 0, Y_\xi)$ to the integrand, and using assumption \eqref{standinghypothesisLipschitzbound} inside the \cref{standinghypotheses1} along with \cref{LtLXbound}, we have almost surely that
    \[
      \left\|X_s - X_{\tau^N(s)}\right\| = \left\|\int_{\tau^N(s)}^s f(\xi, X_\xi, Y_\xi)\;\mathrm{d}\xi\right\| \leq \int_{\tau^N(s)}^s \left(\|f(\xi, 0, Y_\xi)\| + c_L\|X_\xi\|\right)\;\mathrm{d}\xi.
    \]
    Integrating over $[0, t_j]$ and using Fubini's theorem to exchange the order of integration,
    \begin{align*}
        \int_0^{t_j} \left\|X_s - X_{\tau^N(s)}\right\| \;\mathrm{d}s & \leq \int_0^{t_j}\int_{\tau^N(s)}^s (\|f(\xi, 0, Y_\xi)\| + c_L\|X_\xi\|) \;\mathrm{d}\xi \;\mathrm{d}s \\
        & \leq \int_0^{t_j} (\tau^N(\xi) + \Delta t_N - \xi) (\|f(\xi, 0, Y_\xi)\| + c_L\|X_\xi\|) \;\mathrm{d}\xi.
    \end{align*}
    Using that $\tau^N(\xi) \leq \xi,$ we have $\tau^N(\xi) + \Delta t_N - \xi \leq \Delta t_N,$ and we obtain, upon taking the expectation
    \[
        \mathbb{E}\left[\int_0^{t_j}\left\|X_s - X_{\tau^N(s)}\right\| \;\mathrm{d}s\right] \leq \Delta t_N \int_0^{t_j} \left(\mathbb{E}\left[\|f(s, 0, Y_s)\|\right] + c_L\mathbb{E}\left[\|X_s\|\right]\right) \;\mathrm{d}s.
    \]
    Thanks to \cref{EX0strongbound}, \cref{CMmeanbound} and \cref{LtLXbound}, \cref{lemstrongbound} applies and inequality \cref{EXtstrongbound} holds, so that we obtain
    \begin{align*}
        \mathbb{E}\bigg[\int_0^{t_j}  \|X_s - X_{\tau^N(s)}\| \;\mathrm{d}s\bigg] & \leq \Delta t_N  \bigg(\int_0^{t_j} \mathbb{E}[\|f(s, 0, Y_s)\|] \;\mathrm{d}s \\
        & \quad + c_L\int_0^{t_j}\left(\mathbb{E}[\|X_0\|] + \int_0^s \mathbb{E}[\|f(\xi, 0, Y_\xi)\|]\;\mathrm{d}\xi\right)e^{c_L s} \;\mathrm{d}s\bigg) \\
        & \leq \Delta t_N  \bigg(\int_0^{t_j} \mathbb{E}[\|f(s, 0, Y_s)\|] \;\mathrm{d}s \\
        & \quad + \left(\mathbb{E}[\|X_0\|] + \int_0^{t_j} \mathbb{E}[\|f(\xi, 0, Y_\xi)\|]\;\mathrm{d}\xi\right)\left(e^{c_L t_j} - 1\right) \bigg).
    \end{align*}
    This simplifies to \cref{estimatesecondterminglobalstrongerrorintegral}.
\end{proof}

\section{Strong convergence for semi-martingale noises}
\label{secsubmartingale}

We now prove the desired strong order 1 convergence.

\begin{theorem}
    \label{thmsemimartingale}
    Consider the Random ODE \eqref{rodeeq} and the Euler approximation \eqref{emscheme}, with the function $f=f(t, x, y),$ the initial condition $X_0,$ and the noise $\{Y_t\}_{t\in I}$ satisfying the \cref{standinghypotheses1}. Suppose  that \cref{EX0strongbound}, \cref{CMmeanbound} and \cref{LtLXbound} hold; that $X_0^N$ is a first order approximation of $X_0,$ in the sense that 
    \begin{equation}
        \label{EX0X0Nunif}
        \mathbb{E}[\|X_0 - X_0^N\|] \leq c_0 \Delta t_N, \qquad \forall N\in \mathbb{N},
    \end{equation}
    for some constant $c_0 \geq 0;$ and that $f=f(t, x, y)$ is uniformly globally Lipschitz continuous in $x$, continuously differentiable with respect to $t$ and twice continuously differentiable with respect to $y$, with partial differentials $\partial_t f$, $D_y f$ and $D_{yy} f$ satisfying the bounds
    \begin{align}
        \label{ftfyunifboundcadlagfvpartialt}
        \left\|\partial_t f(t, x, y)\right\| & \leq c_1 + c_2\|x\| + c_3\|y\|, \\
        \label{ftfyunifboundcadlagfvpartialy}
        \left\|D_y f(t, x, y)\right\| & \leq c_4, \\
        \label{ftfyunifboundcadlagfvpartialyy}
        \left\|D_{yy} f(t, x, y)\right\| & \leq c_7,
    \end{align}
    in $(t, x, y)\in I\times \mathbb{R}^d\times \mathbb{R}^k$, for suitable constants $c_1, c_2, c_3, c_4, c_7 \geq 0$. Assume the semi-martingale noise $\{Y_t\}_{t\in I}$ is decomposed as $Y_t = Y_0 + F_t + Z_t$, where $\{F_t\}_{t\in I}$ is an FV process, $\{Z_{t \wedge T}\}_{t\geq 0}$ is a martingale, and $F_0 = Z_0 = 0,$ and assume they satisfy
    \begin{align}
        \label{Y0finiteexpectation}
        & \mathbb{E}\left[\|Y_0\|\right] < \infty, \\
        \label{expectFVFtfinite}
        & \mathbb{E}\left[V(\{F_t\}_{t\in I}; I)^{2}\right] < \infty, \\
        \label{Zmartingale}
        & \mathbb{E}\left[\|Z_T\|^2\right] < \infty.
    \end{align}   
    Then, the Euler scheme is strongly convergent of order 1, i.e.
    \begin{equation}
        \label{ordersemimartingaleunif}
        \max_{j=0, \ldots, N} \mathbb{E}\left[\left\| X_{t_j} - X_{t_j}^N \right\|\right] \leq c \Delta t_N, \qquad \forall N \in \mathbb{N},
    \end{equation}
    for a suitable constant $c\geq 0$ independent of $N.$
\end{theorem}

\begin{proof}
    We apply \cref{propbasicstrongestimate}. The troublesome term is the last one in \cref{expectedestimateglobalerrorintegral}. Using \cref{changeofvariablesformulasemimartingale}, the integral in this last term reads
    \begin{equation}
        \label{smlasttermerror}
        \begin{aligned}
            \int_0^{t_j} & \left( f(s, X_{\tau^N(s)}^N, Y_s) - f(\tau^N(s), X_{\tau^N(s)}^N, Y_{\tau^N(s)}) \right)\;\mathrm{d}s \\
            & \quad = \int_0^{t_j} \int_{\tau^N(s)}^s \partial_\xi f(\xi, X_{\tau^N(s)}^N, Y_{\xi^-})\;\mathrm{d}\xi \;\mathrm{d}s \\
            & \qquad + \int_0^{t_j} \int_{\tau^N(s)^+}^s D_y f(\xi, X_{\tau^N(s)}^N, Y_{\xi^-}) \;\mathrm{d}Y_\xi \;\mathrm{d}s \\
            & \qquad + \int_0^{t_j} \sum_{\tau^N(s) < \xi \leq s} \left(f(\xi, X_{\tau^N(s)}^N, Y_\xi) - f(\xi, X_{\tau^N(s)}^N, Y_{\xi^{-}}) \right. \\
            & \hspace{2.4in} \left. - D_y f(\xi, X_{\tau^N(s)}^N, Y_{\xi^-})\Delta Y_\xi\right) \;\mathrm{d}s \\ 
            & \qquad + \frac{1}{2} \int_0^{t_j} \int_{\tau^N(s)}^s D_{yy}f(\xi, X_{\tau^N(s)}^N, Y_{\xi^-})\;\mathrm{d}[Y, Y]_\xi^c\;\mathrm{d}s.
        \end{aligned}
    \end{equation}
    
    The bound on the expectation of the first term in \cref{smlasttermerror} is obtained using \cref{ftfyunifboundcadlagfvpartialt}, so that
    \begin{equation}
        \label{strongorderterm1}
        \begin{aligned}
            \mathbb{E}\bigg[\bigg\| \int_0^{t_j} \int_{\tau^N(s)}^s & \partial_\xi f(\xi, X_{\tau^N(s)}^N, Y_{\xi^-})\;\mathrm{d}\xi\;\mathrm{d}s\bigg\|\bigg] \\
            & \leq \mathbb{E}\left[\int_0^{t_j} \int_{\tau^N(s)}^s \left(c_1 + c_2 \|X_{\tau^N(s)}^N\| + c_3\|Y_{\xi^-}\|\right)\;\mathrm{d}\xi\;\mathrm{d}s\right] \\
            & \leq t_j \Delta t_N \left(c_1 + c_2 \sup_{0\leq s \leq t_j} \mathbb{E}\left[\|X_{\tau^N(s)}^N\|\right] + c_3\sup_{0 \leq s \leq t_j}\mathbb{E}\left[\|Y_s\|\right]\right).
        \end{aligned}
    \end{equation}

    For the second term in \cref{smlasttermerror}, we use the decomposition $Y_t = Y_0 + F_t + Z_t.$ With that, the second term can be written as
    \begin{equation}
        \begin{aligned}
            \label{strongorderterm2}
            \int_0^{t_j} \int_{\tau^N(s)^+}^s D_y f(\xi, X_{\tau^N(s)}^N, Y_{\xi^-}) \;\mathrm{d}Y_\xi \;\mathrm{d}s & = \int_0^{t_j} \int_{\tau^N(s)^+}^s D_y f(\xi, X_{\tau^N(s)}^N, Y_{\xi^-}) \;\mathrm{d}F_\xi \;\mathrm{d}s \\  
            & + \int_0^{t_j} \int_{\tau^N(s)^+}^s D_y f(\xi, X_{\tau^N(s)}^N, Y_{\xi^-}) \;\mathrm{d}Z_\xi \;\mathrm{d}s,
        \end{aligned}
    \end{equation}
    where the first inner integral is a Lebesgue-Stieltjes integral and the second inner integral is a stochastic integral with respect to a local martingale. Since $\{F_t\}_{t\in I}$ is of finite variation, the bound on the first integral on the right hand side of \cref{strongorderterm2} is obtained using \cref{ftfyunifboundcadlagfvpartialy} and switching the order of integration, so that
    \begin{align*}
        & \left\|\int_0^{t_j} \int_{\tau^N(s)^+}^s D_y f(\xi, X_{\tau^N(s)}^N, Y_{\xi^-}) \;\mathrm{d}F_\xi\;\mathrm{d}s\right\| \leq \int_0^{t_j} \int_{\tau^N(s)}^s c_4\;\|\mathrm{d}F_\xi\|\;\mathrm{d}s \\
        & \qquad \leq c_4\int_0^{t_j} \int_{\xi}^{\tau^N(\xi) + \Delta t_N} \;\mathrm{d}s\;\|\mathrm{d}F_\xi\|  \leq c_4\Delta t_N\int_0^{t_j} \;\|\mathrm{d}F_\xi\| \leq c_4\Delta t_N V(\{F_t\}_{t\in I}; [0, t_j]).
    \end{align*}
    Taking the expectation and using \cref{expectFVFtfinite}, we obtain
    \begin{equation}
        \label{strongorderterm2a}
        \mathbb{E}\left[\left\|\int_0^{t_j} \int_{\tau^N(s)^+}^s D_y f(\xi, X_{\tau^N(s)}^N, Y_{\xi^-}) \;\mathrm{d}F_\xi\;\mathrm{d}s\right\|\right] \leq c_4 \Delta t_N \mathbb{E}\left[ V(\{F_t\}_{t\in I}; [0, t_j]) \right].
    \end{equation}
    
    For the second integral on the right hand side of \cref{strongorderterm2}, we once again switch the order of integration, then use the Lyapunov inequality, the It\^o isometry \cref{eqitoisometrylocalmartingale} for local martingales, and then again the bound \cref{ftfyunifboundcadlagfvpartialy}, so that
    \begin{align*}
        \mathbb{E}\bigg[\bigg\|\int_0^{t_j} \int_{\tau^N(s)^+}^s & D_y f(\xi, X_{\tau^N(s)}^N, Y_{\xi^-}) \;\mathrm{d}Z_\xi\;\mathrm{d}s\bigg\|\bigg] \\
        & = \mathbb{E}\left[\left\|\int_0^{t_j} \int_{\xi}^{\tau^N(\xi) + \Delta t_N} D_y f(\xi, X_{\tau^N(s)}^N, Y_{\xi^-}) \;\mathrm{d}s \;\mathrm{d}Z_\xi\right\|\right] \\
        & \leq \mathbb{E}\left[\left\|\int_0^{t_j} \int_{\xi}^{\tau^N(\xi) + \Delta t_N} D_y f(\xi, X_{\tau^N(\xi)}^N, Y_{\xi^-}) \;\mathrm{d}s \;\mathrm{d}Z_\xi\right\|^2\right]^{1/2} \\
        & \leq c_4\Delta t_N\mathbb{E}\left[\sum_{l_1,l_2=1}^k\int_0^{t_j} \left|\mathrm{d}[Z^{(l_1)}, Z^{(l_2)}]_\xi\right|\right]^{1/2},
    \end{align*}
    where we used that, for any $0 \leq \xi \leq T,$ and for $\xi \leq s \leq \tau^N(\xi) + \Delta t_N,$ we have $\tau^N(s) = \tau^N(\xi),$ so we can write $X_{\tau^N(s)}^N = X_{\tau^N(\xi)}^N.$ Now, using that each quadratic variation and covariation is of finite variation, we obtain
    \begin{multline}
        \label{strongorderterm2b}
        \mathbb{E}\bigg[\bigg\|\int_0^{t_j} \int_{\tau^N(s)^+}^s D_y f(\xi, X_{\tau^N(s)}^N, Y_{\xi^-}) \;\mathrm{d}Z_\xi\;\mathrm{d}s\bigg\|\bigg] \\ 
        \leq c_4\Delta t_N\mathbb{E}\left[\sum_{l_1,l_2=1}^k V(\{[Z^{(l_1)}, Z^{(l_2)}]_t\}_{t\in I}; [0, t_j])\right]^{1/2}.
    \end{multline}
    
    Now we move to the third term in \cref{smlasttermerror}. We use the bound \cref{ftfyunifboundcadlagfvpartialyy} on the second derivative of $D_{yy}f$ to write
    \begin{equation*}
        \left\| f(\xi, X_{\tau^N(s)}^N, Y_\xi) - f(\xi, X_{\tau^N(s)}^N, Y_{\xi^{-}}) - D_y f(\xi, X_{\tau^N(s)}^N, Y_{\xi^-}) \Delta Y_\xi \right\| \leq \frac{c_7}{2}\|\Delta Y_\xi\|^2.
    \end{equation*}
    Recall $\tau^N(s) = \tau^N(\xi)$ for $\xi \leq s < \tau^N(\xi) + \Delta t_N$ and $\tau^N(\xi) + \Delta t_N- \xi \leq \Delta t_N$. Thus,
    \begin{align*}
        \Bigg\|\int_0^{t_j} \sum_{\tau^N(s) < \xi \leq s} \bigg(f(\xi, X_{\tau^N(s)}^N, Y_\xi) - & f(\xi, X_{\tau^N(s)}^N, Y_{\xi^{-}}) - D_y f(\xi, X_{\tau^N(s)}^N, Y_{\xi^-})\Delta Y_\xi\bigg) \;\mathrm{d}s\Bigg\| \\
        & \qquad \leq \frac{c_7}{2}\int_0^{t_j} \sum_{\tau^N(s) < \xi \leq s} \left\| \Delta Y_\xi \right\|^2\;\mathrm{d}s \\
        & \qquad \leq \frac{c_7}{2}\sum_{0 < \xi \leq t_j}\int_\xi^{\tau^N(\xi) + \Delta t_N} \left\| \Delta Y_\xi \right\|^2\;\mathrm{d}s \\
        & \qquad \leq \frac{c_7}{2} \Delta t_N  \sum_{0 < \xi \leq t_j}\left\|\Delta Y_\xi \right\|^2.
    \end{align*}
    Taking the expectation,
    \begin{multline}
        \label{strongorderterm3}
        \mathbb{E}\Bigg[\Bigg\|\int_0^{t_j} \sum_{\tau^N(s) < \xi \leq s} \bigg(f(\xi, X_{\tau^N(s)}^N, Y_\xi) - f(\xi, X_{\tau^N(s)}^N, Y_{\xi^{-}}) \\
        - D_y f(\xi, X_{\tau^N(s)}^N, Y_{\xi^-})\Delta Y_\xi\bigg) \;\mathrm{d}s\Bigg\|\Bigg]
        \leq \frac{c_7}{2} \Delta t_N \mathbb{E}\left[\left\|\sum_{0 < \xi \leq t_j}\left\| \Delta Y_\xi \right\|^2\right\|\right].
    \end{multline}

    The final term in \cref{smlasttermerror} involves a Lebesgue-Stieltjes integral and is handled in a similar way, except that instead of the It\^o isometry we just bound the Euclidian norm by the sum norm and obtain
    \begin{align*}
        \bigg\|\frac{1}{2} \int_0^{t_j} & \int_{\tau^N(s)}^s D_{yy}f(\xi, X_{\tau^N(s)}^N, Y_{\xi^-})\;\mathrm{d}[Y, Y]_\xi^c\;\mathrm{d}s \bigg\| \\
        & = \bigg\|\frac{1}{2} \int_0^{t_j} \int_{\xi}^{\tau^N(\xi) + \Delta t_N} D_{yy}f(\xi, X_{\tau^N(s)}^N, Y_{\xi^-})\;\mathrm{d}s \;\mathrm{d}[Y, Y]_\xi^c \bigg\| \\
        & \leq \frac{1}{2} \sum_{i=1}^n \left|\int_0^{t_j} \int_{\xi}^{\tau^N(\xi) + \Delta t_N} \left(D_{yy}f(\xi, X_{\tau^N(s)}^N, Y_{\xi^-})\;\mathrm{d}s \;\mathrm{d}[Y, Y]_\xi^c\right)^{(i)}\right| \\
        & \leq \frac{1}{2} \sum_{i=1}^n \left|\int_0^{t_j} \int_{\xi}^{\tau^N(\xi) + \Delta t_N} \sum_{l_1, l_2 = 1, \ldots, d} \partial_{y_{l_1}}\partial_{y_{l_2}} f_i(s, Y_{s^-})\;\mathrm{d}s\;\mathrm{d}[Y^{(l_1)}, Y^{(l_2)}]_\xi^c\right| \\
        & \leq \frac{1}{2} \sum_{i=1}^n \sum_{l_1, l_2 = 1, \ldots, d} \int_0^{t_j} \left(\int_{\xi}^{\tau^N(\xi) + \Delta t_N} \left|\partial_{y_{l_1}}\partial_{y_{l_2}} f_i(s, Y_{s^-})\right|\;\mathrm{d}s\right) \left|\;\mathrm{d}[Y^{(l_1)}, Y^{(l_2)}]_\xi^c\right| \\
        & \leq \frac{1}{2} \sum_{i=1}^n \sum_{l_1, l_2 = 1, \ldots, d} \int_0^{t_j} \left(\int_{\xi}^{\tau^N(\xi) + \Delta t_N} c_7 \;\mathrm{d}s\right)\left|\;\mathrm{d}[Y^{(l_1)}, Y^{(l_2)}]_\xi^c\right| \\
        & \leq \frac{c_7}{2}\Delta t_N \sum_{i=1}^n \sum_{l_1, l_2 = 1, \ldots, d} V(\{[Y^{(l_1)}, Y^{(l_2)}]_t\}_{t\in I}; [0, t_j])
    \end{align*}
    Taking the expectation,
    \begin{multline}
        \label{strongorderterm4}
        \mathbb{E}\left[\bigg\|\frac{1}{2} \int_0^{t_j} \int_{\tau^N(s)}^s D_{yy}f(\xi, X_{\tau^N(s)}^N, Y_{\xi^-})\;\mathrm{d}[Y, Y]_\xi^c\;\mathrm{d}s \bigg\|\right] \\
        \leq \frac{c_7}{2}\Delta t_N \sum_{i=1}^n \sum_{l_1, l_2 = 1, \ldots, d}\mathbb{E}\left[V(\{[Y^{(l_1)}, Y^{(l_2)}]_t\}_{t\in I}; [0, t_j])\right].
    \end{multline}

    Putting the bounds \cref{strongorderterm1}, \cref{strongorderterm2}, \cref{strongorderterm2a}, \cref{strongorderterm2b}, \cref{strongorderterm3}, and \cref{strongorderterm4} together, we obtain, from \cref{smlasttermerror}, that
    \begin{equation}
        \label{lasttermestimateforstrongconv}
        \mathbb{E}\bigg[\bigg\| \int_0^{t_j} \left( f(s, X_{\tau^N(s)}^N, Y_s) - f(\tau^N(s), X_{\tau^N(s)}^N, Y_{\tau^N(s)}) \right)\;\mathrm{d}s \bigg\|\bigg] \leq C_T \Delta t_N,
    \end{equation}
    where
    \begin{multline*}
            C_T = c_1 T + c_2T\sup_{0\leq s \leq T} \mathbb{E}\left[\|X_{\tau^N(s)}^N\|\right] + c_3T\sup_{0 \leq s \leq T}\mathbb{E}\left[\|Y_s\|\right] \\
            + c_4\mathbb{E}\left[ V(\{F_t\}_{t\in I}; I) \right] + c_4\mathbb{E}\left[\sum_{l_1,l_2=1}^k V(\{[Z^{(l_1)}, Z^{(l_2)}]_t\}_{t\in I}; I)\right]^{1/2} \\
            + \frac{c_7}{2} \mathbb{E}\left[\sum_{0 < s \leq T}\left\| \Delta Y_s \right\|^2\right] + \frac{c_7}{2} \sum_{i=1}^n \sum_{l_1, l_2 = 1, \ldots, d} V(\{[Y^{(l_1)}, Y^{(l_2)}]_t\}_{t\in I}; [0, t_j]).
    \end{multline*}
    Thanks to \cref{EX0strongbound}, \cref{CMmeanbound}, \cref{LtLXbound}, \cref{XNtmeanboundLXMt}, \cref{Y0finiteexpectation}, \cref{expectFVFtfinite}, \cref{Zmartingale}, \cref{expectationYtfinite}, \cref{expectationZtsquaredfinite}, \cref{eqquadraticvarcovarfiniteexpectation}, \cref{expectationofquadraticvariation}, and \cref{expectationsquaredjumps}, the constant $C_T$ above is finite and independent of $N.$ Plugging \cref{lasttermestimateforstrongconv} into \cref{expectedestimateglobalerrorintegral} and applying \cref{estimatesecondterminglobalstrongerrorintegral} and \cref{EX0X0Nunif} yield the desired strong order 1 convergence \cref{ordersemimartingaleunif}.
\end{proof}

\begin{remark}
    If the initial condition $X_0$ is bounded almost surely, then \cref{XtboundLXMt} shows that $X_t$ is also bounded almost surely, uniformly on a bounded interval $t\in [0, T].$ In this case, the conditions on the derivative of $f=f(t, x, y)$ can be relaxed to have an arbitrary growth on $x$. Similarly, if the noise has bounded support almost surely, uniformly on $t\in [0, T],$ then the derivatives are allowed to grow with $y$. If the initial condition has no bounded support but higher moments are bounded, then appropriate power growths on $x$ are allowed. Even the uniform Lipschitz condition \cref{LtLXbound} may be relaxed, provided higher-order moments of $X_0$ are finite.
\end{remark}

\section{Numerical example}
\label{secnumex}

We have briefly mentioned, in the Introduction, several numerical examples confirming the strong order one convergence and therefore showing this rate is optimal. Here, we describe how the order of convergence is estimated numerically and then we give more details on one of the examples, namely the non-homogeneous linear example with different types of noises. The remaining examples can be found in the github repository \cite{RODEConvEM2023}.

\subsection{Estimating the order of convergence}
\label{secconvergence}

For estimating the order of convergence, we use the Monte Carlo method, computing a number of numerical approximations $\{X_{t_j}^{N_i}(\omega_m)\}_{j=0, \ldots, N_i}$, of sample path solutions $\{X_t(\omega_m)\}_{t\in I}$, for samples $\omega_m,$ with $m = 1, \ldots, M$, and mesh resolutions $N_i,$ with $i=1, \ldots, i_{\max},$ and taking the maximum in time of the average of their absolute differences at the mesh points:
\begin{equation}
    \epsilon^{N_i} = \max_{j=0, \ldots, N_i} \frac{1}{M}\sum_{m=1}^M \left|X_{t_j}(\omega_m) - X_{t_j}^{N_i}(\omega_m)\right|.
\end{equation}
Then we fit the errors $\epsilon^{N_i}$ to the power law $C\Delta t_{N_i}^p$, in order to find $p$, along with a 95\% confidence interval. The fitting process and the estimation of the confidence intervals follow classical methods in Statistics \cite{DeGrootSchervish2018,HastieTibshiraniFriedman2009,JamesWittenHastieTibshirani2021,JohnsonWichern2007}. Reporting confidence intervals in the context of approximations of stochastic differential equations, in particular, can also be seen in \cite{KloedenPlatenSchurz2012}.

Here are the main parameters for the error estimate:
\begin{enumerate}
    \item The number $M\in\mathbb{N}$ of samples for the Monte Carlo estimate of the strong error.
    \item The time interval $[0, T], $ $T > 0,$ for the initial-value problem.
    \item The distribution law for the random initial condition $X_0$.
    \item A series of time steps $\Delta t_{N_i} = T/N_i$, with $N_i=2^{n_i},$ for some $n_i\in\mathbb{N},$ $i=1, \ldots, i_{\max},$ so that finer meshes are refinements of coarser meshes.
    \item A number $N_{\mathrm{tgt}}=2^{n_\mathrm{tgt}}$ for a finer resolution to compute a target solution path, typically $N_{\mathrm{tgt}} = \max_i\{N_i^2\}$, unless an exact pathwise solution is available, in which case a coarser mesh of the order of $\max_i\{N_i\}$ can be used.
\end{enumerate}

With this setup, we proceed as follows.
\begin{enumerate}
    \item For each sample $\omega_m,$ $m=1, \ldots, M$, we first generate a discretization of a sample path of the noise, $\{Y_{t_j}\}_{j=0, \ldots, N_{\mathrm{tgt}}},$ on the finest grid $t_j = j \Delta t_{N_{\textrm{tgt}}}$, $j = 0, \ldots, N_{\mathrm{tgt}},$ using an exact distribution for the noise.
    \item Next, we use the values of the noise at the target time mesh to generate the target solution $\{X_{t_j}\}_{j=0, \ldots, N_{\mathrm{tgt}}}$, still on the fine mesh. This is constructed either using the Euler approximation itself, keeping in mind that the mesh is sufficiently fine, or an exact distribution of the solution, when available.
    \item Then, for each time resolution $N_i,$ we compute the Euler approximation using the computed noise values at the corresponding coarser mesh $t_j = j\Delta t_{N_i},$ $j=0, \ldots, N_i.$
    \item We then compare each approximation $\{X_{t_j}^{N_i}\}_{j=0, \ldots, N_i}$ to the values of the target path on that coarse mesh and update the strong error
    \[
        \epsilon_{t_j}^{N_i} = \frac{1}{M}\sum_{m=1}^M \left|X_{t_j}(\omega_m) - X_{t_j}^{N_i}(\omega_m)\right|,
    \]
    at each mesh point.
    \item At the end of all the simulations, we take the maximum in time, on each corresponding coarse mesh, to obtain the error for each mesh,
    \[
        \epsilon^{N_i} = \max_{j=0, \ldots, N_i} \epsilon_{t_j}^{N_i}.
    \]
    \item We fit $(\Delta t_{N_i}, \epsilon^{N_i})$ to the power law $C\Delta t_{N_i}^p$, via linear least-square regression in log scale, so that $\ln \epsilon^{N_i} \approx \ln C + p \ln \Delta t_{N_i}$, for suitable $C$ and $p$, with $p$ giving the order of convergence. This amounts to solving the normal equation $(A^\tr A)\mathbf{v} = A^\tr\ln(\boldsymbol{\epsilon})$, where $\mathbf{v}$ is the vector $\mathbf{v} = (\ln(C), p)$, $A$ is the Vandermonde matrix associated with the logarithm of the mesh steps $(\Delta t_{N_i})_i$, and $\ln(\boldsymbol{\epsilon})$ is the vector $\ln(\boldsymbol{\epsilon}) = (\ln(\epsilon^{N_i}))_i$.
    \item We also compute the standard error of the Monte-Carlo sample at each time step,
    \[
        s_{t_j}^{N_i} = \frac{\sigma_{t_j}^{N_i}}{\sqrt{M}},
    \]
    where $\sigma_{t_j}^{N_i}$ is the sample standard deviation given by
    \[
        \sigma_{t_j}^{N_i} = \sqrt{\frac{1}{M-1}\sum_{m=1}^M \left(\left|X_{t_j}(\omega_m) - X_{t_j}^{N_i}(\omega_m) \right|- \epsilon_{t_j}^{N_i}\right)^2},
    \]
    and compute the 95\% confidence interval $[\epsilon_{\min}^{{N_i}}, \epsilon_{\max}^{{N_i}}]$ for the strong error at each mesh resolution $N_i$ with
    \[
        \epsilon_{\min}^{N_i} = \max_{j=0, \ldots, N_i} (\epsilon_{t_j}^{N_i} - {1.96}\sigma_{t_j}^{N_i}), \quad \epsilon_{\max}^{N_i} = \max_{j=0, \ldots, N_i} (\epsilon_{t_j}^{N_i} + {1.96}\sigma_{t_j}^{N_i}).
    \]
    \item Next we estimate the 95\% confidence interval for the joint distribution $\boldsymbol{\epsilon} = (\epsilon^{N_i})_{i=1, \ldots, i_{\max}}$ using the Bonferroni inequality \cite{JohnsonWichern2007}. More precisely, we find the $z$-score value for the $\alpha = 1 - 0.05/i_{\max}$ confidence level for a univariate normal distribution and use that to obtain the $\alpha$ confidence intervals $[\epsilon_{\alpha, \min}^{N_i}, \epsilon_{\alpha, \max}^{N_i}]$ for the marginals $\epsilon^{N_i}$ with 
    \[
        \epsilon_{\alpha, \min}^{N_i} = \max_{j=0, \ldots, N_i} (\epsilon_{t_j}^{N_i} - z\sigma_{t_j}^{N_i}), \quad \epsilon_{\alpha,\max}^{N_i} = \max_{j=0, \ldots, N_i} (\epsilon_{t_j}^{N_i} + z\sigma_{t_j}^{N_i}).
    \]
    Then, the probability of the region $\Pi_{i=1}^{i_{\max}} [\epsilon_{\alpha, \min}^{N_i}, \epsilon_{\alpha, \max}^{N_i}]$ is bounded from below by $1 - i_{\max}(1 - \alpha) = 1 - 0.05 = 0.95,$ which means this is the confidence region of confidence level (at least) 95\%, for $\boldsymbol{\epsilon}$.
    
    \item Finally, from the normal equations and the confidence region for $\boldsymbol{\epsilon}$, we compute the (at least) 95\% confidence interval $[p_{\min}, p_{\max}]$ for $p$ via the transformation method. We do so by computing the minimum and maximum values of $p$ in the image of the confidence region $\Pi_{i=1}^{i_{\max}} [\epsilon_{\alpha, \min}^{N_i}, \epsilon_{\alpha, \max}^{N_i}]$ by the map $\mathbf{e} \mapsto (\ln(C), p) = (A^\tr A)^{-1}A^\tr \ln(\mathbf{e})$. The map $\mathbf{e} \mapsto \ln(\mathbf{e})$ takes the hyperrectangle $\Pi_{i=1}^{i_{\max}} [\epsilon_{\alpha, \min}^{N_i}, \epsilon_{\alpha, \max}^{N_i}]$ into the hyperrectangle $\Pi_{i=1}^{i_{\max}} [\ln(\epsilon_{\alpha, \min}^{N_i}), \ln(\epsilon_{\alpha, \max}^{N_i})].$ This hyperrectangle is taken onto a polygon in the plane $(C, p)$ by the linear map ${(A^\tr A)^{-1}A^\tr },$ which is then projected onto an interval $[p_{\min}, p_{\max}]$ in the $p$ axis. Since the pseudo-inverse ${(A^\tr A)^{-1}A^\tr }$ and the projection are not bijective, the corresponding confidence level of this interval is at least 95\%, but probably higher. In order to find $p_{\min}$ and $p_{\max},$ we just need to look at the vertices of the polygon in $(C, p)$, which are the set of points obtained by all combinations of the indices in $(\epsilon_{{\alpha,}\min}^{N_i})_i$ and $(\epsilon_{{\alpha,}\max}^{N_i})_i.$ We take the minimum of the $p$-coordinate of such points to be $p_{\min}$ and the maximum to be $p_{\max}.$
\end{enumerate}

\subsection{Non-homogeneous linear system of RODEs with different types of noises}
\label{subseclinearrodenumerics}

We examine a linear equation with noise present in both the homogeneous and nonhomogeneous terms, approached in two ways: first, as a series of scalar equations, each with a unique univariate noise in a family of different types of noises; and secondly, as a system of equations that combines all these noise types. For most of the noises, the current knowledge expects a lower order of strong convergence than the strong order 1 we prove here. The aim of this section is to illustrate this improvement for all such noises.

The equations take the form
\begin{equation}
    \label{supp:allnoisesRODEsystem}
    \begin{cases}
        \displaystyle \frac{\mathrm{d}{X}_t}{\mathrm{d} t} = - \left\|{Y}_t\right\|^2 {X}_t + {Y}_t, \qquad 0 \leq t \leq T, \\
        \left. {X}_t \right|_{t = 0} = {X}_0,
    \end{cases}
\end{equation}
where $\{X_t\}_{t\in I}$ is either a scalar or a vector-valued process, and $\{Y_t\}_{t\in I}$ is, correspondingly, either a scalar or a vector-valued noise process with the same dimension as ${X}_t$. In the case of a system, each coordinate of $\{Y_t\}_{t\in I}$ is a scalar noise process independent of the noises in the other coordinates. The scalar noises used in the simulations are the following, in the order of coordinates of $Y_t$:
\begin{enumerate}
    \item A standard Wiener process (W);
    \item An Ornstein-Uhlenbeck process (OU) with drift $\nu = 0.3,$ diffusion $\sigma = 0.5,$ and initial condition $y_0 = 0.2$;
    \item A geometric Brownian motion process (gBm) with drift $\mu = 0.3,$ diffusion coefficient $\sigma = 0.5,$ and initial condition $y_0 = 0.2$;
    \item A non-autonomous homogeneous linear It\^o process (hlp) $\{H_t\}_{t\in I}$ given by $\mathrm{d}H_t = (\mu_1 + \mu_2\sin(\vartheta t))H_t\;\mathrm{d}t + \sigma\sin(\vartheta t)H_t\;\mathrm{d}W_t$ with $\mu_1 = 0.5,$ $\mu_2 = 0.3,$ $\sigma = 0.5,$ $\vartheta=3\pi,$ and initial condition $H_0 = 0.2;$
    \item A compound Poisson process (cP) with rate $\lambda = 5.0$ and jump law following an exponential distribution with scale $\theta = 0.5;$
    \item A Poisson step process (sP) with rate $\lambda = 5.0$ and step law following a uniform distribution within the unit interval;
    \item An exponentially decaying Hawkes process with initial rate $\lambda_0 = 3.0$, base rate $a = 2.0$, exponential decay rate $\delta = 3.0$, and jump law following an exponential distribution with scale $\theta = 0.5;$
    \item A transport process of the form $t \mapsto \sum_{i=1}^{6} \sin^{1/3}(\omega_i t)$, where the frequencies $\omega_i$ are independent random variables following a Gamma distribution with shape parameter $\alpha = 7.5$ and scale $\theta = 2.0;$
    \item A fractional Brownian motion (fBm) process with Hurst parameter $H=0.6$ and initial condition $y_0 = 0.2$.
\end{enumerate}

\begin{table}
    \begin{center}
        \begin{tabular}[htb]{|r|l|l|l|}
            \hline N & dt & error & std err \\
            \hline \hline
            64 & 0.0156 & 0.205 & 0.0093 \\
            128 & 0.00781 & 0.0993 & 0.00451 \\
            256 & 0.00391 & 0.049 & 0.00223 \\
            512 & 0.00195 & 0.0243 & 0.00111 \\
            \hline
        \end{tabular}
    \end{center}

    \bigskip

    \caption{Mesh points (N), time steps (dt), strong error (error), and standard error (std err) of the Euler method for $\mathrm{d}{X}_t/\mathrm{d}t = - \| {Y}_t\|^2 {X}_t + {Y}_t$ for each mesh resolution $N$, with initial condition ${X}_0 \sim \mathcal{N}({0}, \mathrm{I})$ and vector-valued noise $\{{Y}_t\}_t$ with all the implemented noises, on the time interval $I = [0.0, 1.0]$, based on $M = 400$ sample paths for each fixed time step, with the target solution calculated with $262144$ points. The order of strong convergence is estimated to be $p = 1.024$, with the 95\% confidence interval $[0.8926, 1.1556]$.}
    
        \label{supp:taballnoises}    
\end{table}

\begin{figure}[htb]
    \centerline{\includegraphics[width=0.8\textwidth]{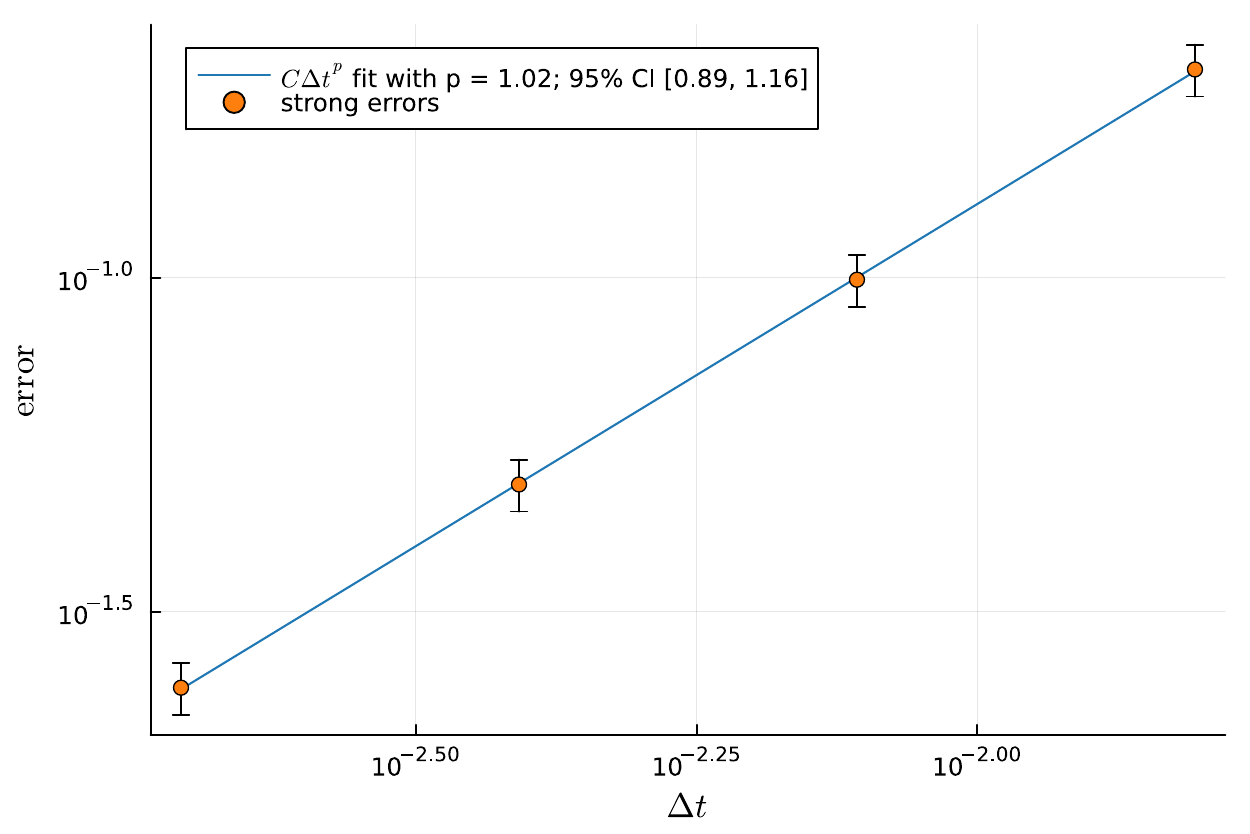}}
    \caption{Estimated order of convergence of the strong error of the Euler method for $\mathrm{d}{X}_t/\mathrm{d}t = - \left\|{Y}_t\right\|^2 {X}_t + {Y}_t,$ based on \cref{supp:taballnoises}.}
    \label{supp:figallnoises}
\end{figure}

\Cref{supp:taballnoises} shows the estimated strong error obtained from the Monte Carlo simulations, with the target solution obtained via Euler approximation on a much finer mesh. 

\begin{figure}[htb]
    \centerline{\includegraphics[width=0.8\textwidth]{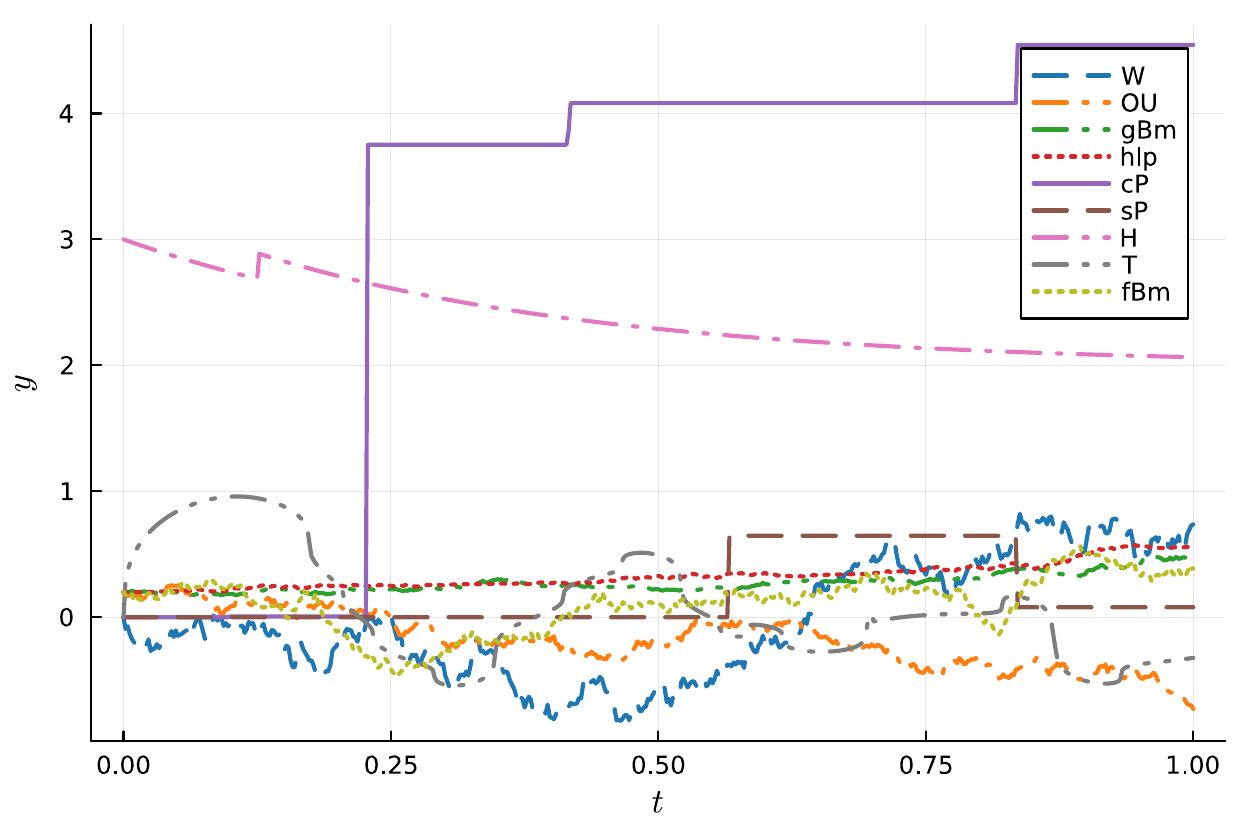}}
    \caption{Sample paths of all the noises used in the linear system \cref{supp:allnoisesRODEsystem}, mixing all different types of implemented noises.}
    \label{supp:figsamplepathsallnoises}
\end{figure}

\Cref{supp:figallnoises} illustrates the obtained order of convergence, while \cref{supp:figsamplepathsallnoises} illustrates some sample paths of all the noises used in this system.

\begin{figure}[htb]
    \centerline{\includegraphics[width=0.8\textwidth]{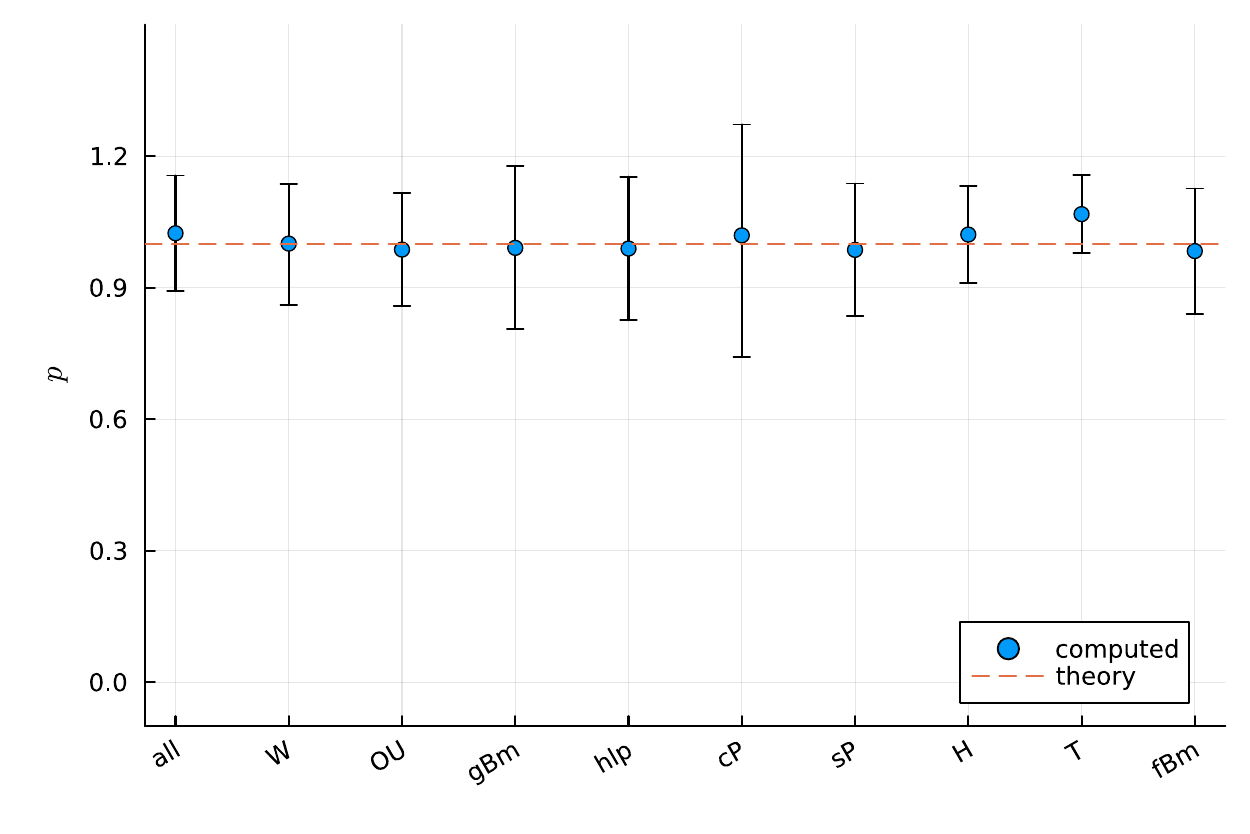}}
    \caption{Estimated order of convergence for the equation \cref{supp:allnoisesRODEsystem} as it varies among different types of noises, starting with all the noises combined in a system of equations (all) and continuing with the scalar cases, each with a different type of noise, according to \cref{supp:tabeachnoise}.}
    \label{supp:figordereachnoise}
\end{figure}

\Cref{supp:tabeachnoise} and \cref{supp:figordereachnoise} show the order of convergence as it varies among different types of noises, starting with all the noises combined in a system of equations and continuing with the scalar cases, each with a different type of noise, as described above.

\begin{table}
    \begin{center}
        \begin{tabular}[htb]{|l|l|l|l|}
            \hline noise & $p$ & $p_{\textrm{min}}$ & $p_{\textrm{max}}$ \\
            \hline \hline
            all noises combined (all) & 1.02409 & 0.892626 & 1.15557 \\
            Wiener (W) & 1.0005 & 0.860131 & 1.1361 \\
            OrnsteinUhlenbeck (OU) & 0.98683 & 0.858187 & 1.11547 \\
            GeometricBrownianMotion (gBm) & 0.990724 & 0.806347 & 1.17738 \\
            HomogeneousLinearIto (hlp) & 0.989188 & 0.825703 & 1.15296 \\
            CompoundPoisson (cP) & 1.01923 & 0.742455 & 1.2723 \\
            PoissonStep (sP) & 0.98615 & 0.835225 & 1.1373 \\
            ExponentialHawkes (H) & 1.02111 & 0.910741 & 1.13141 \\
            Transport (T) & 1.0677 & 0.97865 & 1.15692 \\
            FractionalBrownianMotion (fBm) & 0.983487 & 0.840285 & 1.12595 \\
            \hline
        \end{tabular}
    \end{center}

    \bigskip

    \caption{Type of noise, order $p$ of strong convergence, and 95\% confidence interval endpoints $p_{\textrm{min}}$ and $p_{\textrm{max}}$, for the linear equation \cref{supp:allnoisesRODEsystem} with various types of noise.}

    \label{supp:tabeachnoise}
\end{table}

The strong order 1 convergence is not a surprise in the case of the Wiener and Ornstein-Uhlenbeck process since the corresponding RODE can be turned into an SDE with an additive noise. In such a case, the Euler-Maruyama approximation for the noise part of the SDE is distributionally exact, and the Euler method for the RODE becomes equivalent to the Euler-Maruyama method for the SDE, and it is known that the Euler-Maruyama method for an SDE with additive noise is of strong order 1 \cite{HighamKloeden2021}. The geometric Brownian motion is also expected to yield strong order 1 convergence thanks to the results in \cite{WangCaoHanKloeden2021} and the fact that it is an It\^o diffusion process. For the remaining noises, however, previous works would show the order of convergence to be below the order 1 attained here.

Notice we chose the Hurst parameter of the fractional Brownian motion process to be between 1/2 and 1, so that our results for fBm indicate that the strong convergence should also be of order 1, just like for the other types of noises in $\{Y_t\}_{t\in I}$. Previous knowledge would expect a strong convergence of order $H$, with $1/2 < H < 1,$ instead.

\section*{Data Availability Statement}

The research code associated with this article is available at the open-access, github repository \url{https://github.com/rmsrosa/rode_convergence_euler}, and detailed documentation for the repository is available at the website \url{https://rmsrosa.github.io/rode_convergence_euler} \cite{RODEConvEM2023}.

\section*{Acknowledgments}

The authors thank Glauco Valle (IM/UFRJ), for helpful discussions about semi-martingales, Luan Lima Freitas (IM/UFRJ), for many helpful comments on the manuscript, and Ralph dos Santos Silva (IM/UFRJ) for helpful discussions about the estimates for the confidence intervals. We also like to thank the reviewers for a careful reading of the manuscript, with several comments and suggestions that helped improve the paper. The second author was partly supported by Laborat\'orio de Matem\'atica Aplicada, Instituto de Matem\'atica, Universidade Federal do Rio de Janeiro (LabMA/ IM/UFRJ), Coordena\c{c}\~ao de Aperfei\c{c}oamento de Pessoal de N\'{\i}vel Superior (CAPES), Brasil, grant 001, and Conselho Nacional de Desenvolvimento Cient\'{\i}fico e Tecnol\'ogico (CNPq), Bras\'{\i}lia, Brasil, grant no. 408751/2023-1.

\appendix

\section{Numerical setup and examples}
\label{secsupplementary}

In this supplementary material, we describe in details the remaining numerical examples illustrating the main work. In all the examples driven by semi-martingale noises, we numerically verify the strong order 1 convergence of the Euler method. In the example with fractional Brownian motion, we verify the strong convergence of order $\min\{H+1/2, 1\},$ where $H$ is the Hurst parameter. The order of convergence is displayed with a 95\% confidence interval that perfectly matches the theory developed in the paper. These estimates show that the theoretical estimate is sharp.

The code for all the examples is written in the Julia Programming Language \cite{Julia2017} and is available in the github repository \url{https://github.com/rmsrosa/rode_convergence_euler/} \cite{RODEConvEM2023}. The Julia language itself is available for download on the website \url{https://julialang.org/}.

One of the distinguishing powers of the Julia language is to work both as a high-level, interactive, scripting language (with its dynamically typed and just-in-time compilation system), and as a fast, high-performance language (with its aggressive compilation capabilities). It allows for both easy prototyping and the production of efficient end code. As such, it is growing quickly as a language suitable for scientific computing and computationally-demanding real-world applications.

Our code has been tested in the currently most recent Julia version 1.11 as well as in previous versions down to version 1.8. The code is set up in a way that any change to the repository triggers an automatic rerun of all the documentation code and all the unit tests in a github server, to make sure it is always up to date. The result of the github documentation run is displayed online, on the website \url{https://rmsrosa.github.io/rode_conv_em/}.

\subsection{Numerical setup}
\label{secnumericalsetup}

Julia has a feature-rich \texttt{DifferentialEquations.jl} ecosystem of packages for solving differential equations \cite{DifferentialEquations.jl-2017}, including random and stochastic differential equations, as well as delay equations, differential-algebraic equations, jump diffusions, partial differential equations, neural differential equations, and so on. It also has packages to seemlessly compose such equations in optimization problems, Bayesian parameter estimation, global sensitivity analysis, uncertainty quantification, and domain specific applications.

Although all the source code for \texttt{DifferentialEquations.jl} is publicly available, it involves a quite large ecosystem of packages, with an intricate interplay between them. Hence, for the numerical results presented here, we chose to implement our own routines, with a minimum set of methods necessary for the convergence estimates. This is done mostly for the sake of transparency, so that the reviewing process for checking the accuracy of the implementation becomes easier. All the source code for the numerical simulations presented below are in the github repository \cite{RODEConvEM2023}.

\subsubsection{Simulating a Random ODE}

For the simulation of an approximation of a pathwise solution of a random ODE
\begin{equation}
    \begin{cases}
      \displaystyle \frac{\mathrm{d}X_t}{\mathrm{d} t} = f(t, X_t, Y_t), & 0 \leq t \leq T, \\
      \left. X_t \right|_{t = 0} = X_0,
    \end{cases}
\end{equation}
we define the following variables
\begin{itemize}
    \item We fix the initial time, usually \texttt{t0=0.0}, and set the final time \texttt{tf} as the desired time $T$ for the time interval $[0, T].$
    \item The initial condition \texttt{x0} as a scalar or vector, depending on whether it is a scalar equation or a system of equations.
    \item A noise sample \texttt{yt} as an \texttt{AbstractVector}, indexed by a time index \texttt{j}, where each \texttt{yt[j]} represents the state of the noise sample at the $j$-th time in the time mesh. The time mesh is computed from the initial and final times and from the length \texttt{n} of \texttt{yt}, starting at time \texttt{t0}, at index \texttt{j=1}, and advancing at equal time steps \texttt{dt = ( tf - t0 ) / (n - 1)}, up to time \texttt{tf}, at index \texttt{j=n}. The value \texttt{yt[j]} is either a scalar or a vector, depending on the dimension of the noise.
    \item A function representing $f(t, X_t, Y_t)$ which can either be in the out-of-place form \texttt{f(t, x, y, p)}, in the case of a scalar equation, or in the in-place form \texttt{f!(dx, t, x, y, p)}, in the case of a system of equations, where \texttt{dx} is a vector to store the value of the numerical representation of $f(t, X_t, Y_t).$ In either case, \texttt{p} is used to store extra parameters in the equation.
    \item A set of parameters \texttt{p = params} to be passed on to either \texttt{f(t, x, y, p)} or \texttt{f!(dx, t, x, y, p)}, when solving the equation.
    \item An \texttt{AbstractVector}, denoted here as \texttt{xt}, with the same length as \texttt{yt} and with each \texttt{xt[j]} of the same type as \texttt{x0}, pre-allocated in memory to hold the solution sample path.
\end{itemize}
Then, one computes the Euler approximation via
\[
    \texttt{solve!(xt, t0, tf, x0, f, yt, params, ::RandomEuler)}.
\]
This fills up \texttt{xt} with the sample solution.

\subsubsection{Calculating the order of convergence}

For the implementation of the calculation of the order of convergence as described in \cref{secconvergence}, we also define the following variables.
\begin{itemize}
    \item A random seed \texttt{rng} as an \texttt{AbstractRNG}, from the \texttt{Random} standard Julia library, so the computed samples are reproducible. In the applications, we use the  Julia implementation of the Xoshiro256++ pseudorandom number generator \cite{BlackmanVigna2021}.
    \item The initial distribution law \texttt{x0law} as a univariate or multivariate distribution from the \texttt{Distributions.jl} package \cite{JSSv098i16}, depending on whether $X_t$ is a scalar or vector valued.
    \item A \texttt{noise} defining a univariate or multivariate process. There are several implemented noises, with the corresponding parameter options. The implemented univariate processes are 
    \begin{itemize}
        \item \texttt{WienerProcess(...)}
        \item \texttt{OrnsteinUhlenbeckProcess(...)}
        \item \texttt{HomogeneousLinearItoProcess(...)}
        \item \texttt{GeometricBrownianMotionProcess(...)}
        \item \texttt{CompoundPoissonProcess(...)}
        \item \texttt{PoissonStepProcess(...)}
        \item \texttt{ExponentialHawkesProcess(...)}
        \item \texttt{TransportProcess(...)}
        \item \texttt{FractionalBrownianMotionProcess(...)}.
    \end{itemize}
    Then there is the option to create a multivariate process
    \begin{itemize}
        \item \texttt{ProductProcess(noises...)}
    \end{itemize}
    with any combination of the univariate process as components. For a thorough explanation of the noises, see the github repository \cite{RODEConvEM2023}.
    \item A \texttt{target} method, which is either \texttt{RandomEuler()} to solve the equation via Euler approximation on a finer mesh or a custom method to yield an exact sample of the solution, when a closed-form solution for the exact process $\{X_t\}_t$ is available.
    \item The method \texttt{method = RandomEuler()} for solving the approximations at each resolution $N_i$ via Euler method. The method \texttt{RandomEuler()} is just an instance of a \texttt{RandomEuler} type, which indicates the solver to use the routine for computing the Euler approximation with the given parameters and stores cached vectors for the efficient computation in the case of systems of equations.
    \item A target resolution \texttt{ntg} representing $N_{\textrm{tgt}}.$
    \item A vector $\texttt{ns}$ containing the desired resolutions $(N_i)_{i=1, \ldots, I}$ for the coarsers meshes.
    \item An integer \texttt{m} representing the number $M$ of Monte Carlo simulations.
    \item When approximating a PDE, a vector \texttt{ks} of steps used to reduce the spatial resolution may be provided. It is a vector of the same length as \texttt{ns}, to accompany each time mesh resolution $N_i$. If the target approximation $X_{t_j}^{N_{\textrm{tgt}}}$ has dimension/spatial resolution \texttt{d = length(xt[i])}, then each $X_{t_j}^{N_i}$ has dimension/resolution \texttt{div(d-1,ks[i])+1}. If \texttt{ks} is not provided, it defaults to a vector of the same length as $\texttt{ns}$ and all components equal to \texttt{1} for no change in resolution.
    \item With all these variables set up, the ``convergence suite'' is defined by
    \begin{multline*}
        \texttt{suite = ConvergenceSuite(t0, tf, x0law, f, noise,} \\ 
        \texttt{params, target, \allowbreak method, \allowbreak ntgt, \allowbreak ns, \allowbreak m, ks)}.
    \end{multline*}
\end{itemize}

Then, we compute the convergence order via 
\[ \texttt{result = solve(rng, suite)}.
\]
The returned \texttt{result} is of a type \texttt{ConvergenceResult}, which stores all the necessary info about the errors at each resolution, their confidence intervals, the order of convergence and its confidence interval, and so on. The error table can be accessed via \texttt{generate\_error\_table(results)}. The order of convergence and the bounds for the 95\% confidence interval are stored in \texttt{result.p}, \texttt{result.pmin}, \texttt{result.pmax}. A plot of the convergence errors, the confidence intervals, and the fitted error estimate $C\Delta t^p$ can be visualized via \texttt{plot(results)}, which contains a plotting recipe based on the \texttt{Plots.jl} package \cite{PlotsJL2023}.

\subsection{Numerical examples}
\label{secnumericalexamples}

We now describe in more details the examples used to illustrate the order of convergence, except the non-homogeneous linear system, which is described in \cref{secnumex}. We start with a linear equation driven by a Wiener process, then we illustrate the $\min\{H + 1/2, 1\}$ order of convergence in the case of a fractional Brownian motion (fBm) noise with Hurst parameter $0 < H < 1.$ Next we move to a number of linear and nonlinear models, with different semi-martingale noises, as mentioned in the \cref{secintro}.

\subsubsection{Homogeneous linear equation with Wiener noise}
\label{seclinearhomogeneousrode}

We start with one of the simplest Random ODEs, that of a linear homogenous equation with a Wiener process as the coefficient:
\begin{equation}
    \label{linearhomogeneousrode}
    \begin{cases}
        \displaystyle \frac{\mathrm{d}X_t}{\mathrm{d} t} = W_t X_t, \qquad 0 \leq t \leq T, \\
        \left. X_t \right|_{t = 0} = X_0.
      \end{cases}
\end{equation}

Since the noise is simply a Wiener process, the corresponding RODE can be turned into an SDE with an additive noise. In this case, the Euler-Maruyama approximation for the noise part of the SDE is distributionally exact and the Euler method for the RODE becomes equivalent to the Euler-Maruyama method for the SDE. Moreover, it is known that the Euler-Maruyama method for an SDE with additive noise is of strong order 1 \cite{HighamKloeden2021}. Nevertheless, we illustrate the strong convergence directly for the Euler method for this RODE equation, for the sake of completeness.

Equation \cref{linearhomogeneousrode} has the explicit solution
\begin{equation}
    \label{Xtlinearhomogeneousrodesolution}
    X_t = e^{\int_0^t W_s \;\mathrm{d}s}X_0.
\end{equation}

When we compute an approximate solution via Euler's method, however, we only draw the realizations $\{W_{t_i}\}_{i=0}^N$ of a sample path, on the mesh points. We cannot compute the exact integral $\int_0^{t_j} W_s\;\mathrm{d}s$ just from these values, and, in fact, an exact solution is not uniquely defined from them. We can, however, find its exact distribution and use that to draw feasible exact solutions and use them to estimate the error.

The exact distribution of $\int_0^\tau W_s\;\mathrm{d}s$ given the step $\Delta W = W_\tau - W_0 = W_\tau$ is computed in \cite[Section 14.2]{HanKloeden2017} as
\begin{equation}
    \label{exactintWt}
    \int_0^{\tau} W_s\;\mathrm{d}s \sim \frac{\tau}{2}\Delta W + \sqrt{\frac{\tau^3}{12}}\mathcal{N}(0, 1) = \frac{\tau}{2}\Delta W + \mathcal{N}\left(0, \frac{\tau^3}{12}\right).
\end{equation}
As for the distribution of the integral over a mesh interval $[t_i, t_{i+1}]$ when given the endpoints $W_{t_i}$ and $W_{t_{i+1}},$ we use that $s \mapsto W_{t_i+s} - W_{t_i}$ is a standard Wiener process to find, from \cref{exactintWt}, that
\begin{align*}
    \int_{t_i}^{t_{i+1}} W_s\;\mathrm{d}s & = W_{t_i}(t_{i+1} - t_i) + \int_{t_i}^{t_{i+1}} (W_s - W_{t_i})\;\mathrm{d}s \\
    & = W_{t_i}(t_{i+1} - t_i) + \int_{0}^{t_{i+1} - t_i} (W_{t_i+s} - W_{t_i})\;\mathrm{d}s \\
    & = W_{t_i}(t_{i+1} - t_i) + \frac{(t_{i+1} - t_i)}{2}(W_{t_{i+1}}-W_{t_{i}}) + Z_i \\
    & = \frac{(W_{t_{i+1}}+W_{t_{i}})}{2}(t_{i+1} - t_i) + Z_i,
\end{align*}
where
\begin{equation}
    \label{linearhomogeneousZidistribution}
    Z_i \sim \mathcal{N}\left(0, \frac{(t_{i+1}- t_i)^3}{12}\right) = \frac{\sqrt{(t_{i+1} - t_i)^3}}{\sqrt{12}}\mathcal{N}(0, 1).
\end{equation}

Thus, breaking down the sum over the mesh intervals,
\begin{equation}
    \label{exactintW}
    \int_0^{t_j} W_s\;\mathrm{d}s = \sum_{i = 0}^{j-1} \int_{t_i}^{t_{i+1}} W_s\;\mathrm{d}s = \sum_{i=0}^{j-1} \left( \frac{(W_{t_{i+1}}+W_{t_{i}})}{2}(t_{i+1} - t_i) + Z_i\right),
\end{equation}
where $Z_0, \ldots, Z_{N-1}$ are independent random variables with distribution \cref{linearhomogeneousZidistribution}.

Notice the first term in \cref{exactintW} is the trapezoidal rule, with the second one being a correction term.

\begin{figure}
    \centerline{\includegraphics[width=0.8\textwidth]{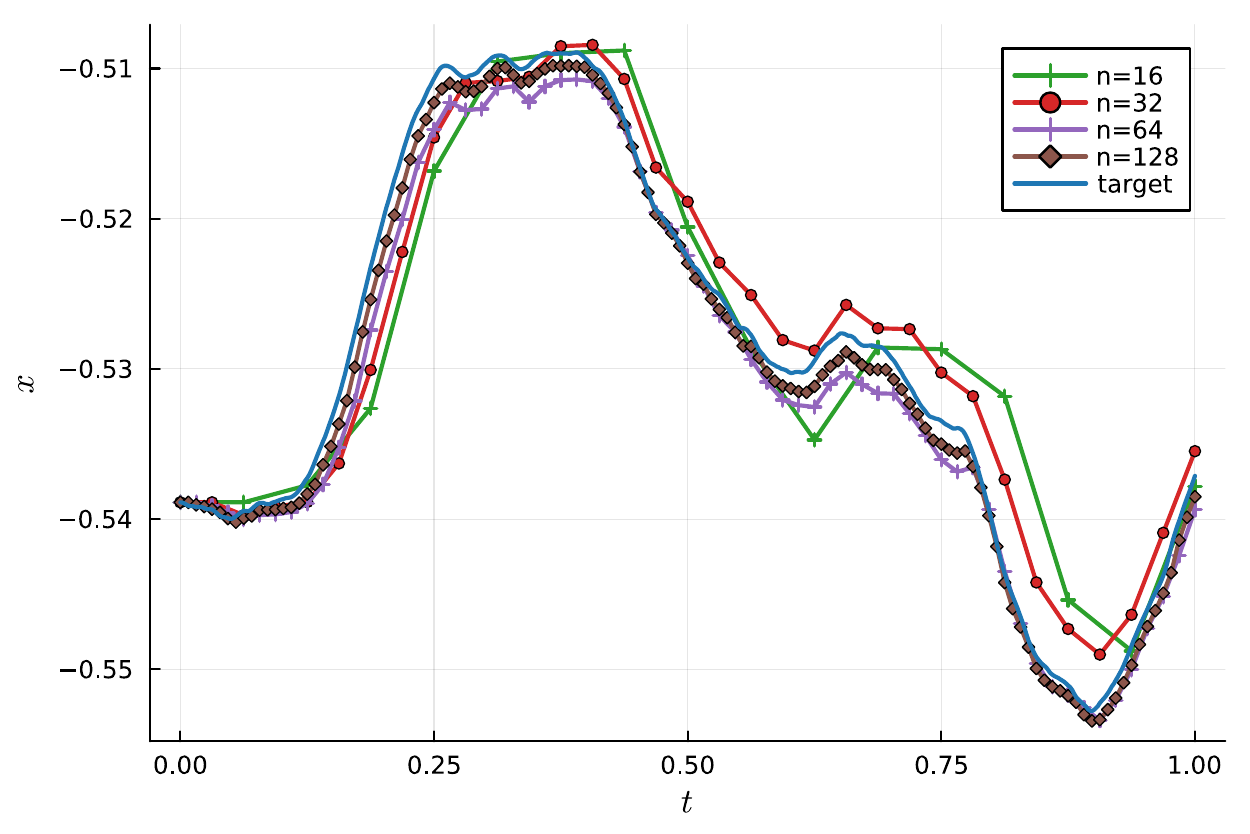}}
    \caption{Euler approximation of a sample solution of $\mathrm{d}X_t/\mathrm{d}t = W_t X_t$ with $X_0 \sim \mathcal{N}(0, 1)$, with a few mesh resolutions and with the target resolution.}
    \label{samplepathslinearhomogeneousrode}
\end{figure}

For a normal variable $N \sim \mathcal{N}(\mu, \sigma)$, the expectation of the random variable $e^N$ is $\mathbb{E}[e^N] = e^{\mu + \sigma^2/2}$. Hence,
\begin{equation}
    \label{Xtlinearhomogeneousrodeestimateexponentialofnomral}
    \mathbb{E}[e^{Z_i}] = e^{((t_{i+1}- t_i)^3)/24}.
\end{equation}
This is the contribution of this random variable to the mean of the exact solution. But we actually draw directly $Z_i$ and compute $e^{\sum_i Z_i}$.

\begin{table}
    \begin{center}
        \begin{tabular}[htb]{|r|l|l|l|}
            \hline N & dt & error & std err \\
            \hline \hline
            16 & 0.0625 & 0.0402 & 0.00602 \\
            32 & 0.0312 & 0.0214 & 0.00351 \\
            64 & 0.0156 & 0.0102 & 0.00176 \\
            128 & 0.00781 & 0.00519 & 0.000931 \\
            256 & 0.00391 & 0.00253 & 0.000431 \\
            512 & 0.00195 & 0.00123 & 0.000201 \\
            1024 & 0.000977 & 0.000621 & 9.97e-5 \\
            2048 & 0.000488 & 0.00031 & 4.93e-5 \\
            4096 & 0.000244 & 0.000166 & 2.57e-5 \\
            8192 & 0.000122 & 8.22e-5 & 1.33e-5 \\
            16384 & 6.1e-5 & 3.9e-5 & 6.94e-6 \\
            \hline
        \end{tabular}
    \end{center}

    \bigskip

    \caption{Mesh points (N), time steps (dt), strong error (error), and standard error (std err) of the Euler method for $\mathrm{d}X_t/\mathrm{d}t = W_t X_t$ for each mesh resolution $N$, with initial condition $X_0 \sim \mathcal{N}(0, 1)$ and a standard Wiener process noise $\{W_t\}_t$, on the time interval $I = [0.0, 1.0]$, based on $M = 200$ sample paths for each fixed time step, with the target solution calculated with $65536$ points. The order of strong convergence is estimated to be $p = 1.001$, with the 95\% confidence interval $[0.936, 1.0664]$.}

    \label{tablinearhomogeneousrode}
\end{table}

Hence, once an Euler approximation of \cref{linearhomogeneousrode} is computed, along with realizations $\{W_{t_i}\}_{i=0}^N$ of a sample path of the noise, we consider an exact pathwise solution given by
\begin{equation}
    \label{Xtlinearhomogeneousrode}
    X_{t_j} = X_0 e^{\sum_{i = 0}^{j-1}\left(\frac{1}{2}\left(W_{t_i} + W_{t_{i+1}}\right)(t_{i+1} - t_i) + Z_i\right)},
\end{equation}
for realizations $Z_i$ drawn from a normal distributions given by \cref{linearhomogeneousZidistribution}. \cref{samplepathslinearhomogeneousrode} shows an approximate solution and a few sample paths of exact solutions associated with the given realizations of the noise on the mesh points.

The function $f(t, x, y) = yx$ associated with equation \cref{linearhomogeneousrode} does not satisfy the uniform Lipschitz condition \cref{LtLXbound}. Nevertheless, thanks to the form of the exact solution \cref{Xtlinearhomogeneousrodesolution} the form of the approximation \cref{Xtlinearhomogeneousrode} and estimates such as \cref{Xtlinearhomogeneousrodeestimateexponentialofnomral}, the proof of \cref{thmsemimartingale} can be adapted to show strong order $1$ convergence in this case, as well.

\begin{figure}[htb]
    \centerline{\includegraphics[width=0.8\textwidth]{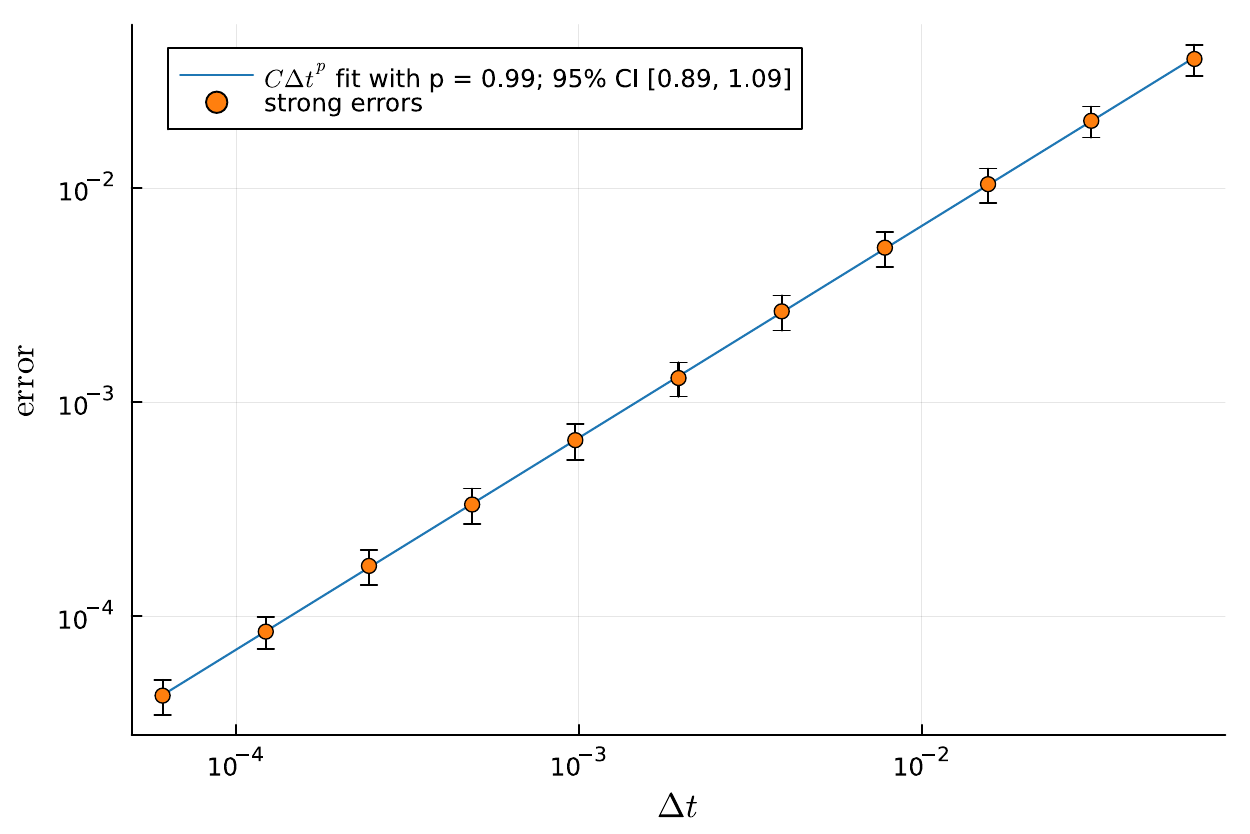}}
    \caption{Estimated order of convergence of the strong error of the Euler method for $\mathrm{d}X_t/\mathrm{d}t = W_t X_t$, based on \cref{tablinearhomogeneousrode}.}
    \label{figlinearhomogeneousrode}
\end{figure}

\Cref{tablinearhomogeneousrode} shows the estimated strong error obtained in the numerical experiments, based on a target solution computed with the exact distribution. \Cref{figlinearhomogeneousrode} illustrates the order of convergence.

\subsubsection{Fractional Brownian motion noise}
\label{secfBmnoise}

Now, we consider again a linear equation, this time of the form
\begin{equation}
    \label{linearnonhomogeneousfbm}
    \begin{cases}
        \displaystyle \frac{\mathrm{d}X_t}{\mathrm{d} t} = -X_t + B^H_t, \qquad 0 \leq t \leq T, \\
        \left. X_t \right|_{t = 0} = X_0,
      \end{cases}
\end{equation}
and where the noise $\{B^H_t\}_t$ is assumed to be a fractional Brownian motion (fBm) with Hurst parameter $0 < H < 1$. It turns out that, for $0 < H < 1/2$, the order of convergence is $H + 1/2,$ while for $1/2\leq H < 1,$ the order of convergence is 1. In any case, this is higher than the order $H$ estimated in previous works. The same seems to hold for a nonlinear dependency on the fBm, but the proof is more involved, depending on a fractional It\^o formula (see \cite[Theorem 4.2.6]{BHOB2008}, \cite[Theorem 4.1]{BENDER200381}and \cite[Theorem 2.7.4]{Mishura2008}), based on the Wick It\^o Skorohod (WIS) integral (see \cite[Chapter 4]{BHOB2008}). A corresponding WIS isometry is also needed (see e.g. \cite[Theorem 4.5.6]{BHOB2008}), involving Malliavin calculus and fractional derivatives. For these reasons, we leave the nonlinear case to a subsequent work and focus on this simple linear example, which suffices to illustrate the peculiarity of the dependence on $H$ of the order of convergence. For this linear equation, the proof of convergence is done rigorously below, with the framework developed in the main article.

We need to estimate the last term in \cref{expectedestimateglobalerrorintegral} of \cref{propbasicstrongestimate}, involving the steps of the term $f(t, x, y) = -x + y$, which in this case reduce to
\begin{equation}
    \label{stepfBm}
    f(s, X_{\tau^N(s)}^N, Y_s) - f(\tau^N(s), X_{\tau^N(s)}^N, Y_{\tau^N(s)}) = B^H_s - B^H_{\tau^N(s)},
\end{equation}
for $0 \leq s \leq T$. There are several ways to represent an fBm (see e.g. \cite{BHOB2008, Mishura2008}). We use the formula (see \cite[eq. (2.1)]{MandelbrotVanNess1968} or \cite[eq. (1.1)]{BHOB2008})
\begin{multline}
    \label{BHtintegralformula}
    B^H_t = \frac{1}{\Gamma(H + 1/2)}\left(\int_{-\infty}^0 \left( (t-s)^{H-1/2} - (-s)^{H-1/2}\right) \;\mathrm{d}W_s \right. \\
    \left. + \int_0^t (t - s)^{H-1/2} \;\mathrm{d}W_s\right),
\end{multline}
where $\Gamma(\cdot)$ is the well-known Gamma function. For the step, \cref{BHtintegralformula} means that
\begin{multline}
    \label{BHtintegralformulastep}
    B^H_s - B^H_{\tau^N(s)} = \frac{1}{\Gamma(H + 1/2)}\left(\int_{-\infty}^{\tau^N(s)} \left( (s-\xi)^{H-1/2} - (\tau^N(s)-\xi)^{H-1/2}\right) \;\mathrm{d}W_\xi \right. \\
    \left. + \int_{\tau^N(s)}^s (s - \xi)^{H-1/2} \;\mathrm{d}W_\xi\right).
\end{multline}
Then, using the stochastic Fubini Theorem to exchange the order of integration (see, again, \cite[Section IV.6]{Protter2005}),
\begin{equation}
    \label{integralofstepfBm}
    \begin{aligned}
        \int_0^{t_j} & \left( f(s, X_{\tau^N(s)}^N, Y_s) - f(\tau^N(s), X_{\tau^N(s)}^N, Y_{\tau^N(s)}) \right)\;\mathrm{d}s \\
        & = \frac{1}{\Gamma(H + 1/2)}\int_0^{t_j} \int_{-\infty}^{\tau^N(s)} \left( (s-\xi)^{H-1/2} - (\tau^N(s)-\xi)^{H-1/2}\right) \;\mathrm{d}W_\xi \;\mathrm{d}s \\
        & \qquad + \frac{1}{\Gamma(H + 1/2)}\int_0^{t_j} \int_{\tau^N(s)}^s (s - \xi)^{H-1/2} \;\mathrm{d}W_\xi \;\mathrm{d}s \\
        & = \frac{1}{\Gamma(H + 1/2)}\int_{-\infty}^{0} \int_{0}^{t_j} \left( (s-\xi)^{H-1/2} - (\tau^N(s)-\xi)^{H-1/2}\right) \;\mathrm{d}s \;\mathrm{d}W_\xi \\
        & \qquad + \frac{1}{\Gamma(H + 1/2)}\int_{0}^{t_j} \int_{\tau^N(\xi)+\Delta t_N}^{t_j} \left( (s-\xi)^{H-1/2} - (\tau^N(s)-\xi)^{H-1/2}\right)  \;\mathrm{d}s \;\mathrm{d}W_\xi\\
        & \qquad + \frac{1}{\Gamma(H + 1/2)}\int_0^{t_j} \int_\xi^{\tau^N(\xi) + \Delta t_N} (s - \xi)^{H-1/2} \;\mathrm{d}s \;\mathrm{d}W_\xi. \\
    \end{aligned}
\end{equation}

For the first term, notice $\sigma \mapsto 1/(\sigma - \xi)^{H-1/2}$ is continuously differentiable on the interval $\sigma > \xi$, so that
\[
    (s-\xi)^{H-1/2} - (\tau^N(s)-\xi)^{H-1/2} = - (H-1/2)\int_{\tau^N(s)}^s (\sigma - \xi)^{H - 3/2} \;\mathrm{d}\sigma.
\]
Thus,
\begin{multline*}
    \int_{0}^{t_j} \left( (s-\xi)^{H-1/2} - (\tau^N(s)-\xi)^{H-1/2}\right) \;\mathrm{d}s \\
    = (H-1/2)\int_{0}^{t_j} \int_{\tau^N(s)}^s (\sigma - \xi)^{H - 3/2} \;\mathrm{d}\sigma \;\mathrm{d}s.
\end{multline*}
Exchanging the order of integration yields
\begin{align*}
    \int_{0}^{t_j} \left( (s-\xi)^{H-1/2} \right. & \left. - (\tau^N(s)-\xi)^{H-1/2} \right) \;\mathrm{d}s \\
    & = (H-1/2)\int_{0}^{t_j} \int_{\sigma}^{\tau^N(\sigma) + \Delta t_N} (\sigma - \xi)^{H - 3/2} \;\mathrm{d}s \;\mathrm{d}\sigma \\
    & = (H-1/2)\int_{0}^{t_j} \left(\tau^N(\sigma) + \Delta t_N - \sigma\right) (\sigma - \xi)^{H - 3/2} \;\mathrm{d}\sigma.
\end{align*}
Hence,
\begin{multline*}
    \left|\int_{0}^{t_j} \left( (s-\xi)^{H-1/2} - (\tau^N(s)-\xi)^{H-1/2} \right) \;\mathrm{d}s\right| \\
    \leq (1/2 - H)\int_{0}^{t_j} \Delta t_N (\sigma - \xi)^{H - 3/2} \;\mathrm{d}\sigma.
\end{multline*}
Now, using the Lyapunov inequality and the It\^o isometry, and using the same trick as above,

\begin{align*}
    & \mathbb{E}\left[\left|\int_{-\infty}^{0} \int_{0}^{t_j} \left( (s-\xi)^{H-1/2} - (\tau^N(s)-\xi)^{H-1/2}\right) \;\mathrm{d}s \;\mathrm{d}W_\xi \right|\right] \\
    & \qquad\qquad \leq \left(\int_{-\infty}^{0} \left(\int_{0}^{t_j} \left( (s-\xi)^{H-1/2} - (\tau^N(s)-\xi)^{H-1/2}\right) \;\mathrm{d}s \right)^2 \;\mathrm{d}\xi \right)^{1/2} \\
    & \qquad\qquad \leq \Delta t_N \left(\int_{-\infty}^{0} \left( (1/2 - H)\int_0^{t_j} (\sigma - \xi)^{H-3/2} \;\mathrm{d}\sigma \right)^2 \;\mathrm{d}\xi \right)^{1/2} \\
    & \qquad\qquad \leq (1/2 - H)\Delta t_N \left(\int_{-\infty}^{0} \left(\int_0^T (\sigma - \xi)^{H-3/2} \;\mathrm{d}\sigma \right)^2 \;\mathrm{d}\xi \right)^{1/2}.
\end{align*}
Therefore,
\begin{multline}
    \label{firsttermfBm}
    \frac{1}{\Gamma(H + 1/2)}\Delta t_N \mathbb{E}\left[\left|\int_{-\infty}^{0} \int_{0}^{t_j} \left( (s-\xi)^{H-1/2} - (\tau^N(s)-\xi)^{H-1/2}\right) \;\mathrm{d}s \;\mathrm{d}W_\xi \right|\right] \\
    \leq c_H^{(1)}\Delta t_N,
\end{multline}
for a suitable constant $c_H^{(1)}$. We see this term is of order 1 in $\Delta t_N.$

The second term is similar,
\begin{align*}
    \int_{\tau^N(\xi)+\Delta t_N}^{t_j} & \left( (s-\xi)^{H-1/2} - (\tau^N(s)-\xi)^{H-1/2}\right) \;\mathrm{d}s \\ 
    & = (H-1/2)\int_{\tau^N(\xi)+\Delta t_N}^{t_j} \int_{\tau^N(s)}^s (\sigma - \xi)^{H - 3/2} \;\mathrm{d}\sigma \;\mathrm{d}s \\
    & = (H-1/2)\int_{\tau^N(\xi)+\Delta t_N}^{t_j} \int_\sigma^{\tau^N(\sigma) + \Delta t_N} (\sigma - \xi)^{H - 3/2} \;\mathrm{d}s \;\mathrm{d}\sigma \\
    & = (H-1/2)\int_{\tau^N(\xi)+\Delta t_N}^{t_j} \left(\tau^N(\sigma) + \Delta t_N - \sigma\right) (\sigma - \xi)^{H - 3/2} \;\mathrm{d}\sigma.
\end{align*}
Thus,
\begin{multline*}
    \left| \int_{\tau^N(\xi)+\Delta t_N}^{t_j} \left( (s-\xi)^{H-1/2} - (\tau^N(s)-\xi)^{H-1/2}\right) \;\mathrm{d}s \right| \\
    \leq (1/2 - H)\Delta t_N \int_{\tau^N(\xi)+\Delta t_N}^{t_j} (\sigma - \xi)^{H - 3/2} \;\mathrm{d}\sigma.
\end{multline*}
Hence,
\begin{align*}
    & \mathbb{E}\left[\left|\int_{0}^{t_j} \int_{\tau^N(\xi)+\Delta t_N}^{t_j} \left( (s-\xi)^{H-1/2} - (\tau^N(s)-\xi)^{H-1/2}\right) \;\mathrm{d}s \;\mathrm{d}W_\xi\right|\right] \\
    & \qquad\qquad \leq \left(\int_{0}^{t_j} \left(\int_{\tau^N(\xi)+\Delta t_N}^{t_j} \left( (s-\xi)^{H-1/2} - (\tau^N(s)-\xi)^{H-1/2}\right) \;\mathrm{d}s \right)^2 \;\mathrm{d}\xi \right)^{1/2} \\
    & \qquad\qquad \leq \Delta t_N (1/2 - H)\left(\int_{0}^{t_j} \left( \int_{\tau^N(\xi)+\Delta t_N}^{T} (\sigma - \xi)^{H-3/2} \;\mathrm{d}\sigma \right)^2 \;\mathrm{d}\xi \right)^{1/2}.
\end{align*}
Therefore,
\begin{multline}
    \label{secondtermfBm}
    \frac{1}{\Gamma(H + 1/2)}\mathbb{E}\left[\left|\int_{0}^{t_j} \int_{\tau^N(\xi)+\Delta t_N}^{t_j} \left( (s-\xi)^{H-1/2} - (\tau^N(s)-\xi)^{H-1/2}\right)  \;\mathrm{d}s \;\mathrm{d}W_\xi\right|\right] \\
    \leq c_H^{(2)}\Delta t_N,
\end{multline}
for a possibly different constant $c_H^{(2)}$. This term is also of order 1.

For the last term, we have
\begin{multline*}
    0 \leq \int_\xi^{\tau^N(\xi) + \Delta t_N} (s - \xi)^{H-1/2} \;\mathrm{d}s = \frac{1}{H + 1/2} (\tau^N(\xi) + \Delta t_N - \xi)^{H + 1/2} \\
    \leq \frac{1}{H + 1/2} \Delta t_N^{H + 1/2}.
\end{multline*}
so that, using the Lyapunov inequality and the It\^o isometry
\begin{multline*}
    \mathbb{E}\left[\left|\int_0^{t_j} \int_\xi^{\tau^N(\xi) + \Delta t_N} (s - \xi)^{H-1/2} \;\mathrm{d}s \;\mathrm{d}W_\xi\right|\right] \\
    \leq \left( \int_0^{t_j} \left(\int_\xi^{\tau^N(\xi) + \Delta t_N} (s - \xi)^{H-1/2} \;\mathrm{d}s\right)^2 \;\mathrm{d}\xi\right)^{1/2} \\ 
    \leq \left( \int_0^{t_j} \Delta t_N^{2H + 1} \;\mathrm{d}\xi\right)^{1/2} \leq t_j^{1/2} \Delta t_N^{H + 1/2}.
\end{multline*}
Therefore,
\begin{equation}
    \label{thirdtermfBm}
    \frac{1}{\Gamma(H + 1/2)}\mathbb{E}\left[\left|\int_0^{t_j} \int_\xi^{\tau^N(\xi) + \Delta t_N} (s - \xi)^{H-1/2} \;\mathrm{d}s \;\mathrm{d}W_\xi\right|\right] \leq c_H^{(3)} \Delta t_N^{H + 1/2},
\end{equation}
for a third constant $c_H^{(3)}$.

Putting the three estimates \cref{firsttermfBm}, \cref{secondtermfBm}, \cref{thirdtermfBm} in \cref{integralofstepfBm} we find that
\begin{multline}
    \mathbb{E}\left[\left|\int_0^{t_j} \left( f(s, X_{\tau^N(s)}^N, Y_s) - f(\tau^N(s), X_{\tau^N(s)}^N, Y_{\tau^N(s)}) \right)\;\mathrm{d}s\right|\right] \\
    \leq c_H^{(4)} \Delta t_N + c_H^{(3)} \Delta t_N^{H + 1/2},
\end{multline}
where $c_H^{(4)} = c_H^{(1)} + c_H^{(2)}$. Using this estimate in \cref{propbasicstrongestimate} shows that the Euler method is of order $H + 1/2$, when $0 < H < 1/2$, and is of order 1, when $1/2 \leq H < 1$, having in mind that $H=1/2$ corresponds to the classical Wiener process.

In summary, we have proved the following result.
\begin{theorem}
    Consider the equation \cref{linearnonhomogeneousfbm} where $\{B^H_t\}_t$ is a fractional Brownian motion (fBm) with Hurst parameter $0 < H < 1$. Suppose the initial condition $X_0$ satisfies
    \begin{equation}
        \label{EX0square2b}
        \mathbb{E}[\|X_0\|^2] < \infty.
    \end{equation}
    Then, the Euler scheme for this initial value problem is of strong order $H+1/2$, for $0 < H < 1/2$, and is of order $1$, for $1/2 \leq H < 1$. More precisely,
    \begin{equation}
        \max_{j=0, \ldots, N}\mathbb{E}\left[ \left| X_{t_j} - X_{t_j}^N \right| \right] \leq c_1 \Delta t_N + c_2 \Delta t_N^{H + 1/2}, \qquad \forall N \in \mathbb{N},
    \end{equation}
    for suitable constants $c_1, c_2 \geq 0$.
\end{theorem}

As for illustrating numerically the order of strong convergence, although the above linear equation has the explicit solution
\begin{equation}
    X_t = e^{-t}X_0 + \int_0^t e^{-(t-s)}B^H_s\;\mathrm{d}s,
\end{equation}
computing a distributionally exact solution of this form is a delicate process. Thus we check the convergence numerically by comparing the approximations with another Euler approximation on a much finer mesh.

\begin{figure}[htb]
    \centerline{\includegraphics[width=0.8\textwidth]{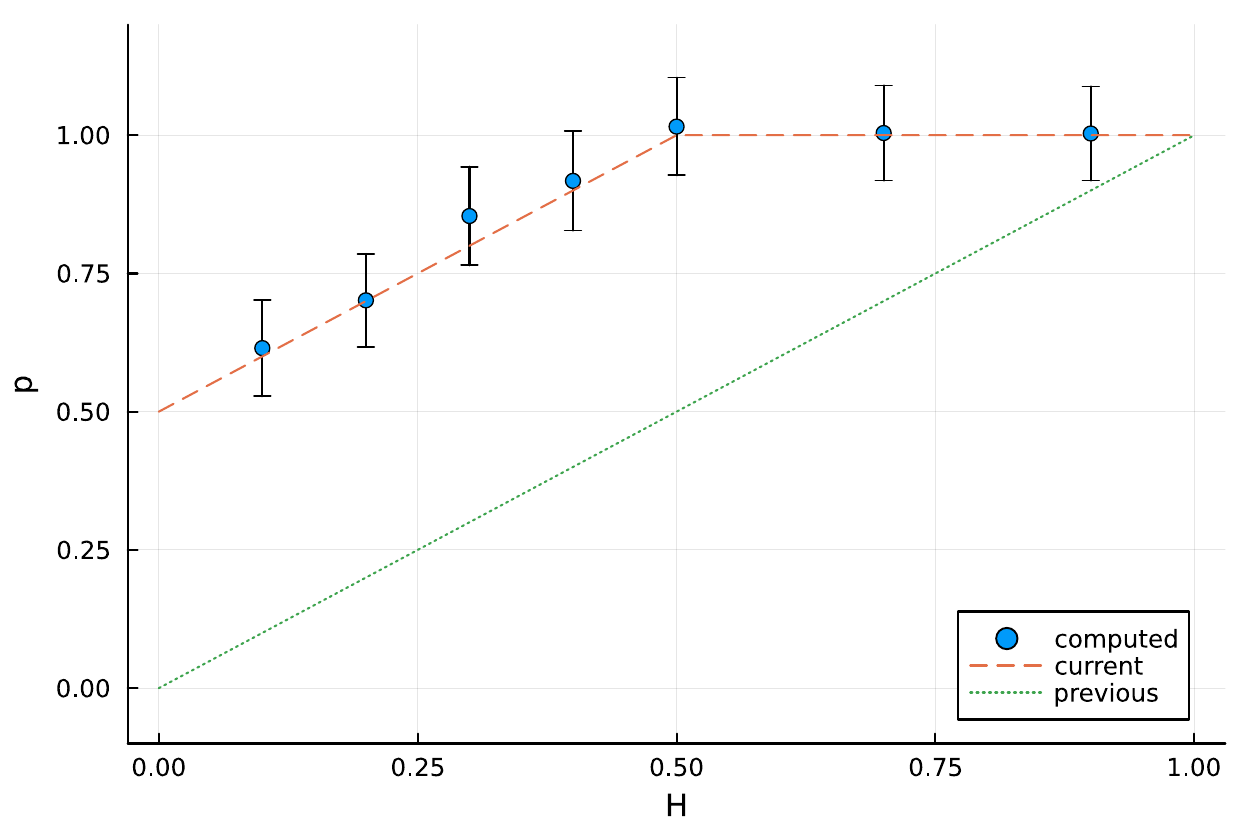}}
    \caption{Order $p$ of strong convergence of the Euler method for \cref{eqlinearfbm} at different values of the Hurst parameter $H$ (scattered plot) along with the theoretical value $p=\min\{H + 1/2, 1\}$ proved here (dashed line), and the previously known theoretical value $p=H$ (dotted).}
    \label{figorderdepHfBm}
\end{figure}

The fractional Browninan motion noise term is generated with the $\mathcal{O}(N\log N)$ fast Fourier transform (FFT) method of Davies and Harte, as presented in \cite{DiekerMandjes2003} (see also \cite[Section 14.4]{HanKloeden2017}). \cref{taborderdepHfBm} shows the obtained convergence estimates, for a number of Hurst parameters, which is illustrated in \cref{figorderdepHfBm}, matching the theoretical estimate of $p = \min\{H+1/2, 1\}.$

\begin{table}
    \begin{center}
        \begin{tabular}[htb]{|c|c|c|c|}
            \hline $H$ & $p$ & $p_{\textrm{min}}$ & $p_{\textrm{max}}$ \\
            \hline \hline
            0.1  & 0.640103 &  0.582705 &  0.697128 \\
            0.2  & 0.742881 &  0.681001 &  0.804453 \\
            0.3  & 0.794755 &  0.735716 &  0.853437 \\
            0.4  & 0.921899 &  0.865631 &  0.978194 \\
            0.5  & 1.00777  &  0.94886  &  1.06662 \\
            0.7  & 1.00637  &  0.945425 &  1.06735 \\
            0.9  & 1.00559  &  0.94014  &  1.07104 \\
            \hline
        \end{tabular}
    \end{center}
    \bigskip

    \caption{Hurst parameter $H$, order $p$ of strong convergence, and 95\% confidence interval $p_{\textrm{min}}$ and $p_{\textrm{max}}$, for a number of Hurst values.}
    \label{taborderdepHfBm}
\end{table}

\subsubsection{Population dynamics with harvest}
\label{secpopdyn}

Here, we consider a population dynamics problem modeled by a logistic equation with random coefficients, loosely inspired by \cite[Section 15.2]{HanKloeden2017}, with an extra term representing harvest:
\begin{equation}
    \label{rodepopulationdynamics}
    \frac{\mathrm{d}X_t}{\mathrm{d}t} = \Lambda_t X_t (1 - \frac{X_t}{r}) - \alpha H_t \frac{X_t}{r + X_t},
\end{equation}
where $r, \alpha > 0$ are constants and $\{\Lambda_t\}_{t \geq 0}$ and $\{H_t\}_{t \geq 0}$ are stochastic processes. The first term is a logistic growth, with carrying capacity $r$ and random specific growth $\Lambda_t$. The second term is the harvest term, with $\alpha$ being an overall modulation factor, the term $X_t / (r + X_t)$ accounting for the scarcity of the population when $X_t \ll r$, and with $0 \leq H_t \leq 1$ being a random harvest factor. 

More specifically, $\{\Lambda_t\}_{t \geq 0}$ is given by
\[
    \Lambda_t = \gamma (1 + \varepsilon \sin(G_t)),
\]
where $\gamma > 0,$ $0 < \varepsilon < 1,$ and $\{G_t\}_{t\geq 0}$ is a geometric Brownian motion (gBm) process. A Wiener process is a natural choice, but we choose a gBm process instead since it is a multiplicative noise.

The harvest term $\{H_t\}_{t\geq 0}$ is a ``Poisson step'' process of the form
\[
    H_t = S_{N_t},
\]
where $\{N_t\}_{t\geq 0}$ is a Poisson point-process with rate $\lambda$, $S_0 = 0$, and the $S_i$, for $i=1, 2, \ldots$, are independent and identically distributed random variables with nonnegative values within the interval $[0, 1]$, independent also of the Poisson counter $\{N_t\}_{t\geq 0}$.

We suppose the initial condition is nonnegative and bounded almost surely, i.e.
\[
    0 \leq X_0 \leq R,
\]
for some $R > r$.

The noise process $\{\Lambda_t\}_{t \geq 0}$ itself satisfies
\[
    0 < \gamma - \varepsilon \leq \Lambda_t \leq \gamma + \varepsilon < 2\gamma, \qquad \forall t \geq 0.
\]

Define $f:\mathbb{R} \times \mathbb{R} \times \mathbb{R}^2 \rightarrow \mathbb{R}$ by
\[
    f(t, x, y) = \begin{cases}
        \displaystyle \gamma (1 + \varepsilon \sin(y_1)) x(r - x) - \alpha \frac{xy_2}{r + x}, & x \geq 0, \\
        0, & x < 0.
    \end{cases}
\]
The equation \cref{rodepopulationdynamics} becomes
\[
    \frac{\mathrm{d}X_t}{\mathrm{d}t} = f(t, X_t, Y_t),
\]
where $\{Y_t\}_{t\geq 0}$ is the vector-valued process given in coordinates by $Y_t = (G_t, H_t)$.

Notice that $f(t, x, y) = 0$, for $x < 0$, for arbitrary $y=(y_1, y_2)$, while $f(t, x, y) < 0$, for $x \geq r$, $y_2 \geq 0$, and for arbitrary $y_1$. Since the noise $y_2 = H_t$ is always nonnegative, we see that the interval $0 \leq x \leq R$ is positively invariant and attracts the orbits with a nonnegative initial condition. Thus, the pathwise solutions of the initial-value problem under consideration are almost surely bounded as well.

The function $f=f(t, x, y)$ is continuously differentiable infinitely many times and with bounded derivatives within the positively invariant region. Hence, within the region of interest, all the conditions of \cref{thmsemimartingale} hold and the Euler method is of strong order 1.

We simulate numerically the solutions of the above problem, with $\gamma = 0.8,$ $\varepsilon = 0.3,$ $r = 1.0,$ and $\alpha = \gamma r^2.$ The geometric Brownian motion process $\{G_t\}_{t\geq 0}$ is taken with drift coefficient $\mu = 1.0,$ diffusion coefficient $\sigma = 0.8,$ and initial condition $y_0 = 1.0.$ The Poisson process $\{N_t\}_{t \geq 0}$ is taken with rate $\lambda = 15.0$. And the step process $\{H_t\}_{t \geq 0}$ is taken with steps following a Beta distribution with shape parameters $\alpha = 5.0$ and $\beta = 7.0$. The initial condition $X_0$ is taken to be a Beta distribution with shape parameters $\alpha = 7.0$ and $\beta = 5.0$, hence we can take $R = 1$.

Notice that we can write
\[
    \frac{\mathrm{d}X_t}{\mathrm{d}t} = \frac{X_t}{r + X_t} \left(r\Lambda_t - \alpha H_t - \frac{\Lambda_t}{r} X_t^2\right).
\]
Hence, when $\alpha H_t \geq r\Lambda_t$, the population decays for arbitrary $X_t > 0$, leading to an extinction of the population. The parameters chosen above keep the population from extinction but may often get close to the critical values.

\begin{figure}[htb]
    \centerline{\includegraphics[width=0.8\textwidth]{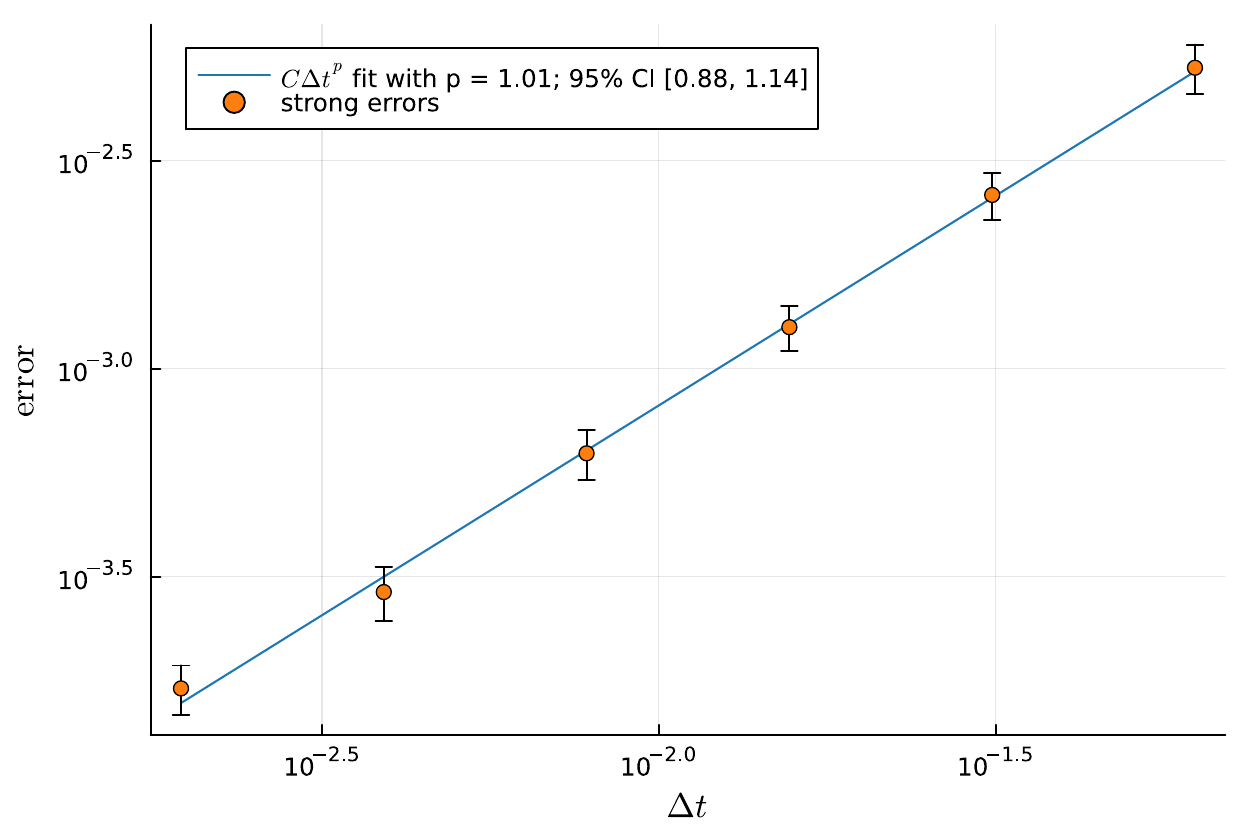}}
    \caption{Order of convergence of the strong error of the Euler method for the population dynamics equation \cref{rodepopulationdynamics}, based on \cref{tabpopdyn}.}
    \label{figpopdyn}
\end{figure}

\begin{figure}[htb]
    \centerline{\includegraphics[width=0.8\textwidth]{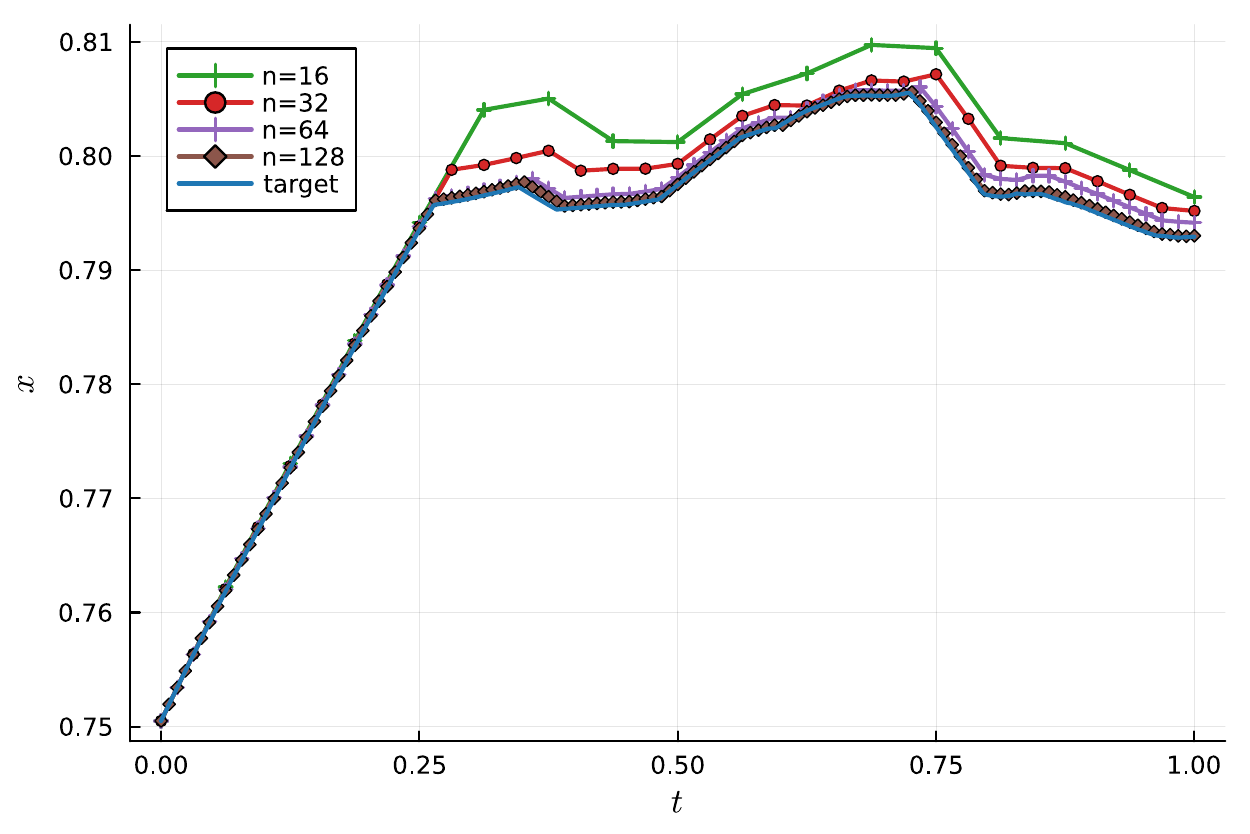}}
    \caption{An approximation sequence of a sample path solution of the population dynamics equation \cref{rodepopulationdynamics}.}
    \label{figsamplepopdyn}
\end{figure}

\cref{tabpopdyn} shows the estimated strong error obtained for each mesh resolution, while \cref{figpopdyn} illustrates the order of convergence, estimated to be close enough to the theoretical value of strong order 1. Finally, \cref{figsamplepopdyn} shows an approximation sequence of a sample path.

\begin{table}
    \begin{center}
        \begin{tabular}[htb]{|r|l|l|l|}
            \hline N & dt & error & std err \\
            \hline \hline
            16 & 0.0625 & 0.00529 & 0.000367 \\
            32 & 0.0312 & 0.00262 & 0.000172 \\
            64 & 0.0156 & 0.00126 & 7.86e-5 \\
            128 & 0.00781 & 0.000626 & 4.35e-5 \\
            256 & 0.00391 & 0.00029 & 2.2e-5 \\
            512 & 0.00195 & 0.00017 & 1.18e-5 \\
            \hline
        \end{tabular}
    \end{center}

    \bigskip

    \caption{Mesh points (N), time steps (dt), strong error (error), and standard error (std err) of the Euler method for population dynamics for each mesh resolution $N$, with initial condition $X_0 \sim \mathrm{Beta}(7.0, 5.0)$ and gBm and step process noises, on the time interval $I = [0.0, 1.0]$, based on $M = 100$ sample paths for each fixed time step, with the target solution calculated with $262144$ points. The order of strong convergence is estimated to be $p = 1.009$, with the 95\% confidence interval $[0.9603, 1.0571]$.}

    \label{tabpopdyn}
\end{table}

\subsubsection{A toggle-switch model for gene expression}

Here, we consider the toggle-switch model in \cite[Section 7.8]{Asai2016}, originated from \cite{VerdCrombachJaeger2014}; see also \cite{StrasserTheisMarr2012}.

Toogle switches in gene expression consist of genes that mutually repress each other and exhibit two stable steady states of ON and OFF. It is a regulatory mechanism which is active during cell differentiation and is believed to act as a memory device, able to choose and maintain cell fate decisions.

We consider the following simple model as discussed in \cite[Section 7.8]{Asai2016}, of two interacting genes, $X$ and $Y$, with the concentration of their corresponding protein products at each time $t$ denoted by $X_t$ and $Y_t$. These are stochastic processes defined by the system of equations
\begin{equation}
    \label{toggleswitchsystem}
   \begin{cases}
   \frac{\displaystyle \mathrm{d}X_t}{\displaystyle \mathrm{d} t} = \left( A_t + \frac{\displaystyle X_t^4}{\displaystyle a^4 + X_t^4}\right)\left(\frac{\displaystyle b^4}{\displaystyle b^4 + Y_t^4}\right) - \mu X_t, \\
   \frac{\displaystyle \mathrm{d}Y_t}{\displaystyle \mathrm{d} t} = \left( B_t + \frac{\displaystyle Y_t^4}{\displaystyle c^4 + Y_t^4}\right)\left(\frac{\displaystyle d^4}{\displaystyle d^4 + X_t^4}\right) - \nu Y_t,
   \end{cases}
\end{equation}
on $t \geq 0$, with initial conditions $X_0$ and $Y_0$, where $\{A_t\}_{t\geq 0}$ and $\{B_t\}_{t\geq 0}$ are given stochastic process representing the external activation on each gene; $a$ and $c$ determine the auto-activation thresholds; $b$ and $d$ determine the thresholds for mutual repression; and $\mu$ and $\nu$ are protein decay rates. In this model, the external activations $A_t$ is taken to be a compound Poisson process, while $B_t$ is non-autonomous homogeneous linear It\^o process of the form $\mathrm{d}B_t = (\mu_1 + \mu_2\sin(\vartheta t))B_t\;\mathrm{d}t + \sigma\sin(\vartheta t)B_t\;\mathrm{d}W_t$.

In the simulations below, we use parameters similar to those in \cite[Section 7.8]{Asai2016}. We fix $a = c = 0.25$; $b = d = 0.4$; and $\mu = \nu = 0.75$. The initial conditions are set to $X_0 = Y_0 = 4.0$. The external activation $\{A_t\}$ is a compound Poisson process with Poisson rate $\lambda = 5.0$ and jumps uniformly distributed within $[0.0, 0.5]$, while $\{B_t\}_t$ is taken with $\mu_1 = 0.7,$ $\mu_2 = 0.3,$ $\sigma = 0.3,$ $\vartheta=3\pi,$ and initial condition $B_0 = 0.2.$

We do not have an explicit solution for the equation so we use as target for the convergence an approximate solution via Euler method at a much higher resolution.

\begin{figure}[htb]
    \centerline{\includegraphics[width=0.8\textwidth]{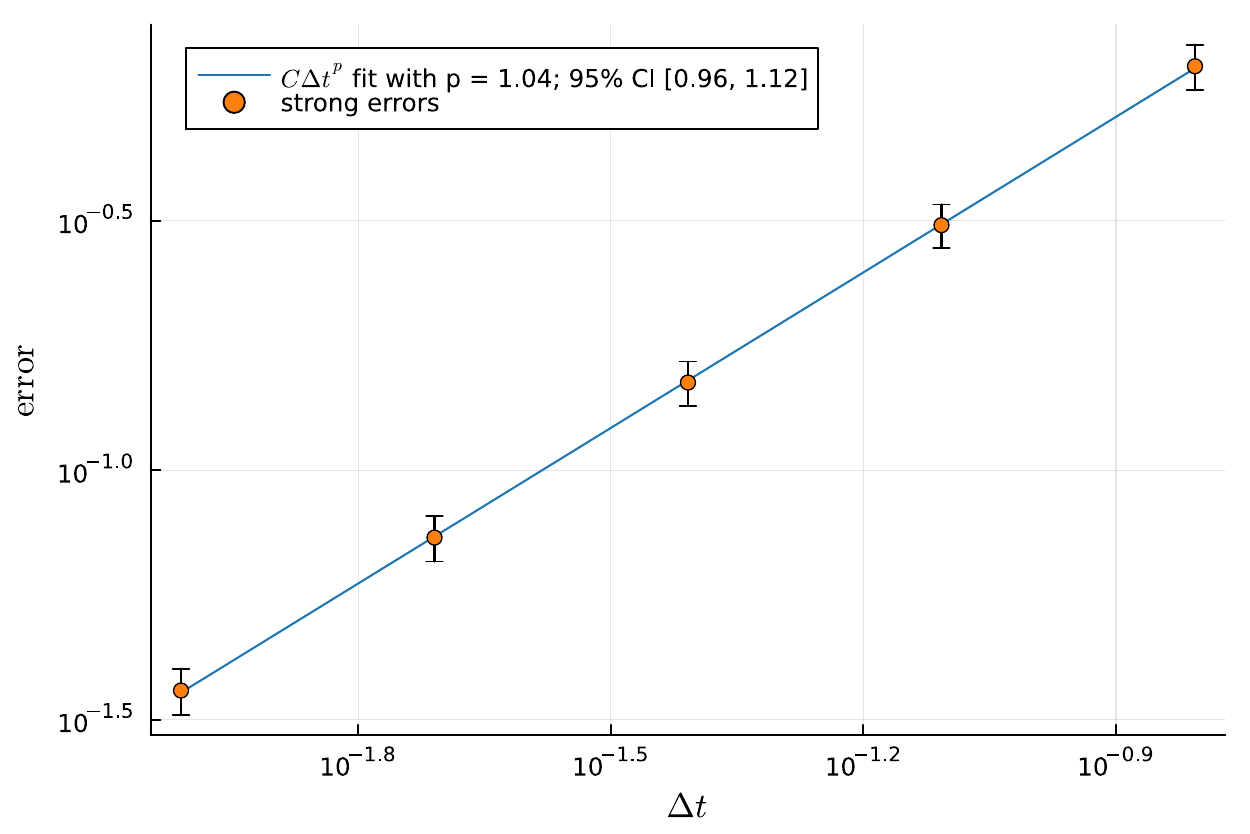}}
    \caption{Order of convergence of the strong error of the Euler method for the toggle-switch model \cref{toggleswitchsystem}, based on \cref{tabletoggleswitch}.}
    \label{figtoggleswitch}
\end{figure}

\begin{figure}[htb]
    \centerline{\includegraphics[width=0.8\textwidth]{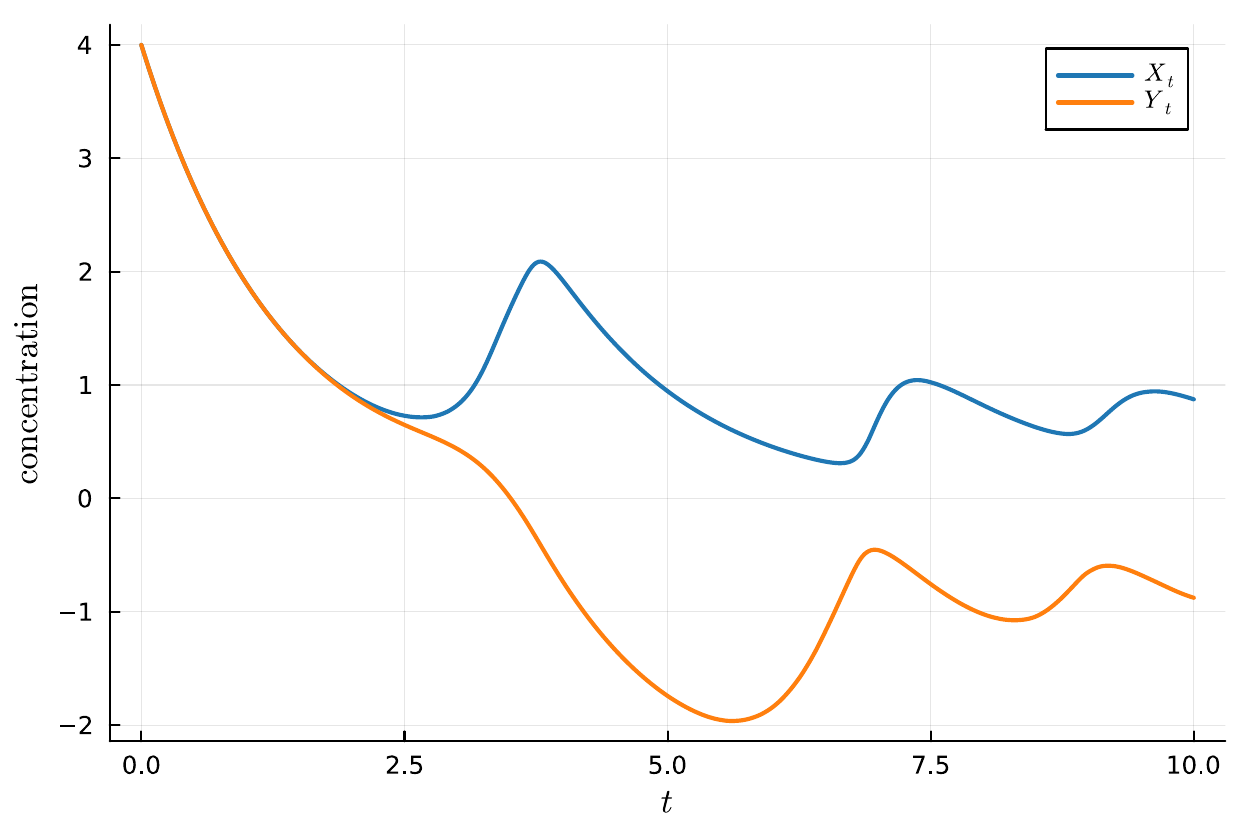}}
    \caption{Evolution of a sample solution of the toggle-switch model \cref{toggleswitchsystem}.}
    \label{figtoggleswitchevolution}
\end{figure}

\begin{figure}[htb]
    \centerline{\includegraphics[width=0.8\textwidth]{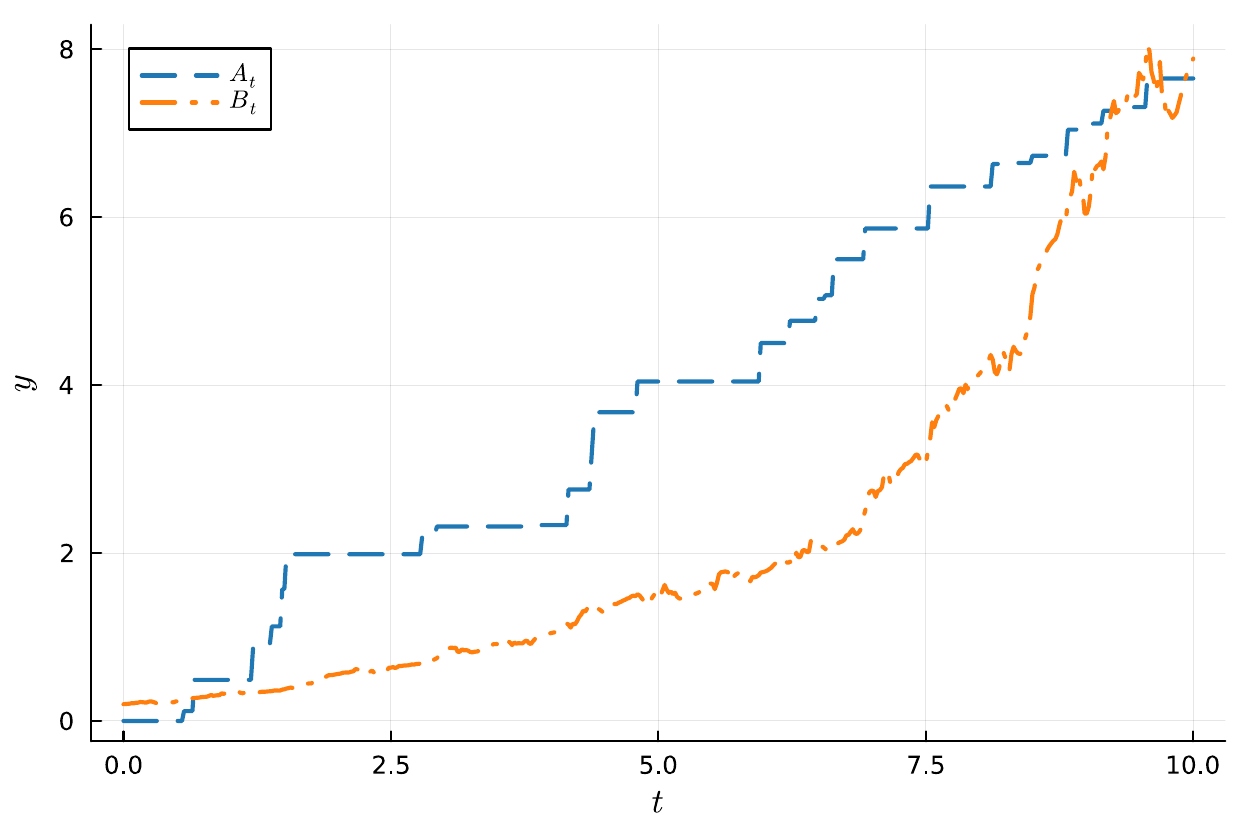}}
    \caption{Noise sample paths for the toggle-switch model \cref{toggleswitchsystem}.}
    \label{figtoggleswitchnoise}
\end{figure}

\Cref{tabletoggleswitch} shows the estimated strong error obtained with the Monte Carlo method, while \cref{figtoggleswitch} illustrates the order of convergence. \Cref{figtoggleswitchevolution} shows a sample solution, while \cref{figtoggleswitchnoise} illustrates the two components $(A_t, B_t)$ of a sample noise.

\begin{table}
    \begin{center}
        \begin{tabular}[htb]{|r|l|l|l|}
            \hline N & dt & error & std err \\
            \hline \hline
            32 & 0.156 & 0.645 & 0.0344 \\
            64 & 0.0781 & 0.31 & 0.0158 \\
            128 & 0.0391 & 0.15 & 0.00774 \\
            256 & 0.0195 & 0.0733 & 0.0039 \\
            512 & 0.00977 & 0.0362 & 0.00195 \\
            \hline
        \end{tabular}
    \end{center}

    \bigskip

    \caption{Mesh points (N), time steps (dt), strong error (error), and standard error (std err) of the Euler method for a toggle-switch model of gene regulation for each mesh resolution $N$, with initial condition $X_0 = 4.0, Y_0 = 4.0$ and coupled compound Poisson process and geometric Brownian motion noises, on the time interval $I = [0.0, 5.0]$, based on $M = 100$ sample paths for each fixed time step, with the target solution calculated with $262144$ points. The order of strong convergence is estimated to be $p = 1.039$, with the 95\% confidence interval $[1.0103, 1.0681]$.}

    \label{tabletoggleswitch}
\end{table}

\subsubsection{Mechanical structure under random Earthquake-like seismic disturbances}

Now we consider a simple mechanical structure model driven by a random disturbance in the form of a transport process, simulating seismic ground-motion excitations, inspired by the model in \cite{BogdanoffGoldbergBernard1961} (see also \cite[Chapter 18]{NeckelRupp2013} and \cite{HousnerJenning1964} with this and other models).

There are a number of models for earthquake-type forcing, such as the ubiquitous Kanai-Tajimi and the Clough-Penzien models, where the noise has a characteristic spectral density, determined by the mechanical properties of the ground layer. The ideia, from \cite{Kanai1957} and \cite{Tajimi1960}, is that the spectrum of the noise at bedrock is characterized by a constant pattern, while at the ground surface it is modified by the vibration property of the ground layer. This interaction between the bedrock and the ground layer is modeled as a stochastic oscillator driven by a zero-mean Gaussian white noise, and whose solution leads to a noise with a characteristic power spectrum.

We follow, however, the Bogdanoff-Goldberg-Bernard model, which takes the form of a transport process noise and models qualitatively well the main shock and subsequent aftershocks. We chose the later not only because of its realistic modeling but also so we can illustrate the improved convergence for such type of transport noise, complementing the other examples. The noise is described in more details shortly. Let us first introduce the model for the vibrations of the mechanical structure.

A single-storey building is considered, with its ground floor centered at position $M_t$ and its ceiling at position $M_t + X_t$. The random process $X_t$ refers to the motion relative to the ground. The ground motion $M_t$ affects the motion of the relative displacement $X_t$ as an excitation force proportional to the ground acceleration $\ddot M_t$. The damping and elastic forces are in effect within the structure. In this framework, the equation of motion for the relative displacement $X_t$ of the ceiling of the single-storey building takes the form

\begin{equation}
    \label{mechanicalstructuremodel}
    \ddot X_t + 2\zeta_0\omega_0\dot X_t + \omega_0^2 X_t = - \ddot M_t.
\end{equation}
where $\zeta_0$ and $\omega_0$ are damping and elastic model parameters depending on the structure of the building.

For the numerical simulations, the second-order equation is written as a system of first-order equations,
\[
    \begin{cases}
        \dot X_t = V_t, \\
        \dot V_t = - \omega_0^2 X_t - 2\zeta_0\omega_0 X_t - Y_t,
    \end{cases}
\]
where $\{V_t\}_t$ is the random velocity of the celing relative to the ground and where $\{Y_t\}_t$ is the stochastic noise excitation term given as the ground acceleration, $Y_t = \ddot M_t$, generated by an Earthquake and its aftershocks, or any other type of ground motion.

The structure is originally at rest, so we have the conditions
\[
    X_0 = 0, \quad V_0 = \dot X_0 = 0.
\]

In the Bogdanoff-Goldberg-Bernard model \cite{BogdanoffGoldbergBernard1961}, 
the excitation $\ddot M_t$ is made of a composition of oscillating signals with random frequencies, modulated by a linear attack rate followed by an exponential decay. This can be written, more precisely, as
\[
    \sum_{j=1}^n a_j t e^{-\delta_j t}\cos(\omega_j t + \theta_j).
\]

In order to simulate the start of the first shock-wave and the subsequent aftershocks, we modify this model sligthly to be a combination of such terms but at different incidence times. We also remove the attack rate from the excitation to obtain a rougher instantaneous, discontinuous excitation. This jump-like excitation is connected with a square power attack rate for the displacement itself. Finally, for simulation purposes, we model the displacement $M_t$ instead of modeling directly the excitation $\ddot M_t$, but in such a way that the ground-motion excitation follows essentially the proposed signal.

Thus, with this framework in mind, we model the ground displacement as a transport process composed of a series of time-translations of a square-power ``attack" front, with an exponentially decaying tail and an oscillating background wave:
\begin{equation}
    M_t = \sum_{i=1}^k \gamma_i (t - \tau_i)_+^2 e^{-\delta_i (t - \tau_i)}\cos(\omega_i (t - \tau_i)),
\end{equation}
where $k\in \mathbb{N}$ is given, $(t-\tau_i)_+ = \max\{0, t - \tau_i\}$ is the positive part of the function, and the parameters $\gamma_i,$ $\tau_i,$ $\delta_i,$ and $\omega_i$ are all random variables, with $\tau_i$ being exponentially distributed, and $\gamma_i$, $\delta_i$, and $\omega_i$ being uniformly distributed, each with different support values, and all of them independent of each other.

The excitation itself becomes
\begin{align*}
    \ddot M(t) = & 2\sum_{i=1}^k\gamma_i H(t - \tau_i) e^{-\delta_i (t - \tau_i)}\cos(\omega_i (t - \tau_i)) \\
        & + \sum_{i=1}^k\gamma_i (\delta_i^2 - \omega_i^2)(t - \tau_i)_+^2 e^{-\delta_i (t - \tau_i)}\cos(\omega_i (t - \tau_i)) \\
        & -2\sum_{i=1}^k\gamma_i (\delta_i + \omega_i) (t - \tau_i)_+ e^{-\delta_i (t - \tau_i)}\cos(\omega_i (t - \tau_i)) \\
        & +\delta_i\sum_{i=1}^k\omega_i\gamma_i (t - \tau_i)_+^2 e^{-\delta_i (t - \tau_i)}\sin(\omega_i (t - \tau_i)),
\end{align*}
where $H = H(s)$ is the Heaviside function, where, for definiteness, we set $H(s) = 1,$ for $s \geq 1,$ and $H(s) = 0$, for $s < 0$.

More specifically, for the numerical simulations, we use $\zeta_0 = 0.6$ and $\omega_0 = 15$ as the structural parameters. For the transport process, we set $k=8$ and define the random parameters as the following exponential and uniform distributions: $\tau_i \sim \textrm{Exponential}(0.5),$ $\gamma_i \sim \textrm{Uniform}(16.0, 32.0),$ $\delta_i \sim \textrm{Uniform}(12.0, 16.0),$ and $\omega_i \sim \textrm{Uniform}(16\pi, 32\pi).$

\begin{figure}[htb]
    \centerline{\includegraphics[width=0.8\textwidth]{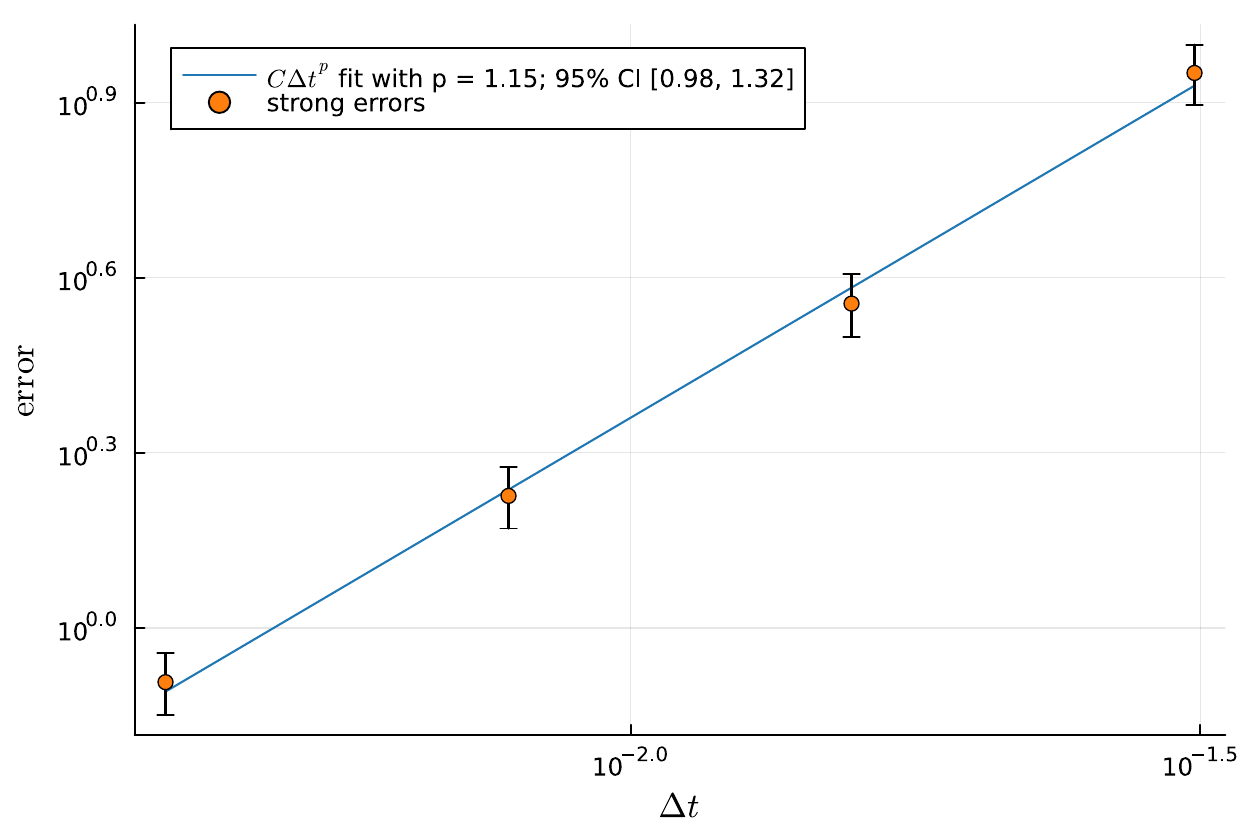}}
    \caption{Order of convergence of the strong error of the Euler method for the mechanical structure model \cref{mechanicalstructuremodel}, based on \cref{tableearthquake}.}
    \label{figearthquake}
\end{figure}

\begin{figure}[htb]
    \centerline{\includegraphics[width=0.8\textwidth]{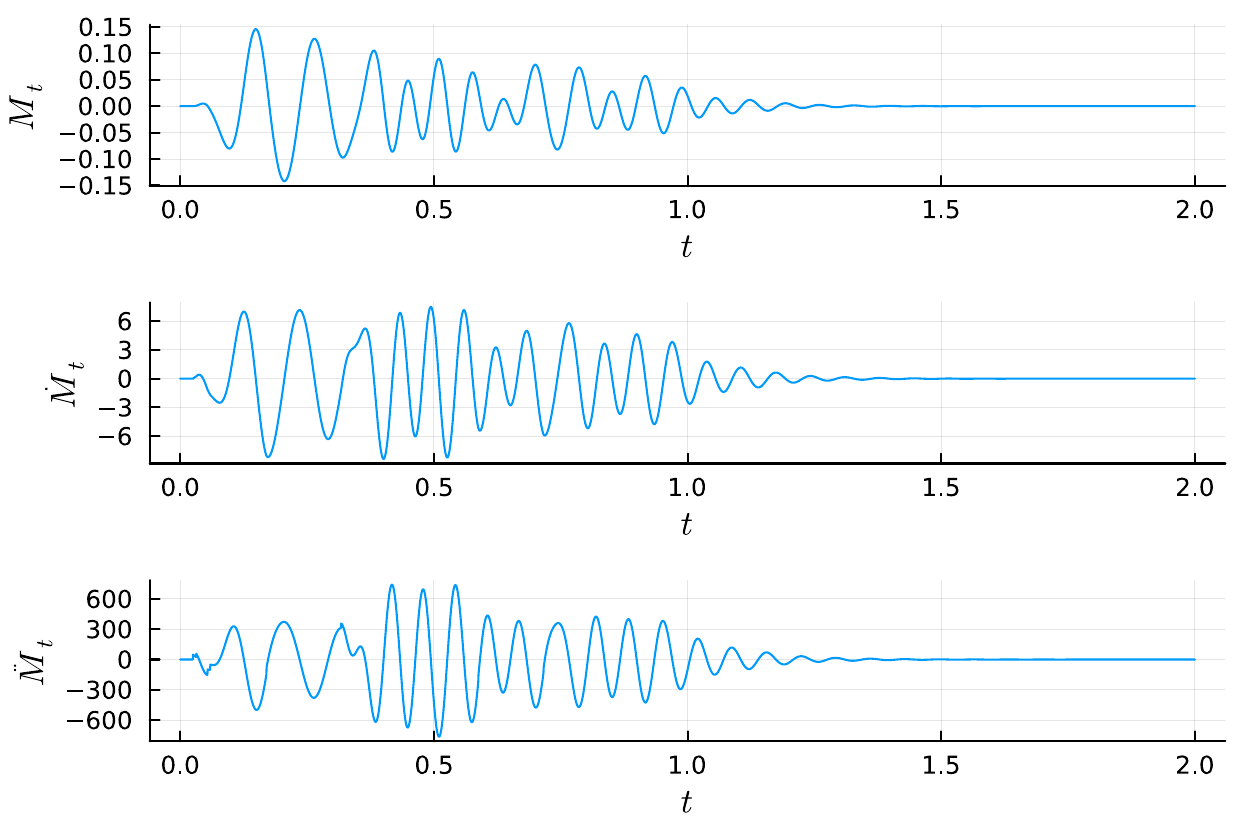}}
    \caption{Sample ground motion, ground velocity, and ground acceleration for the mechanical structure model \cref{mechanicalstructuremodel}.}
    \label{figearthquakenoise}
\end{figure}

\Cref{tableearthquake} shows the estimated strong error obtained with the Monte-Carlo method, while \cref{figearthquake} illustrates the order of convergence. \cref{figearthquakenoise} shows a sample ground motion and the corresponding ground acceleration and its envelope.

\begin{table}
    \begin{center}
        \begin{tabular}[htb]{|r|l|l|l|}
            \hline N & dt & error & std err \\
            \hline \hline
            64 & 0.0312 & 9.27 & 0.721 \\
            128 & 0.0156 & 3.66 & 0.308 \\
            256 & 0.00781 & 1.75 & 0.142 \\
            512 & 0.00391 & 0.838 & 0.0672 \\
            \hline
        \end{tabular}
    \end{center}

    \bigskip

    \caption{Mesh points (N), time steps (dt), strong error (error), and standard error (std err) of the Euler method for mechanical structure model under ground-shaking random excitations for each mesh resolution $N$, with initial condition $X_0 = \mathbf{0}$ and transport process noise, on the time interval $I = [0.0, 2.0]$, based on $M = 100$ sample paths for each fixed time step, with the target solution calculated with $262144$ points. The order of strong convergence is estimated to be $p = 1.147$, with the 95\% confidence interval $[1.0595, 1.2341]$.}

    \label{tableearthquake}
\end{table}

\subsubsection{An actuarial risk model}

A classical model for the surplus $U_t$ at time $t$ of an insurance company is the Cram\'er-Lundberg model (see e.g. \cite{GerberShiu1998}) given by
\[
  U_t = U_0 + \gamma t - \sum_{i=1}^{N_t} C_i,
\]
where $U_0$ is the initial capital, $\gamma$ is a constant premium rate received from the insurees, $C_i$ is a random variable representing the value of the $i$-th claim paid to a given insuree, and $N_t$ is the number of claims up to time $t$. The process $\{N_t\}_t$ is modeled as a Poisson counter, so that the accumulated claims form a compound Poisson process. It is also common to use inhomogeneous Poisson processes and Hawkes self-exciting process, or combinations of such processes for the incidence of the claim, but the classical model uses a homogeneous Poisson counter.

The model above, however, does not take into account the variability of the premium rate received by the company, nor the investiment of the accumulated reserves, among other things. Several diffusion type models have been proposed to account for these and other factors. We consider here a simple model, with a randomly perturbed premium and with variable rentability.

More precisely, we start by rewriting the above expression as the following jump (or impulse) differential equation
\[
  \mathrm{d}U_t = \gamma\;\mathrm{d}t - \mathrm{d}C_t,
\]
where
\[
  C_t = \sum_{i=1}^{N_t} C_i.
\]

The addition of an interest rate $r$ leads to
\[
  \mathrm{d}U_t = r U_t \mathrm{d}t + \gamma\;\mathrm{d}t - \mathrm{d}C_t.
\]

Assuming a premium rate perturbed by a white noise and assuming the interest rate as a process $\{R_t\}_t$, we find
\[
  \mathrm{d}U_t = R_t U_t\;\mathrm{d}t + \gamma\;\mathrm{d}t + \varepsilon\;\mathrm{d}W_t - \mathrm{d}C_t,
\]
so the equation becomes
\[
  \mathrm{d}U_t = (\gamma + R_t U_t)\;\mathrm{d}t + \varepsilon\;\mathrm{d}W_t - \mathrm{d}C_t.
\]

Since we can compute exactly the accumulated claims $C_t$, we subtract it from $U_t$ to get rid of the jump term. We also subtract an Ornstein-Uhlenbeck process, in the classical way to transform an SDE into a RODE. So, defining
\[
  X_t = U_t - C_t - O_t,
\]
where $\{O_t\}_t$ is given by
\[
  \mathrm{d}O_t = -\nu O_t\;\mathrm{d}t + \varepsilon\;\mathrm{d}W_t,
\]
we obtain
\[
  \mathrm{d}X_t = (\gamma + R_t U_t)\;\mathrm{d}t + \nu O_t\;\mathrm{d}t = (\gamma + R_t (X_t + C_t + O_t))\;\mathrm{d}t + \nu O_t\;\mathrm{d}t.
\]

This leads us to the linear random ordinary differential equation
\begin{equation}
    \label{riskmodel}
    \frac{\mathrm{d}X_t}{\mathrm{d}t} = R_t X_t + R_t (C_t + O_t) + \nu O_t + \gamma.
\end{equation}

This equation has the explicit solution
\[
  X_t = X_0 e^{\int_0^t R_s\;\mathrm{d}s} + \int_0^t e^{\int_s^t R_\tau\;\mathrm{d}\tau} (R_s (C_s + O_s) + \nu O_s + \gamma)\;\mathrm{d}s.
\]

As for the interest rate process $\{R_t\}_t$, there is a vast literature with models for it, see e.g. \cite[Chapter 3]{BrigoMercurio2006} in particular Table 3.1. Here, we consider the Dothan model (\cite[Section 3.2.2]{BrigoMercurio2006} of the aforementioned reference), which consists simply of a geometric Brownian motion process
\[
  \mathrm{d}R_t = \mu R_t \;\mathrm{d}t + \sigma R_t\;\mathrm{d}t,
\]
with $R_t = r_0$, where $\mu, \sigma, r_0$ are positive constants. This has an explicit solution
\[
  R_t = r_0 e^{(\mu - \sigma^2/2)t + \sigma W_t},
\]
so that the equation \cref{riskmodel} for $\{X_t\}_t$ is a genuine random ODE.

Once the solution of $\{X_t\}_t$ is obtained, we find an explicit formula for the surplus $U_t = X_t + C_t + O_t$, namely
\[
  U_t = C_t + O_t + X_0 e^{\int_0^t R_s\;\mathrm{d}s} + \int_0^t e^{\int_s^t R_\tau\;\mathrm{d}\tau} (R_s (C_s + O_s) + \nu O_s + \gamma)\;\mathrm{d}s,
\]
with $\{R_t\}_t$ as above.

For the numerical simulations, we use $O_0 = 0.0$, $\nu = 5.0$ and $\varepsilon = 0.8$, for the Ornstein-Uhlenbeck process $\{O_t\}_t$; $\lambda = 8.0$ and $C_i \sim \mathrm{Uniform}(0, 0.2)$, for the compound Poisson process $\{C_t\}$; $R_0 = 0.2$, $\mu = 0.02$ and $\sigma = 0.4$, for the interest rate process $\{R_t\}_t$; and we take $X_0 = 1.0$, so that $U_0 = X_0 + O_0 + R_0 = 1.2$.

\begin{figure}[htb]
    \centerline{\includegraphics[width=0.8\textwidth]{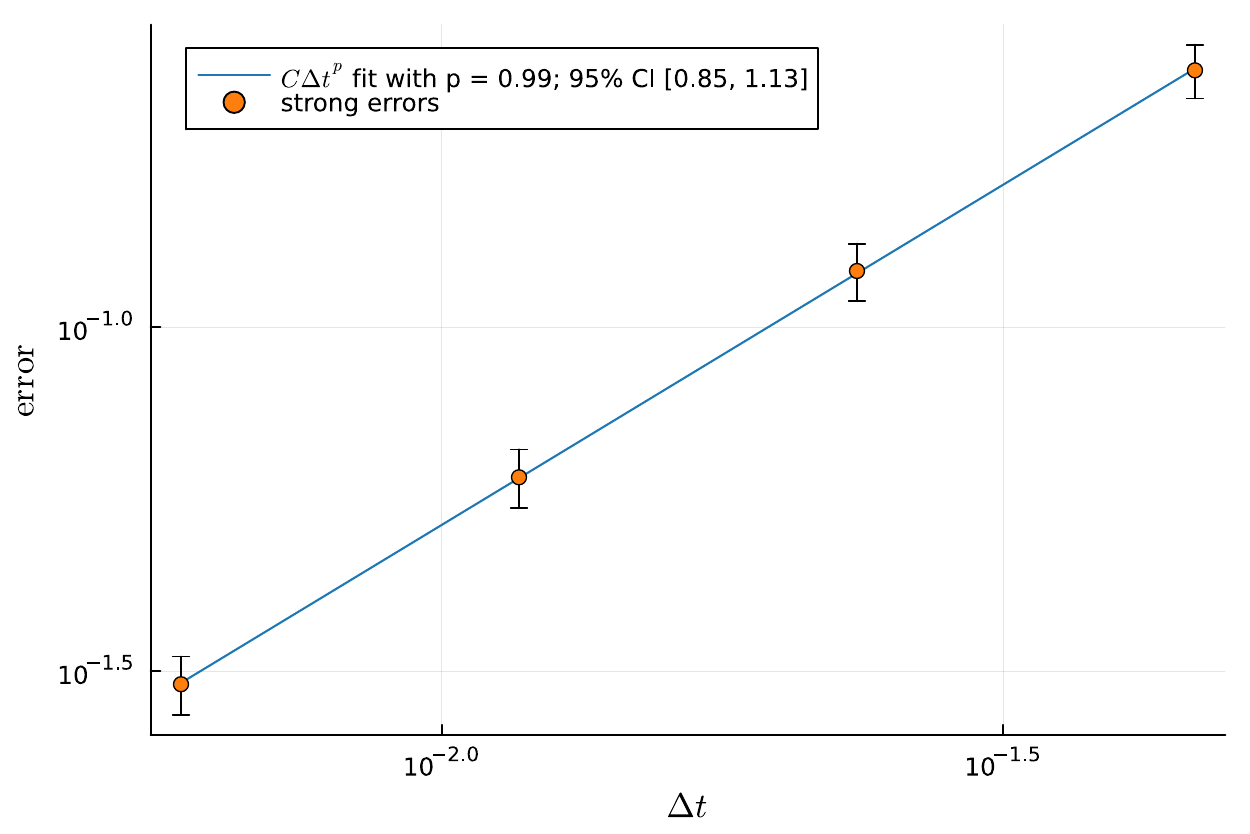}}
    \caption{Order of convergence of the strong error of the Euler method for the actuarial risk model \cref{riskmodel}, based on \cref{tableriskmodel}.}
    \label{figriskmodel}
\end{figure}

\begin{figure}[htb]
    \centerline{\includegraphics[width=0.8\textwidth]{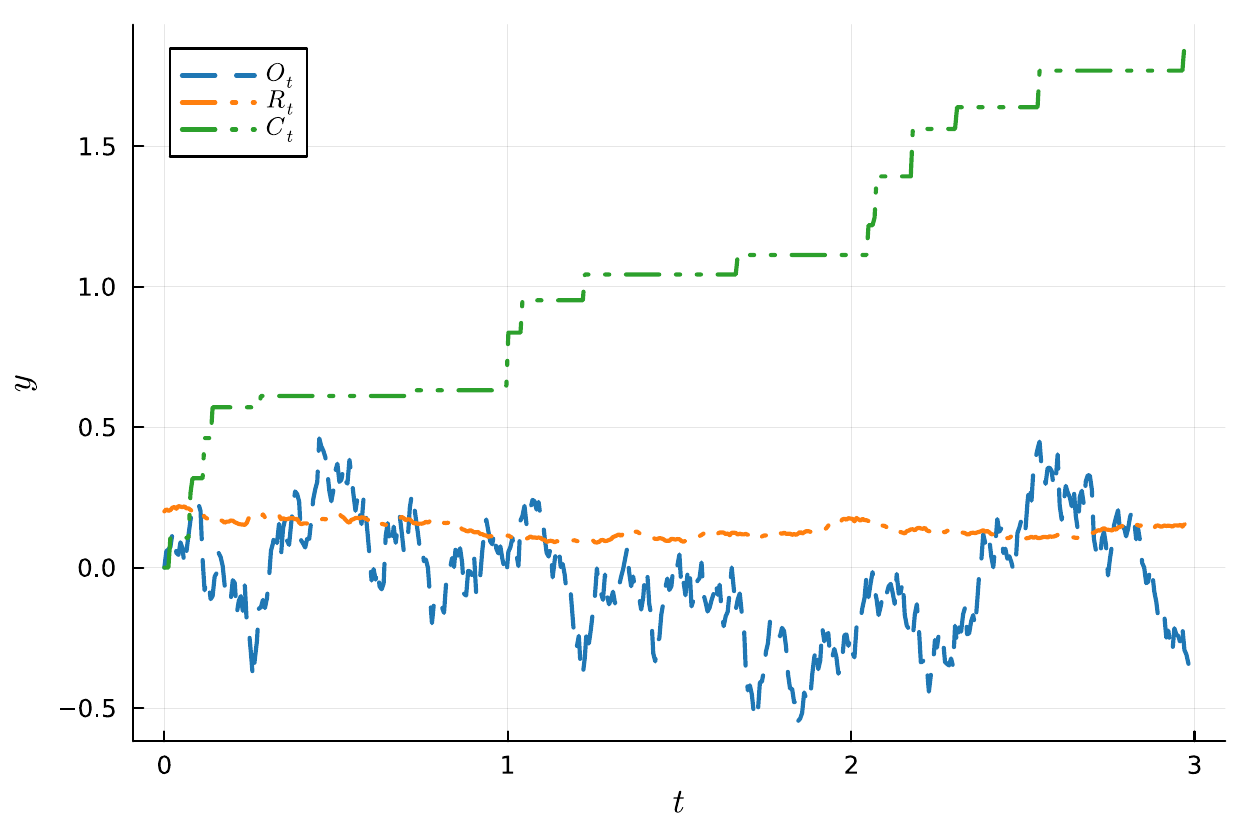}}
    \caption{Sample noises for the risk model \cref{riskmodel}.}
    \label{figriskmodelnoise}
\end{figure}

\begin{figure}[htb]
    \centerline{\includegraphics[width=0.8\textwidth]{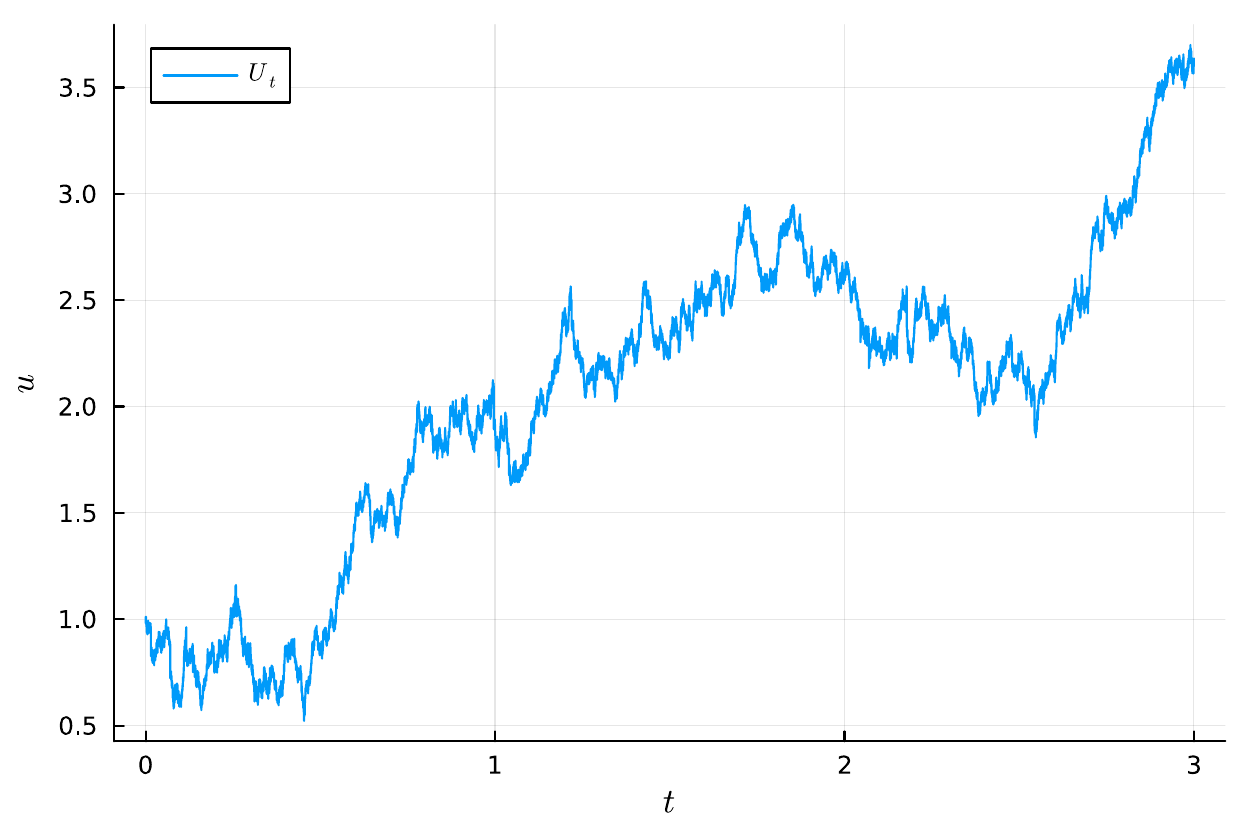}}
    \caption{Sample surplus solution for the risk model \cref{riskmodel}.}
    \label{figriskmodelsurplus}
\end{figure}

Although we have an explicit solution, it is not in closed form, so the target approximation for computing the error is obtained by solving the equation via Euler method in the much finer mesh. \Cref{tableriskmodel} shows the estimated strong error obtained with the Monte-Carlo method, while \cref{figriskmodel} illustrates the order of convergence. \cref{figriskmodelnoise} shows a sample path of the noise, which is composed of three processes, while \cref{figriskmodelsurplus} shows a sample path of the surplus.

\begin{table}
    \begin{center}
        \begin{tabular}[htb]{|r|l|l|l|}
            \hline N & dt & error & std err \\
            \hline \hline
            64 & 0.0469 & 0.224 & 0.0164 \\
            128 & 0.0234 & 0.112 & 0.00847 \\
            256 & 0.0117 & 0.0559 & 0.00397 \\
            512 & 0.00586 & 0.0295 & 0.00199 \\
            \hline
        \end{tabular}
    \end{center}

    \bigskip

    \caption{Mesh points (N), time steps (dt), strong error (error), and standard error (std err) of the Euler method for a risk model for each mesh resolution $N$, with initial condition $X_0 = 1.0$ and coupled Ornstein-Uhlenbeck, geometric Brownian motion, and compound Poisson processes, on the time interval $I = [0.0, 3.0]$, based on $M = 400$ sample paths for each fixed time step, with the target solution calculated with $262144$ points. The order of strong convergence is estimated to be $p = 0.976$, with the 95\% confidence interval $[0.8938, 1.0584]$.}

    \label{tableriskmodel}
\end{table}

\subsubsection{A random Fisher-KPP nonlinear PDE driven by boundary noise}

Finally, we simulate a Fisher-KPP equation with random boundary conditions, as inspired by the works of \cite{SalakoShen2020}  and \cite{FreidlinWentzell1992}. The first work addresses the Fisher-KPP equation with a random reaction coefficient, while the second work considers more general reaction-diffusion equations but driven by random boundary conditions.

The intent here is to illustrate the strong order 1 convergence rate on a discretization of a nonlinear partial differential equation. We use the method of lines (MOL), with finite differences in space, to approximate the random partial differential equation (PDE) by a system of random ODEs.

The deterministic Fisher-KPP equation is a nonlinear parabolic equation of react\-ion-diffusion type, with origins in \cite{Fisher1937} and \cite{KPP1937}. It models inhomogeneous population growth and many other phenomena displaying wave propagation, such as combustion front wave propagation, physiollogy, crystallography pattern formation, and so on.

We consider the Fisher-KPP equation driven by Neumann boundary conditions, with a random influx on the left end point and no flux on the right end point. The equation takes the form
\begin{equation}
    \label{fisherkpprode}
    \frac{\partial u}{\displaystyle \partial t} = \mu\frac{\partial^2 u}{\partial x^2} + \lambda u\left(1 - \frac{u}{u_m}\right), \quad (t, x) \in (0, \infty) \times (0, 1),
\end{equation}
endowed with the boundary conditions
\begin{equation}
    \label{fisherkppbc}
    \frac{\partial u}{\partial x}(t, 0) = - Y_t, \quad \frac{\partial u}{\partial x}(t, 1) = 0,
\end{equation}
and a given initial condition
\[
   u(0, x) = u_0(x).
\]

The unknown $u(t, x)$ represents the population density at time $t$ and point $x$; $D$ is a diffusivity coefficient; $\lambda$ is a reaction, or proliferation, coefficient; and $u_m$ is a carrying capacity density coefficient.

The random process $\{Y_t\}_t$, which drives the flux on the left boundary point, is taken to be a colored noise modulated by an exponentially decaying Hawkes process, representing random wave trains of incoming populations. More precisely, $Y_t = H_t O_t$, where $\{H_t\}_t$ is a Hawkes process with initial rate $\lambda_0 = 5.0$, base rate $a = 1.4$, exponential decay rate $\delta = 8.0$, and jump law following an exponential distribution with scale $\theta = 1/0.4;$ and $\{O_t\}_t$ is an Ornstein-Uhlenbeck process with $O_0 = 0.0,$ time-scale $\tau = 0.005$, drift $\nu = 1/\tau$ and diffusion $\sigma = \zeta/\tau,$ where $\zeta = \sqrt{\tau}.$ Since both $\{H_t\}_t$ and $\{O_t\}_t$ are semi-martingales and the class of semi-martingales is an algebra \cite[Corollary II.3]{Protter2005}, the modulated process $\{Y_t\}_t$ is also a semi-martingale.

This equation displays traveling wave solutions with a minimum wave speed of $2 \sqrt{\lambda \mu}$. We choose $\lambda = 10$ and $\mu= 0.009$, so the limit traveling speed is about $0.6$. The carrying capacity is set to $u_m = 1.0$.

The initial condition is taken to be zero, $u_0(x) = 0$, so all the population originates from the left boundary influx.

The mass within the region $0\leq x \leq 1$ satisfies
\[
   \frac{\mathrm{d}}{\mathrm{d} t} \int_0^1 u(t, x) \;\mathrm{d}x = \mu\int_0^1 u_{xx}(t, x) \;\mathrm{d}x + \lambda \int_0^1 u(t, x)\left(1 - \frac{u(t, x)}{u_m}\right)\;\mathrm{d}x.
\]
Using the boundary conditions, we find that
\[
   \frac{\mathrm{d}}{\mathrm{d}t} \int_0^1 u(t, x) \;\mathrm{d}x = \mu Y_t  + \frac{\lambda}{u_m} \int_0^1 u(t, x)\left(u_m - u(t, x)\right)\;\mathrm{d}x,
\]
which is nonnegative, provided $0 \leq u \leq u_m$ and $Y_t \geq 0$.

The equation involves a nonlinear term which is not globally Lipschitz continuous, but, similiarly to the population dynamics model considered in \cref{secpopdyn}, the region $0 \leq u(t, x) \leq u_m$ is invariant, so that the nonlinear term can be modified outside this region in order to satisfy the required uniform global estimates without affecting the dynamics within this region. The initial condition is chosen to be within this region almost surely. The procedure is the same as that done in \cref{secpopdyn}, so the details are omited.

For the time-mesh parameters, we set $N_{\textrm{tgt}} = 2^{18}$ and $N_i = 2^5, 2^7, 2^9.$ The spatial discretization is done with finite differences, with the number of spatial points depending on the time mesh, for stability and convergence reasons. Indeed, the Von Neumann stability analysis requires that $2\mu\Delta t / \Delta_x^2 \leq 1.$ With that in mind, for each $N_i = 2^5, 2^7, 2^9$, we take the $K_i + 1$ spatial points $0 = x_0 < \ldots x_{K_i},$ with $K_i = 2^3,$ $2^4,$ and $2^5,$ respectively, while for the target solution, we use $K_{\textrm{tgt}} = 2^9.$

\begin{figure}[htb]
    \centerline{\includegraphics[width=0.8\textwidth]{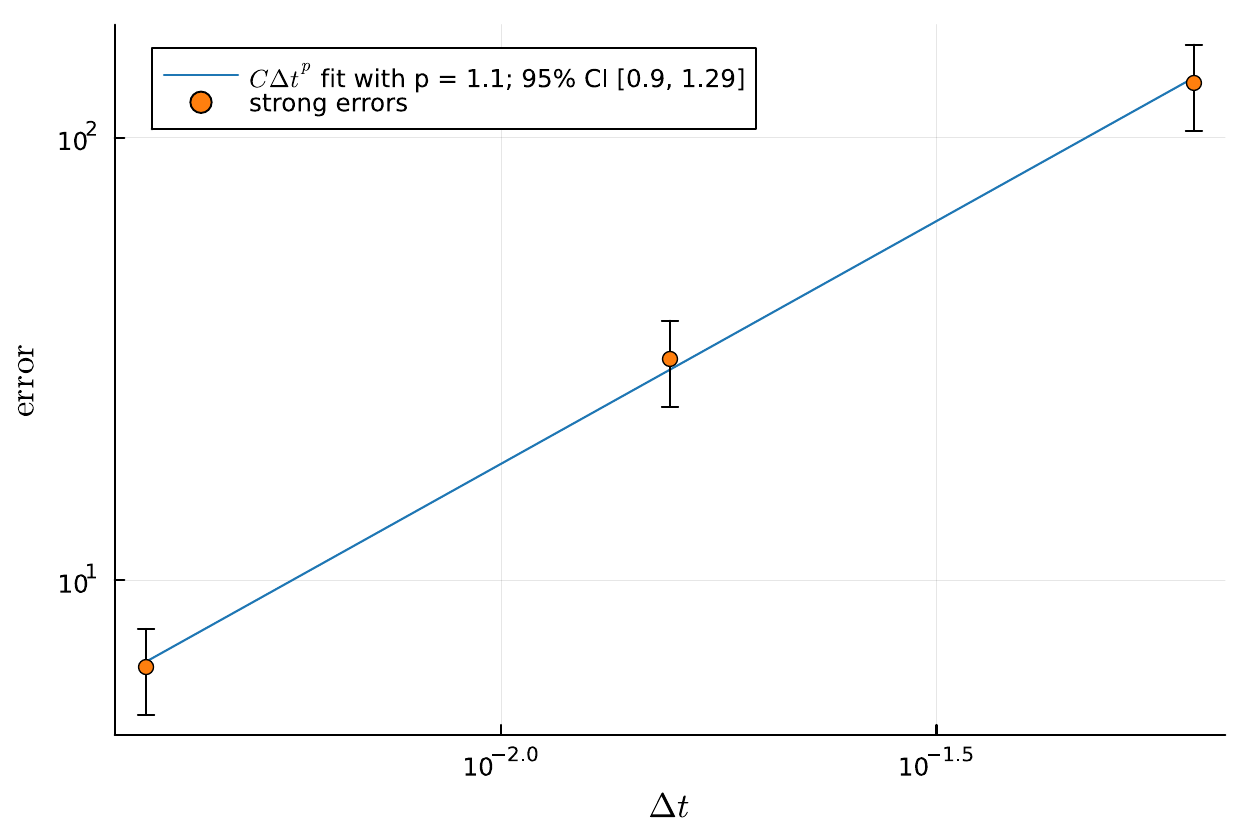}}
    \caption{Order of convergence of the strong error of the Euler method for the Fisher-KPP model \cref{fisherkpprode}-\cref{fisherkppbc}, based on \cref{tablefisherkpp}.}
    \label{figfisherkpp}
\end{figure}

\Cref{tablefisherkpp} shows the estimated strong error obtained with the Monte-Carlo method, while \cref{figfisherkpp} illustrates the order of convergence.

\begin{table}
    \begin{center}
        \begin{tabular}[htb]{|r|l|l|l|}
            \hline N & dt & error & std err \\
            \hline \hline
            32 & 0.0625 & 132.0 & 21.2 \\
            128 & 0.0156 & 31.5 & 5.03 \\
            512 & 0.00391 & 6.01 & 0.959 \\
            \hline
        \end{tabular}
    \end{center}

    \bigskip

    \caption{Mesh points (N), time steps (dt), strong error (error), and standard error (std err) of the Euler method for the Fisher-KPP equation for each mesh resolution $N$, with initial condition $X_0 = 0$ and Hawkes-modulated Ornstein-Uhlenbeck colored noise, on the time interval $I = [0.0, 2.0]$, based on $M = 40$ sample paths for each fixed time step, with the target solution calculated with $262144$ points. The order of strong convergence is estimated to be $p = 1.116$, with the 95\% confidence interval $[0.9995, 1.2318]$.}

    \label{tablefisherkpp}
\end{table}

\bibliographystyle{plain}
\bibliography{rode_convergence_euler}

\begin{thebibliography}{10}

\bibitem{AitSahaliaJacod2014}
Y.~A\"it-Sahalia and J.~Jacod.
\newblock {Chapter 1. From Diffusions to Semimartingales}.
\newblock In {\em High-Frequency Financial Econometrics}. Princeton University
  Press, 2014.

\bibitem{AB2006}
C.~D. Aliprantis and K.~C. Border.
\newblock {\em Infinite Dimensional Analysis: a Hitchhiker's Guide}.
\newblock Springer, Berlin, Heidelberg, 2006.

\bibitem{Asai2016}
Y.~Asai.
\newblock {\em {Numerical Methods for Random Ordinary Differential Equations
  and their Applications in Biology and Medicine}}.
\newblock doctoral thesis, Johann Wolfgang Goethe Universit\"at, Frankfurt am
  Main, 2016.

\bibitem{AsaiKloeden2016}
Y.~Asai and P.~E. Kloeden.
\newblock Numerical schemes for random odes with affine noise.
\newblock {\em Numer. Algorithms}, 72(1):155--171, 2016.

\bibitem{BENDER200381}
C.~Bender.
\newblock An {I}t\^o formula for generalized functionals of a fractional
  {B}rownian motion with arbitrary {H}urst parameter.
\newblock {\em Stochastic Process. Appl.}, 104(1):81--106, 2003.

\bibitem{JSSv098i16}
M.~Besan\c{c}on, T.~Papamarkou, D.~Anthoff, A.~Arslan, S.~Byrne, D.~Lin, and
  J.~Pearson.
\newblock {Distributions.jl: Definition and modeling of probability
  distributions in the JuliaStats ecosystem}.
\newblock {\em Journal of Statistical Software}, 98(16):1--30, 2021.

\bibitem{Julia2017}
J.~Bezanson, A.~Edelman, S.~Karpinski, and V.~B. Shah.
\newblock {Julia: A Fresh Approach to Numerical Computing}.
\newblock {\em SIAM Rev.}, 59(1):65--98, 2017.

\bibitem{BHOB2008}
F.~Biagini, Y.~Hu, B.~{\O}ksendal, and T.~Zhang.
\newblock {\em Stochastic Calculus for Fractional Brownian Motion and
  Applications}.
\newblock Probability and Its Applications. Springer-Verlag, London, 2008.

\bibitem{BlackmanVigna2021}
D.~Blackman and S.~Vigna.
\newblock Scrambled linear pseudorandom number generators.
\newblock {\em ACM Trans. Math. Software}, 47(4), September 2021.

\bibitem{BogdanoffGoldbergBernard1961}
J.~L. Bogdanoff, J.~E. Goldberg, and M.~C. Bernard.
\newblock {Response of a simple structure to a random earthquake-type
  disturbance}.
\newblock {\em Bulletin of the Seismological Society of America},
  51(2):293--310, 1961.

\bibitem{BrigoMercurio2006}
D.~Brigo and F.~Mercurio.
\newblock {\em Interest-Rate Models: Theory and Practice}.
\newblock Springer Finance. Springer-Verlag, Berlin, Heidelberg, 2001.

\bibitem{PlotsJL2023}
S.~Christ, D.~Schwabeneder, C.~Rackauckas, M.~K. Borregaard, and T.~Breloff.
\newblock {Plots.jl -- a user extendable plotting API for the Julia programming
  language}.
\newblock {\em Journal of Open Research Software}, 2023.

\bibitem{CoddingtonLevinson1985}
A.~Coddington and N.~Levinson.
\newblock {\em Theory of Ordinary Differential Equations}.
\newblock International Series in Pure and Applied Mathematics. R.E. Krieger,
  1984.

\bibitem{DeGrootSchervish2018}
M.~H. DeGroot and M.~J. Schervish.
\newblock {\em Probability and Statistics}.
\newblock Pearson Education, 4th edition, 2018.

\bibitem{DiekerMandjes2003}
A.B. Dieker and M.~Mandjes.
\newblock {On spectral simulation of Fractional Brownian motion}.
\newblock {\em Probab. Engrg. Inform. Sci.}, 17(3):417--434, 2003.

\bibitem{Fisher1937}
R.~A. Fisher.
\newblock {The wave of advance of advantageous genes}.
\newblock {\em Annals Eugen.}, 7(4):355--369, 1937.

\bibitem{FreidlinWentzell1992}
M.~I. Freidlin and A.~D. Wentzell.
\newblock {Reaction-diffusion equations with randomly perturbed boundary
  conditions}.
\newblock {\em Ann. Probab.}, 20(2):963--986, 1992.

\bibitem{GerberShiu1998}
H.~U. Gerber and E.~S.~W. Shiu.
\newblock On the time value of ruin.
\newblock {\em N. Am. Actuar. J.}, 2(1):48--72, 1998.

\bibitem{GruneKloeden2001}
L.~Gr{\"u}ne and P.~E. Kloeden.
\newblock Higher order numerical schemes for affinely controlled nonlinear
  systems.
\newblock {\em Numer. Math.}, 89(4):669--690, 2001.

\bibitem{HanKloeden2017}
X.~Han and P.~E. Kloeden.
\newblock {\em Random Ordinary Differential Equations and Their Numerical
  Solution}.
\newblock Springer Publishing Company, Incorporated, 1st edition, 2017.

\bibitem{HastieTibshiraniFriedman2009}
T.~Hastie, R.~Tibshirani, and J.~Friedman.
\newblock {\em The Elements of Statistical Learning}.
\newblock {Springer Series in Statistics}. Springer, New York, NY, 2nd edition,
  2009.

\bibitem{HeWangYan1992}
S.-W. He, J.-G. Wang, and J.-A. Yan.
\newblock {\em Semimartingale Theory and Stochastic Calculus}.
\newblock Taylor \& Francis Group, New York, 1992.

\bibitem{HighamKloeden2021}
D.~Higham and P.~E. Kloeden.
\newblock {\em An Introduction to the Numerical Simulation of Stochastic
  Differential Equations}.
\newblock Society for Industrial and Applied Mathematics, Philadelphia, PA,
  2021.

\bibitem{HousnerJenning1964}
G.~W. Housner and P.~C. Jennings.
\newblock Generation of artificial earthquakes.
\newblock {\em Journal of the Engineering Mechanics Division}, 90(1):113--150,
  1964.

\bibitem{JamesWittenHastieTibshirani2021}
G.~James, D.~Witten, T.~Hastie, and R.~Tibshirani.
\newblock {\em An Introduction to Statistical Learning}.
\newblock {Springer Texts in Statistics}. Springer, New York, NY, 2nd edition,
  2021.

\bibitem{JentzenKloeden2011}
A.~Jentzen and P.~E. Kloeden.
\newblock {\em Taylor Approximations for Stochastic Partial Differential
  Equations}.
\newblock Society for Industrial and Applied Mathematics, 2011.

\bibitem{JohnsonWichern2007}
R.~A. Johnson and D.~W. Wichern.
\newblock {\em Applied Multivariate Statistical Analysis}.
\newblock Pearson Prentice-Hall, 6th edition, 2007.

\bibitem{Kanai1957}
K.~Kanai.
\newblock Semi-empirical formula for the seismic characteristics of the ground.
\newblock {\em Bull. Earthq. Res. Inst.}, 35:309--325, 1957.

\bibitem{KaratzasShreve2014}
I.~Karatzas and S.~E. Shreve.
\newblock {\em {Brownian Motion and Stochastic Calculus}}.
\newblock {Graduate Texts in Mathematics}. Springer, New York, NY, 2nd edition,
  2014.

\bibitem{RODEConvEM2023}
P.~Kloeden and R.~Rosa.
\newblock Numerical examples of strong order of convergence of the euler method
  for random ordinary differential equations.
\newblock \url{https://github.com/rmsrosa/rode_convergence_euler} and
  \url{https://rmsrosa.github.io/rode_convergence_euler}, 2025.

\bibitem{KloedenJentzen2007}
P.~E. Kloeden and A.~Jentzen.
\newblock Pathwise convergent higher order numerical schemes for random
  ordinary differential equations.
\newblock {\em Proc. Roy. Soc. London Sect. A.}, 463:2229--2944, 2007.

\bibitem{KloedenPlatenSchurz2012}
P.~E. Kloeden, E.~Platen, and H.~Schurz.
\newblock {\em Numerical Solution of SDE Through Computer Experiments}.
\newblock {Universitext}. Springer-Verlag, Berlin Heidelberg, 2012.

\bibitem{KPP1937}
A.~Kolmogorov, I.~Petrovskii, and N.~Piscunov.
\newblock A study of the equation of diffusion with increase in the quantity of
  matter, and its application to a biological problem.
\newblock {\em Byul. Moskovskogo Gos. Univ.}, 1(6):1--25, 1937.

\bibitem{MandelbrotVanNess1968}
B.~B. Mandelbrot and J.~W.~Van Ness.
\newblock {Fractional Brownian motions, fractional noises and applications}.
\newblock {\em SIAM Rev.}, 10(4):422--437, 1968.

\bibitem{MARINELLI20161}
C.~Marinelli and M.~R\"ockner.
\newblock {On the maximal inequalities of Burkholder, Davis and Gundy}.
\newblock {\em Expositiones Mathematicae}, 34(1):1--26, 2016.

\bibitem{MerdanBekiryaziciKesemenKhaniyev2017}
M.~Merdan, Z.~Bekiryazici, T.~Kesemen, and T.~Khaniyev.
\newblock Comparison of stochastic and random models for bacterial resistance.
\newblock {\em Adv. Difference Equ.}, 2017(1):133, 2017.

\bibitem{Metivier1982}
M.~M\'etivier.
\newblock {\em Semimartingales}.
\newblock De Gruyter, Berlin, New York, 1982.

\bibitem{Mishura2008}
Y.~S. Mishura.
\newblock {\em Stochastic calculus for fractional Brownian motion and related
  processes}.
\newblock Lecture Notes in Mathematics 1929. Springer-Verlag, Berlin,
  Heidelberg, 2011.

\bibitem{NeckelRupp2013}
T.~Neckel and F.~Rupp.
\newblock {\em Random Differential Equations in Scientific Computing}.
\newblock De Gruyter, London, 2013.

\bibitem{Oksendal2003}
B.~{\O}ksendal.
\newblock {\em Stochastic Differential Equations : An Introduction with
  Applications}.
\newblock Springer-Verlag, Berlin, Heidelberg, 2003.

\bibitem{Protter2005}
P.~E. Protter.
\newblock {\em Stochastic Integration and Differential Equations}.
\newblock Springer-Verlag, Berlin, Heidelberg, 2nd edition, 2005.

\bibitem{DifferentialEquations.jl-2017}
C.~Rackauckas and Q.~Nie.
\newblock {DifferentialEquations.jl -- A performant and feature-rich ecosystem
  for solving differential equations in Julia}.
\newblock {\em The Journal of Open Research Software}, 5(1), 2017.

\bibitem{SalakoShen2020}
R.~B. Salako and W.~Shen.
\newblock {Long time behavior of random and nonautonomous Fisher--KPP
  equations: Part I---stability of equilibria and spreading speeds}.
\newblock {\em J. of Dynam. Differential Equations}, 33(2):1035--1070, 2021.

\bibitem{StrasserTheisMarr2012}
M.~Strasser, F.~J. Theis, and C.~Marr.
\newblock Stability and multiattractor dynamics of a toggle switch based on a
  two-stage model of stochastic gene expression.
\newblock {\em Biophysical Journal}, 102(1):19--29, 2012.

\bibitem{Sugiyama1969}
S.~Sugiyama.
\newblock Stability problems on difference and functional-differential
  equations.
\newblock {\em Proceedings of the Japan Academy}, 45(7):526--529, 1969.

\bibitem{Tajimi1960}
H~Tajimi.
\newblock A statistical method of determining the maximum response of a
  building structure during an earthquake.
\newblock {\em Proceedings of the 2nd World Conference on Earthquake
  Engineering, Tokyo, Japan, 1960}, pages 781--797, 1960.

\bibitem{VerdCrombachJaeger2014}
B.~Verd, A.~Crombach, and J.~Jaeger.
\newblock Classification of transient behaviours in a time-dependent toggle
  switch model.
\newblock {\em BMC Systems Biology}, 8(1):43, 2014.

\bibitem{WangCaoHanKloeden2021}
P.~Wang, Y.~Cao, X.~Han, and P.~E. Kloeden.
\newblock Mean-square convergence of numerical methods for random ordinary
  differential equations.
\newblock {\em Numer. Algorithms}, 87(1):299--333, 2021.

\end{thebibliography}

\end{document}